\pdfoutput=1
\documentclass[11pt]{article}
\usepackage{amssymb}

\usepackage{tikz}
\usetikzlibrary{arrows.meta}
\begin{document}
\topmargin 0.0in
\textheight 8.15in
\hsize 6.7in
\oddsidemargin 0 in
\title{ A scaling limit  of crossing probabilities for critical percolation on the square  lattice
} 
\author{Yu Zhang}
\maketitle
\begin{center}
{\bf Abstract}
\end{center}
We show the existence of a scaling limit for the crossing probabilities on the square lattice  in an equilateral triangle for the critical percolation.
We also show that Cardy's formula does not hold on the square lattice  for the critical  percolation.

\section{ Introduction and statement of results.}
\subsection{Introduction of the percolation model.}
Percolation is one of the fundamental stochastic models studied by mathematicians. Since percolation is
one of the simplest models that exhibit a phase transition, it also becomes  
one of statistical physicists' favorite models for studying critical phenomena. 
In particular, the case of two dimensions is 
special since it links the conformal field theory with the singularity of a phase
transition at criticality.  
We first introduce  the square lattice.  The $\delta$-square lattice is defined by
$${\bf Z}^2_\delta=\{(\delta i, \delta j): \mbox{ $i$ and $j$ are integers and $\delta>0$}\},$$  
and  we denote {\em edges} connecting a
pair of vertices
$u=(u_1, u_2)$ and $v=(v_1, v_2)$ 
with $d(u,v)=\delta$, where $d(u,v)$ is the Euclidean distance between $u$ and $v$. In other words, ${\bf Z}^2_\delta$ is a square lattice scaled to mesh-size $\delta$.
For any two sets ${A}$ and ${B}$, we also define their distance by
$$d({A}, {B})= \min\{ d(u, v): u\in {A}, v\in {B}\}.$$
When $\delta=1$, ${\bf Z}^2_1$ is the  standard square lattice
denoted by ${\bf Z}^2$ and originally considered  by Broadbent and Hammersley (1957).

 Another popular lattice is the triangular lattice with mesh $\delta$. Divide ${\bf R}^2$ into equilateral triangles by means of the horizontal lines
and lines under an angle of $\pi/3$ or $2\pi/3$ with the first coordinate-axis through the points $(\delta k, 0), k\in {\bf Z}$. The vertices of the graph are the vertices of the equilateral   triangles, and the edges are the line segments connecting two vertices
of the same triangle. We divide these line segments by {\em horizontal, north-east, north-west} edges. 
We denote the above graph by ${\bf T}_\delta$.   We also consider the lattice with the same vertices of ${\bf T}_\delta$, keeping all the 
north-east  and the north-west edges,  but removing all the horizontal  edges from ${\bf T}_\delta$. We denote this new lattice by ${\bf N}_\delta$ (see Fig. 1).
$ {\bf N}_\delta $ consists of  the regular parallelograms with side length $\delta$, called $\delta$-{\em parallelograms}.
 The vertices of  ${\bf N}_\delta$ are the vertices of the regular parallelograms, and the edges, called {\em diagonal } edges, are the edges of regular parallelograms. In addition, we also consider the dual graph of  ${\bf N}_\delta$  with the same vertices in ${\bf N}_\delta$ and with the extra
 the horizontal and the  vertical edges besides of these diagonal edges in each $\delta$-parallelogram of ${\bf N}_\delta$. 
 The dual graph is called ${\bf N}_\delta^*$ and
 $({\bf N}_\delta, {\bf N}_\delta^*)$ is called the {\em match pair} (see section 2.2 in Kesten (1982)).
 Two vertices  ${u}=(u_1, u_2)$ 
and ${v}=(v_1, v_2)$ are said to be $\bullet$-{\em adjacent}   if  there is an edge (a diagonal edge) of ${\bf N}_\delta$ with vertices of $u$ and $v$. 
In addition, two vertices  ${u}=(u_1, u_2)$ 
and ${v}=(v_1, v_2)$ are said to be $\circ$-{\em adjacent}  if there is an edge (a diagonal  or a  horizontal or a vertical edge) of  
${\bf N}_\delta^*$ with vertices of $u$ and $v$.  If two vertices are $\bullet$-adjacent , then they are also $\circ$-adjacent.



We introduce site  percolation model on ${\bf N}_\delta$.
For each vertex  in ${\bf N}_\delta$, it is either {\em open} or {\em closed}
 independently with  probability $p$ and $1-p$. 
 The corresponding probability measure on the configurations of open 
and closed  vertices  is denoted by ${\bf P}_{p, \delta}$.  This model is called the  {\em site  percolation} model.
If each edge is either   open or closed, then it is  bond percolation.
In this paper, we mainly work on the site percolation on ${\bf N}_\delta$. 
A $\bullet$-path  from $u$ to $v$ is a sequence $(v_0...v_{i} v_{i+1}...v_n)$
with distinct vertices $v_i$, $0\leq i\leq n$,
 and $v_0=u$, and $v_n=v$ in ${\bf N}_\delta$ such that $v_i$ and $v_{i+1}$ are $\bullet$-adjacent. 
 Similarly,  a $\circ$-path  from $u$ to $v$ is a sequence $(v_0...v_{i} v_{i+1}...v_n)$
with distinct vertices $v_i$, $0\leq i\leq n$,
 and $v_0=u$, and $v_n=v$ in ${\bf N}_\delta$ such that $v_i$ and $v_{i+1}$ are $\circ$-adjacent. 
 Note that a $\bullet$-path is also a $\circ$-path.
If all  the vertices in a path are open, then the path is called an {\em open path}.
Similarly, if  all  the vertices in a path are closed , then the path is called a {\em closed path}. 
The cluster of the vertex $x$, 
${\bf C}(x)$,
consists of all vertices that are connected to $x$ by an open $\bullet$-path.
 A circuit from $v_1$ to $v_n$ is a path  with distinct vertices $v_i$ ($1\leq i\leq n-1$) and $v_0=v_n$. 
For any collection ${A}$ of vertices, $| A |$ 
denotes the cardinality of $A$. We choose ${\bf 0}$ as the origin. The percolation probability is
\begin{eqnarray*}
\theta (p)= {\bf P}_{p}(|{\bf C}({\bf 0})|=\infty),
\end{eqnarray*}
and the {\em critical percolation probability} is 
$$p_c=\sup\{p:\theta (p)=0\}.$$

We want to point out that
for both bond percolation on the square lattice and site percolation on the triangular
lattice, the critical probabilities can be precisely computed (see Kesten (1982)) as
$p_c=1/2$. People take advantage of  symmetry when $p_c=1/2$. In particular for the triangular lattice, its dual lattice is also the triangular lattice.
If we stretch each $\delta$-parallelogram horizontally  in ${\bf N}_\delta$ to be a rotated square, the lattice will be a rotated square lattice ${\bf Z}_\delta^2$.  Thus,  the probability of a stretched open $\bullet$-path 
will remain the same. Thus,  the percolation models of ${\bf N}_\delta$  and ${\bf Z}^2_\delta$  are essentially the same. In particular,  the critical percolation probabilities of ${\bf N}_\delta$  and ${\bf Z}^2_\delta$ are the same for  both site or bond percolation models.

\subsection{ Scaling limits for crossing probability.}
Now we will consider the complex plane ${\bf C}$ and  the lattice ${\bf N}_{\delta}$.
Let ${D} \subset {\bf C}$ be a simple connected and  bounded subset of ${\bf C}$, whose boundary $\partial {D}$ is a Jordan curve. ${D}$ is called 
a {\em Jordan domain}. We select four points $A, X, B, C$ from $\partial D$ with $X\in \overline{AB}$ such that 
$$\partial D=\overline{AX}\cup\overline{XB}\cup\overline{BC}\cup\overline{CA}.$$
With a path from $\overline{AX}$ to $\overline{BC}$, it means that  the  path
$v_0v_1\cdots v_n$ in ${\bf N}_\delta$ is such that $v_1\cdots v_{n-1}$ lie inside  without using the boundary of $D$. Moreover, $v_0$ is  on 
$\overline{AX}$  or $v_0$  is outside $D$ with the edge from $v_0$ to $v_1$ meets $\overline{AX}$, and $v_n$ is  on 
$\overline{BC}$  or $v_n$  is outside $D$ with the edge from $v_{n-1}$ to $v_n$ meets $\overline{BC}$.
With an open $\bullet$-path or a closed $\circ$-path
from $\overline{AX}$ to $\overline{BC}$, it means that all the vertices on the path  are open or closed. Similarly, we define  an  open $\bullet$-path or a  closed $\circ$-path from $\overline{XB}$ to $\overline{AC}$. 
We denote  the {\em crossing probability} by 
$$\pi_{p, \delta} (X)=\pi_{p,\delta}({D}, \overline{AX}, \overline{BC})= {\bf P}_{ p,\delta}\left(\exists \mbox{ an open $\bullet$-path in $D$ from $\overline{AX}$ to $\overline{BC}$}\right).\eqno{}$$
For simplicity,  let
$$\pi_{\delta}({D}, \overline{AX}, \overline{BC})=\pi_{p_c,\delta}({D}, \overline{AX},\overline{BC}).$$

With the definitions, it is known (see Kesten (1982)) that for all ${D}$ and $A, X, B, C$ with $X\neq A, B$,
$$\lim_{\delta\rightarrow 0} \pi_{p,\delta}({D}, \overline{AX}, \overline{BC})=\left \{ \begin{array}{ll}0 & \mbox{ if $p < p_c,$}\\
1 &\mbox{ if $p>p_c$.}
\end{array}\right. \eqno{(1.2.1)}
$$
When $p=p_c$, it is also known   (see Kesten (1982)) that
$$0<\liminf_{\delta \rightarrow 0} \pi_{\delta}({D}, \overline{AX},\overline{BC})\leq \limsup_{\delta \rightarrow 0} 
 \pi_{\delta }({D}, \overline{AX},\overline{BC})<1.\eqno{(1.2.2)}$$
With (1.2.2), the existence of the {\em scaling limit}  was  conjectured by  Langland, Pouliot and Saint-Aubin (1994):
$$\lim_{\delta \rightarrow 0} \pi_{\delta}({D}, \overline{AX},\overline{BC})
=\pi({D}, \overline{AX},\overline{BC}).\eqno{(1.2.3)}$$
  By the Riemann mapping theorem, there is a unique conformal mapping $\phi$ that maps ${D}$ to the unit equilateral triangle $\bigtriangleup abc$ on ${\bf C}$ with 
$$\phi(A)=a=0, \phi(B)=b=1,   \phi(C)=c=e^{i\pi/3}, \phi(X)=x,\eqno{(1.2.4)}$$
for $X\in \overline{AB}$.   
With this map, it was  conjectured by Cardy (1992)
that the limit probability of (1.2.3)  is equal to $x$.
This conjecture is called Cardy's formula with
Carleson version.  \\

{\bf Cardy's formula.} {\em If $D$ is a Jordan domain with four points $A, X, B, C$  for any $X\in \overline{AB}$ defined  above, and $\phi$ is the Riemann mapping such that $ \phi(A)=a=0, \phi(B)=b=1,   \phi(C)=c=e^{i\pi/3}, \phi(X)=x$, then }
$$\lim_{\delta \rightarrow 0} \pi_\delta (D, \overline {AX}, \overline {BC})=x.$$

Smirnov's  celebrated  work (2001) showed that the scaling limit in (1.2.3) indeed holds
and satisfies  Cardy's formula on the triangular lattice. 
Using Cardy's formula has
been quite successful at calculating
 the exponents regarding the behaviors of percolation functions.

We focus on an equilateral triangle  $A=0, B=1,$ and $ C=e^{i\pi/3}$  and  consider the crossing probability on this  equilateral triangle $D=\bigtriangleup ABC$.
From now on, we only consider  this  equilateral triangle $\bigtriangleup ABC$.
 We select a point $X\in \overline{AB}$.  We call D  a four point-domain with boundary points $A, X, B, C$.
For convenience,  we consider $\delta=l^{-1}$ for integer $l$ and  the lattice ${\bf N}_{l^{-1}}$, but the results in this paper can be proved by the same arguments for any mesh-size $\delta$.  
 In this paper, we show that the scaling limit exists in the  above special domain.\\

\begin{figure}
\begin{center}
\setlength{\unitlength}{0.0125in}%
\begin{picture}(50,220)(67,680)
\thicklines

\put(170,680){\line(2,3){100}}
\put(210,680){\line(2,3){80}}
\put(250,680){\line(2,3){60}}
\put(290,680){\line(2,3){40}}
\put(330,680){\line(2,3){20}}

\put(370,680){\line(-2,3){100}}
\put(90,800){\line(0,1){150}}
\put(90,950){\line(1,0){150}}
\put(90,800){\line(1,1){150}}
\put(120,830){\line(0,1){120}}
\put(90,920){\line(1,0){120}}
\put(150,860){\line(0,1){90}}
\put(180,890){\line(0,1){60}}
\put(210,920){\line(0,1){30}}
\put(90,890){\line(1,0){90}}
\put(90,860){\line(1,0){60}}
\put(90,830){\line(1,0){30}}
\put(90,950){\line(1,0){90}}

\put(150,680){\line(-1,1){110}}
\put(-70,680){\line(1,1){110}}
\put(105,680){\line(-1,1){87}}
\put(-25,680){\line(1,1){87}}
\put(63,680){\line(-1,1){66}}
\put(20,680){\line(1,1){65}}
\put(20,680){\line(-1,1){45}}
\put(63,680){\line(1,1){43}}
\put(-25,680){\line(-1,1){22}}
\put(105,680){\line(1,1){22}}

\put(170,680){\line(1,0){200}}
\put(-70,680){\line(1,0){220}}

\put(210,680){\line(-2,3){20}}
\put(250,680){\line(-2,3){40}}
\put(290,680){\line(-2,3){60}}
\put(330,680){\line(-2,3){80}}

\put(61,768){\circle*{8}}
\put(41,745){\circle*{8}}
\put(20,724){\circle*{8}}
\put(42,702){\circle*{8}}
\put(20,680){\circle*{8}}

\put(150,862){\circle*{8}}
\put(150,890){\circle*{8}}
\put(120,920){\circle*{8}}
\put(120,950){\circle*{8}}
\put(120,890){\circle*{8}}

\put(290,800){\circle*{8}}
\put(270,770){\circle*{8}}
\put(250,740){\circle*{8}}
\put(270,710){\circle*{8}}
\put(250,680){\circle*{8}}

\put(-65,660){\mbox{$\alpha=0$}{}}
\put(125,660){\mbox{$\beta=1$}{}}
\put(70,670){\mbox{$X$}{}}
\put(290,670){\mbox{$X$}{}}
\put(185,885){\mbox{$X$}{}}

\put(170,660){\mbox{$A=0$}{}}
\put(360,660){\mbox{$B=1$}{}}
\put(-10,800){\mbox{$\gamma=1/2+i/2$}{}}
\put(260,835){\mbox{$C=e^{i\pi/3}$}{}}
\put(90,790){\mbox{$a=0$}{}}
\put(250,940){\mbox{$b={1+i\over \sqrt{2}}$}{}}
\put(90,960){\mbox{$c=i/\sqrt{2}$}{}}
\end{picture}
\end{center}
\caption{ {\em  The lower-right graph is $\bigtriangleup  ABC$ with vertices  and edges on ${\bf N}_{\delta}$. The lower-left graph is $\bigtriangleup \alpha\beta\gamma$ with vertices by stretching  each parallelogram in ${\bf N}_\delta$ horizontally to be a square.  The upper graph is  $\bigtriangleup  abc$ by rotation $\bigtriangleup \alpha\beta\gamma$ in $\pi/4$ counterclockwise. Each open $\bullet$-path, solid circle path,  from $\overline{BC}$ to $\overline{A X}$ is corresponding to an open $\bullet$-path from $\overline{\beta\gamma}$ to $\overline{\alpha X}$, and from $\overline{bc}$ to $\overline{aX}$ with the same probability.}}
\end{figure}




{\bf Theorem.} {\em If $D= \bigtriangleup ABC$
 is the unit equilateral with four points $A, X, B, C$  for any $X\in \overline{AB}$ defined  above, then, for site percolation on ${\bf N}_\delta$, 
$$\lim_{l \rightarrow \infty} \pi_{l^{-1}}(D, \overline{AX}, \overline{BC}) =\lim_{l \rightarrow \infty} \pi_{l^{-1}}(X)=X.$$
So Cardy's formula holds for $D=\bigtriangleup ABC$.}\\

{\bf Remark 1.}  
The proof of  the Theorem  can be adapted directly for  bond percolation in ${\bf N}_\delta$ and  site percolation in the triangular lattice.
In fact,  site percolation in the triangular lattice is easier to hand, since we do not need to work on an extra dual lattice. Furthermore,  it needs to work on more complicated  topologies if one tries to show  the existence of the scaling limit for any Jordan  domain by using the method of the Theorem.  So we do not show it in this paper.\\

We can stretch  each $\delta$-parallelogram in ${\bf N}_\delta$ horizontally into a rotated square such that $D=\bigtriangleup ABC$ is changed to be $\bigtriangleup\alpha \beta\gamma$ for $\alpha= 0$, $\beta=1$ and $\gamma=1/2+i/2$ (see Fig. 1).  The corresponding crossing probability remains the  same  after stretching (see Fig. 1). 
We rotate $\bigtriangleup\alpha \beta\gamma$ at the origin in $\pi/4$ clockwise to have $\bigtriangleup abc$ (see Fig. 1) with $a=0, b=(1+i )/ \sqrt{2}$ and $c=i/\sqrt{2}$. 
Note that the crossing probabilities,  from $\overline{ BC } $ from $\overline{AX}$  and from $\overline{bc} $ from $\overline{ax}$ are the same (see Fig. 1) in $\bigtriangleup ABC$ and  in $\bigtriangleup abc$. 
Thus, we have the following corollary. \\

{\bf Corollary 1.} {\em If $\bigtriangleup abc$ is the triangle with four points $a, x, b, c$  for any $x\in \overline{ab}$ defined  above, then, for site  percolation on ${\bf Z}_{\delta/\sqrt{2}}^2$, 
$$\lim_{l \rightarrow \infty} \pi_{l^{-1}}(\bigtriangleup abc, \overline{ax}, \overline{bc}) =\lim_{l \rightarrow \infty} \pi_{l^{-1}}(x)=x$$}


Now we show Cardy's formula does not hold on ${\bf Z}^2_\delta$. Let $\phi$ be the unique  conformal map from $\bigtriangleup \alpha\beta\gamma$ to $\bigtriangleup ABC$ with $\phi(\alpha)=A, \phi(\beta)=B$ and $\phi(\gamma)=C$.  
We will show that  there is $X\in \overline{\alpha\beta}$ such that
$$\phi(X) \neq X.\eqno{(1.2.5)}$$
 To show (1.2.5),  
Schwarz-Christoffel's  transformation $\phi_1(w)=z $ maps  from the upper half-plane to $\Delta \alpha\beta\gamma$ such that $\phi_1(0)=\alpha$, $\phi_1(1)=\beta$, and $\phi_1(\infty)=\gamma$ with
$$z=\phi_1(w) =A_1\int^{w} _0 y^{-3/4}(1-y)^{-3/4}dy\mbox{ for }w\in(0, 1) \eqno{(1.2.6)}$$
for a constant $A_1>0$. In addition, 
Schwarz-Christoffel's  transformation $\phi_2(w)=z $ maps from the upper half-plane to $\Delta ABC$ such that
$\phi_2(0)=A$, $\phi_1(1)=B$ and $\phi_1(\infty)=C$ with
$$y=\phi_2(w) =A_2\int^{w} _0 y^{-2/3}(1-y)^{-2/3}dy\mbox{ for }w\in(0, 1)\eqno{(1.2.7)}$$
 for a constant $A_2>0$. By (1.2.5) and (1.2.6), 
$${dz\over dw}=A_1w^{-3/4}(1-w)^{-3/4}\mbox{ for }w\in (0, 1),\eqno{(1.2.8)}$$
and 
$${dy\over dw}=A_2 w^{-2/3}(1-w)^{-2/3}\mbox{ for }w\in (0, 1).\eqno{(1.2.9)}$$
Thus, $\phi_1\neq \phi_2$, otherwise, by (1.2.8) and (1.2.9),
 $$A_1w^{-1/12}(1-w)^{-1/12}=A_2 ,\eqno{(1.2.10)}$$
 but (1.2.10) cannot hold if $w\rightarrow 0$ or $w\rightarrow 1$.
If $\phi_1\neq \phi_2$, then there exists $X\in (0, 1)$ such that
$\phi_2^{-1}(X)\neq \phi_1^{-1}(X)$. In other words,
$$\phi_2\phi^{-1}_1(X)\neq X.\eqno{(1.2.11)}$$
Note that the conformal  map of $\phi$ is unique, so by (1.2.11),
$$\phi(X)=\phi_2\phi^{-1}_1(X) \neq  X.$$ 
Thus, (1.2.5) holds.


 With (1.2.5), we show that Cardy's formula does not hold on ${\bf Z}^2_\delta$.  In fact, we only need to show that Cardy's formula does not hold on 
 ${\bf Z}_{\delta /\sqrt{2}}^2$. 
 We consider $\bigtriangleup abc$ with a mesh $\delta/\sqrt{2}$ defined above.
For each $z\in \bigtriangleup abc$, we use the rotation map $\psi(z)=e^{-i\pi/4} z$ to map 
$\bigtriangleup abc$ into $\bigtriangleup \alpha\beta\gamma$.  
We select $X \in \overline{ab}$ (see Figure 1) for the $X$ in (1.2.5). Note that  $\bigtriangleup abc$ can be viewed to be  $\bigtriangleup  \alpha \beta\gamma$
after rotating $\bigtriangleup abc$ at the origin $\pi/4$ counterclockwise. 
For $a=0$, $b={1+i\over \sqrt{2}}$, and $c=i/\sqrt{2}$, after the rotation $\psi$, each square in 
$\bigtriangleup abc$ will be a rotated square.  We stretch each rotated square vertically to be a $\delta$-parallelogram  in ${\bf N}_\delta$. An open path from $\overline{aX}$ to $\overline{bc}$ will be an open path from $\overline{AX}$ to $\overline{BC}$ and the probability will remain the same (see Fig. 1). With this observation together with the Theorem
 $$\pi_\delta(\bigtriangleup abc, \overline{a X}, \overline{bc})=\pi_\delta(\bigtriangleup ABC,  \overline{AX}, \overline{BC})=X.\eqno{(1.2.12)}$$

Since $\psi$ is a  rotation map,  
$$\psi(X)=X.\eqno{(1.2.13)}$$
Let $\phi$ be the map from $\bigtriangleup \alpha\beta\gamma$ to the unit equilateral triangle $\bigtriangleup  ABC$ defined before.
By (1.2.5),  for the selected  $X$,
$$\phi(X) \neq X\mbox{ for }X\in \overline{\alpha\beta}.\eqno{(1.2.14)}$$
On the other side,  the conformal map $\phi\psi$ is from $\bigtriangleup abc$ to $\bigtriangleup ABC$, so if Cardy's formula holds, then
by (1.2.12) for the $X\in \overline{ab}$,  
$$\phi(X)=\phi\psi(X)=X.\eqno{(1.2.15)}$$
Therefore, (1.2.12) and (1.2.15) cannot hold together, so Cardy's formula does not hold on ${\bf Z}^2_\delta$.  
We summarize the result as the following corollary.\\


{\bf Corollary 2.} {\em Cardy's formula  does not hold for site and bond percolation on ${\bf Z}_\delta^2$.}\\

{\bf Remark 2.}  We can give many examples on periodic graphs such that Cardy's formula  does not hold. In fact, we believe that
Cardy's formula only holds on periodic graphs with the vertices on ${\bf T}_\delta$.

\section{ Separating probability and its differentiability on equilateral.} 
\subsection{Probability estimates of arms. } 
We first introduce $k$-arm paths. 
Let ${\cal Q}_k(r, R)$ be the event that there exist $i$ disjoint  $\bullet$-open paths 
and $j$ disjoint closed $\circ$-paths with $i+j=k$ for all $i\geq j\geq 1$ with $|i-j|\leq 1$
from  $[-r, r]^2$   to $\partial [-R, R]^2$ for $l^{-1}\leq r < R$.
We say that there are $k$-arm-paths. 
Here we assume that any two open $\bullet$-paths are separated by closed $\circ$-paths in the above $k$-arms.
(see a detailed definition for the separation in Kesten (1987)). We also denote by 
${\cal H}_k(r, R)$ 
 the event where there are $k$-arm paths from
$[-r, r]^2$ to $\partial [-R, R]^2$ in the upper half-space.
It is believed by Aizenman, Duplantier, and Aharony (1999) that
$${\bf P}_{ l}({\cal Q}_k(r, R))= \left\{\begin{array}{cc} ({r\over R})^{(k^2-1)/12+o(1),} &{\mbox{for $k\geq 2,$}}\\
 ({r\over R})^{5/48}& { \mbox{for $k=1$}}
\end{array}
\right .\eqno{(2.1.1)}$$
and 
$${\bf P}_{l}({\cal H}_k(r, R))=  \left({r\over R}\right)^{k(k+1)/6+o(1)},
\eqno{(2.1.2)}$$
where $o(1)\rightarrow 0$ as $l \rightarrow \infty$.
It has been proved that (2.1.1) holds for the triangular  lattice (see Werner 2008).

It is well known (see Chapter 3 in Kesten (1982)) that there exists  $d_1 >0$ such that
$${\bf P}_{ l}({\cal Q}_1(r, R))\leq  \left({r\over R}\right)^{d_1}.\eqno{(2.1.3)}$$
We want to remark that  one arm path is either one open $\bullet$-path or one closed $\circ$-path and $d_1$  is a positive constant independent of other parameters. In addition, $d_2$ in (2.1.5) is also positive constants independent of other parameters. 
For  precise estimates, 
(2.1.1) has been proved by Kesten, Sidoravicius, and Zhang (1998) for $i\geq 1$ and $j\geq 1$
with $i+j=5$.
Indeed, they showed that there exist $c_1$ and $c_2$ such that
$$c_1 \left({r\over R}\right)^2\leq {\bf P}_{l}({\cal Q}_5(r, R))\leq c_2\left({r\over R}\right)^2.\eqno{(2.1.4)}$$
In this paper,  $c_i$ for $i=1,2,\cdots$ is always a positive  constant that does not depend on  $l$, $\epsilon$, $\delta$, $r$, and $R$.  We write $c_i=c_i(M)$ if $c_i$ depends on $M$.
The value  of $c_i$ is  not significant, and  it changes from 
appearance to appearance. 
By (2.1.3) and (2.1.4) together with   Reimer's inequality,
$${\bf P}_{ l}({\cal Q}_6(r, R))\leq  c_1 \left({r\over R}\right)^{2+d_2}.\eqno{(2.1.5)}$$
It is also well known that (see Higuchi, Sekei, and Zhang (2012))
$$c_1\left({r\over R}\right)^{k(k+1)/6}\leq {\bf P}_{l}({\cal H}_k(r,R))\leq c_2\left({r\over R}\right)^{k(k+1)/6},
\eqno{(2.1.6)}$$
for $i=1$ and $j=1$, and for $i\geq 1$, $j\geq 1$ with $i+j=3$. They even showed (2.1.6) 
for a dependent percolation model.     For any straight-line  $L$ passing the origin divided $[R, R]^2 \setminus [-r, r]^2$ into two  parts.
We only focus on one of these two parts and consider $k$ arm paths from $[-r, r]^2$ to $\partial [R, R]^2$ in the part. We may
still denote by event $ {\cal H}_k(r,R))$ for the existence of $k$ arm paths in the part above $L$.  The same proof of (2.1.6) shows that for $k=2,3$
$$c_1\left({r\over R}\right)^{k(k+1)/6}\leq {\bf P}_{l}({\cal H}_k(r,R))\leq c_2\left({r\over R}\right)^{k(k+1)/6},
\eqno{(2.1.7)}$$

For any two rays from the origin divide $[R, R]^2 \setminus [-r, r]^2$ into two  parts. We suppose that the angle between two rays is $\theta$. We may
denote by event $ {\cal R}_k(\theta, r,R))$ for the existence of  $k$-arm paths from $[-r, r]^2$ to $\partial [R, R]^2$ in this part.  When $\theta=\pi$, the estimate of two arm paths is known in (2.1.7). By using  Theorem 1 in Kesten and Zhang (1987),  there exists $c_1= c_1(\theta) >0$  and $d_3 =d_3(\theta)$ with $d_3(\theta) <0$ for $\theta > \pi$ or $d_3(\theta) > 0$ for 
$\theta < \pi$ such that
$$c_1\left({r\over R}\right)^{1+d_3 }\leq {\bf P}_{l}({\cal R}_2(\theta, r,R))\leq c_2\left({r\over R}\right)^{1+d_3 },
\eqno{(2.1.8)}$$
We want to point out that Kesten and Zhang (1987) only showed in their Theorem 1 for one arm path. However, their proof implies the inequalities of (2.18) for multiple arm paths. For the triangle lattice, the exponent in (2.1.8) can be precisely   computed.
Finally, we introduce a standard topology result (see  Proposition 2.2 in Kesten (1982) or Proposition 11.2 in Grimmett (1999)).\\

{\bf Circuit lemma.} {\em  If $G$ is a finite  subgraph of ${\bf N}_\delta$ that is $\bullet$-connected. There exists a unique $\circ$-circuit  $\sigma (G)$ of ${\bf N}_\delta^*$,
containing $G$  in its interior such that every vertex of $\sigma (G)$ is $\bullet$-adjacent to  $G$.  Similarly,  If $G$ is a finite  subgraph of ${\bf N}_\delta^*$ that is $\circ$-connected. There exists a unique $\bullet$-circuit  $\sigma (G)$ of ${\bf N}_\delta$
containing $G$  in its interior such that  every vertex of $\sigma (G)$ is $\circ$-adjacent to  $ G$.}\\

\subsection{ Separating probabilities on equilateral.}
Given a vertex set $G$, we  denote by $\tau^{\pm k}(G)$ the vertex set by shifting $G$ horizontally $kl^{-1}$ to the right or to the left. 
In addition, for event ${\cal E}$, we denote by $\tau^{\pm k}({\cal E})$ when we shift $k{l^{-1}}$ for each configuration in ${\cal E}$ to the right or to the left.
For simplicity, we denote by $\tau^{\pm}(G)$  if $k=1$.
Now we consider $\bigtriangleup ABC$ with its three vertices $A=0$, $C=e^{i\pi/3}$,  and $B=1$. 
$\overline{AB}$, $\overline{AC}$ and $\overline{BC}$ are  divided into $(l-1)$ many intervals with length $l^{-1}$ such that the end points of these intervals  are vertices in ${\bf N}_{l^{-1}}$ for a large $l$ (see Fig. 2).   
Thus, $A$, $B$ and $C$ are vertex points.
The vertices of each $l^{-1}$-parallelogram are called {\em lattice points}. For simplicity,  we write $Z_1 < Z_2$  for two lattice points $Z_1=(z_1, y)$ and $Z_2=(z_2, y)$ with $z_1 < z_2$. In other words, we treat $Z_1$ and $Z_2$ as two numbers if they are in a horizontal line. 
\begin{figure}
\begin{center}
\begin{picture}(0,210)(40,70)

\setlength{\unitlength}{0.0125in}%
\begin{tikzpicture}
\thicklines
\begin{scope}[>={Stealth[black]},
              every edge/.style={draw=blue,very thick}]
  \path [-] (5.5,3.1) edge [right=40](6.5,4.6);
   \path [-] (6.4, 4.5) edge [right=40](5.6,5.75);
   \path [-] (5.6, 5.75) edge [right=40](7.4,8.50);
   
\end{scope}
\begin{scope}[>={Stealth[black]},
              every edge/.style={draw=red,very thick}]
 \path [-] (5,1) edge [right=0](5,1);
 \path [-] (3.83,3.1) edge [right=40](4.75,4.5);
 \path [-] (4.75,4.5) edge [right=40](3.97,5.65);
 \path [-] (3.85,5.85) edge [right=40](3.85,8.5);
 \path [-] (3.8,8.5) edge [right=40](3,9.7);

\end{scope}

\put(-0,100){\line(1,0){390}}
\put(-0,100){\line(2,3){195}}
\put(55,100){\line(2,3){167}}
\put(110,100){\line(2,3){139}}
\put(165,100){\line(2,3){112}}
\put(220,100){\line(2,3){85}}
\put(275,100){\line(2,3){58}}
\put(330,100){\line(2,3){30}}

\put(55,100){\line(-2,3){27}}
\put(110,100){\line(-2,3){55}}
\put(165,100){\line(-2,3){82}}
\put(220,100){\line(-2,3){110}}
\put(275,100){\line(-2,3){137}}
\put(330,100){\line(-2,3){165}}

\put(390,100){\line(-2,3){195}}

\put(220,100){\circle*{8}}
\put(247,142){\circle*{8}}
\put(220,180){\circle*{8}}
\put(248,222){\circle*{8}}
\put(278,267){\circle*{8}}

\put(167,100){\circle*{8}}
\put(194,143){\circle*{8}}
\put(166,180){\circle*{8}}
\put(166,267){\circle*{8}}
\put(136,305){\circle*{8}}

\put(-10,90){$A=0$}
\put(400,90){$B=1$}
\put(200,400){$C=e^{i\pi/3}$}
\put(215,85){$Z$}
\put(160,85){$\bar{Z}$}
\put(170,180){$\gamma_2$}
\put(235,190){$\gamma_1$}
\end{tikzpicture}
\end{picture}
\end{center}
\caption {\em  \small{$Z$ is separated by open path $\gamma_1$ from $\overline{AC}$. Event ${\cal H}(D, \bar{Z}, Z)$ is that there are two arm paths: one $\circ$-closed path from  $\bar Z$ to $\overline{AC}$   and  one open $\bullet$-path from $Z$ to $\overline{BC}$. The open $\bullet$-paths and the closed $\circ$-paths  are colored by blue and red, respectively, in this graph and in all the following graphs. }}
\end{figure}
For each $Z\in \overline{AB}\cap {\bf N}_{l^{-1}}$, we denote by ${\cal E}(Z)$  the event that there is an open $\bullet$-path, called separating $Z$ from $\overline {AC}$, from $\overline{BC}$ to $\overline{AZ}$. Let 
$$f_l( Z)={\bf P}_l({\cal E}(Z))\eqno{}$$
 for the above vertices $Z\in \overline{AB}\cap {\bf N}_{l^{-1}}$. 
 It follows from this definition that for any $Z\in \overline{AB}\cap {\bf N}_{l^{-1}}$,
 $$ f_l(Z)= \pi_{l^{-1}}(Z).\eqno{(2.2.1)}$$
 Therefore, $f_l$ is a function defined on all the lattice points of  $Z\in \overline{AB}\cap {\bf N}_{l^{-1}}$.

Let us consider $f(Z)- f(\bar{Z})$ for  $Z,\bar{Z}\in \overline{AB}\cap {\bf N}_{l^{-1}}$  with $d(Z, \bar{Z})=l^{-1}$, where  $Z$ and $\bar{Z}$ are called an  {\em adjacent $\circ$-pair}.  
If $\bar Z$ and $Z$ are $\circ$-adjacent with $\bar Z < Z$, then 
$${\cal E}(\bar{Z})\subset {\cal E}(Z).\eqno{(2.2.3)}$$
By (2.2.3), 
$$f_l(Z)- f_l(\bar{Z})={\bf P}_l({\cal E}(Z)\setminus {\cal E}(\bar{Z})).\eqno{(2.2.4)}$$ 
Let ${\cal H}(D, \bar{Z}, Z)$  for $D=\bigtriangleup AB C$ be the event that there exist an open $\bullet$-path from $Z$ to $\overline{B C}$ and a  $\circ$-closed  path
from $\bar Z$ to $\overline {AC}=$ inside $D$, called {\em two arm paths} in the upper half-space  above $\overline{AB}$ from $\bar{Z}$ and from $Z$ (see Fig. 2).
 It is well known (see Claim 11 in Bollobas and Riordan (2006)) that 
 $$f_l(Z)- f_l(\bar{Z})={\bf P}_l({\cal E}(Z)\setminus {\cal E}(\bar{Z}))={\bf P}_l({\cal H}(D, \bar{Z}, Z)).\eqno{(2.2.5)}$$
 We want to remark that the above Claim 11 is proved for the site percolation on ${\bf T}_{l^{-1}}$. The same proof of Claim 11 can be used to show (2.2.5) on ${\bf N}_{l^{-1}}$.
 It follows from (2.2.6) and (2.1.3) that there exists $d_1>0$  defined in (2.1.3) independent of $l$  such that for any $\bar{Z}, Z\in \overline{ AB}\cap {\bf N}_{l^{-1}}$ with $\bar Z < Z$,
 $$0\leq f_l(Z)-f_l(\bar{Z})\leq l^{-d_1}.\eqno{(2.2.6)}$$
 By (2.2.6), 
 $$ f_l(0)\leq l^{-d_1}\mbox{ and } f_l(1)=1.\eqno{(2.2.7)}$$

\subsection{The  difference of quotient for $f_l(Z)$. }
For a  small constant $\iota>0$ independent of $l$, let
$$[A_\iota, B_\iota]= \{U \in \overline{AB}: d(U, A) \geq \iota \mbox{ and } d(U, B) \geq \iota\}$$
is a closed interval inside $\overline{AB}$.  We select $\iota$ such that  $A_\iota, B_\iota\in {\bf N}_{l^{-1}}$.  
 By (2.1.5),  there exist $c_1=c_1(\iota)$  and $c_2=c_2(\iota)$ such that for any  two $\circ$-adjacent lattice points $\bar Z$ and $Z$ in $[A_\iota, B_\iota]$,
 $$c_1< {f_l(Z)-f_l(\bar{Z})\over l^{-1}}\leq c_2.\eqno{(2.3.1)}$$
 
We have shown that  the first order of  horizontal difference quotient $f_l(Z)$ is bounded on $[A_\iota, B_\iota]$ in (2.3.1).
 Now we investigate the second order of horizontal difference quotient of $f_l(Z)$ for $Z\in [A_\iota, B_\iota]$. 
 We consider $\circ$-adjacent triples $ Z$, $ \bar{Z}$, and $\hat{Z}$ in  $[A_\iota, B_\iota]\cap  {\bf N}_{l^{-1}}$   with $\hat{Z} < \bar{Z} < Z$.  
$\tau^+(D)$  is the graph that $D=\bigtriangleup ABC$ is moved to the right one unit $l^{-1}$, so $\tau^+(\hat{Z})$ and $\tau^+(\bar{Z})$ will be $\bar{Z}$ and $Z$.
Recall that ${\cal H}(D, \bar{Z}, Z)$ is the event that there are two arm paths from $\bar{Z}$ and $Z$ in $D$, so
${\cal H}(\tau^+(D), \bar{Z}, Z)$ is the event that there are  two arm paths from $\tau^+(\hat{Z})=\bar{Z}$ and $\tau^+(\bar{Z})=Z$ in $\tau^+(D)$.
By (2.2.3), (2.2.4), and translation invariance, the second order of  the difference 
\begin{eqnarray*}
&&f_l(Z)-f_l(\bar{Z})-(f_l(\bar{Z})-f_l(\hat{Z})) ={\bf P}_l({\cal H}(D, \bar{Z}, Z))-{\bf P}_l({\cal H}(D, \hat{Z}, \bar{Z}))\\
&=&{\bf P}_l({\cal H}(D, \bar{Z}, Z))-{\bf P}_l({\cal H}(\tau^+(D), \bar{Z}, {Z}))\\
&=&{\bf P}_l({\cal H}(D, \bar{Z}, Z)\cap {\cal H}^C(\tau^+(D), \bar{Z}, {Z}))-{\bf P}_l({\cal H}^C(D, \bar{Z}, Z)\cap {\cal H}(\tau^+(D), \bar{Z}, {Z})),
\hskip 1.5 cm (2.3.2)
\end{eqnarray*}
where ${\cal H}^C(\tau^+(D), \bar{Z}, {Z})$ is the complement of ${\cal H}(\tau^+(D), \bar{Z}, {Z})$.
\begin{figure}
\begin{center}
\begin{picture}(70,200)(40,70)

\setlength{\unitlength}{0.0125in}%
\begin{tikzpicture}
\thicklines
\begin{scope}[>={Stealth[black]},
              every edge/.style={draw=blue,very thick}]
 \path [-] (5,1) edge [right=0](5,1);
  \path [-] (6,3.1) edge [right=40](9.3, 8.5);
\path [-] (4.8, 3.2) edge [right=40](2.3,8.7);
\path [-] (2.3, 8.7) edge [right=40](3.6,10.6);
\path [-] (3.6, 10.6) edge [right=40](6,6.5);
\path [-] (6, 6.5) edge [right=40](8.2, 10.1);

\end{scope}
\begin{scope}[>={Stealth[black]},
              every edge/.style={draw=red,very thick}]
 \path [-] (5,1) edge [right=0](5,1);
 \path [-] (5.5,3.2) edge [right=40](3.1,8.6);

\end{scope}
\put(-0,100){\line(2,3){210}}
\put(30,100){\line(2,3){210}}

\put(420,100){\line(-2,3){210}}
\put(450,100){\line(-2,3){210}}

\put(0,100){\line(1,0){450}}

\put(235,100){\circle*{8}}

\put(218,100){\circle*{8}}

\put(195,100){\circle*{8}}
\put(139,305){\circle*{8}}

\multiput(65, 230)(50,70){1}{\framebox(160,160)}
\multiput(145, 100)(50,70){1}{\framebox(150,110)}

\put(-10,90){$A$}
\put(15,90){$\tau^+(A)$}
\put(450,90){$\tau^+(B)$}
\put(410,90){$B$}
\put(250,410){$\tau^+(C)$}
\put(190,410){$C$}
\put(130,270){$S$}

\put(235,85){$Z$}
\put(210,85){$\bar Z$}
\put(190,180){$\gamma_2$}
\put(285,150){$\gamma_1$}
\end{tikzpicture}
\end{picture}
\end{center}
\caption{ \em \small{On ${\cal H}(\tau^+(D),\bar{Z}, Z)$, there are two arm paths from $\bar Z$ and $Z$ to $\tau^+(\overline{AC})$ and to $\tau^+(\overline{BC})$, respectively. On $\{{\cal H}(D, \bar{Z}, Z)\}^C$,    the  closed  $\circ$-path $\gamma_2$ from $\bar Z$ cannot go further to $\overline{AC}$ such that $\{{\cal H}( D,\bar{Z}, Z)\}^C$ occurs. By the circuit lemma, there is an open $\bullet$-path 
 from $\overline{ \tau^+(A) Z}$ meeting a vertex $S$ at $\tau^+(\overline{AC})$  and then meeting $\tau^+(\overline{BC})$ in $\tau^+(\bigtriangleup ABC)$, denoted by ${\cal D}^-( \bigtriangleup ABC, \tau^+( \bigtriangleup ABC),\bar{Z})$.}}
\end{figure}

We first estimate  the second term in (2.3.2) ${\bf P}_l({\cal H}(\tau^+(D), \bar{Z}, {Z}) \cap {\cal H}^C(D, \bar{Z}, Z))$. 
On ${\cal H}(\tau^+(D), \bar{Z}, {Z})$, $\hat{Z}$ and $ \bar{Z}$ are shifted to be $\bar{Z}$ and $Z$,
and  $A, B,C$ are shifted to be $\tau^+(A), \tau^+(B), \tau^+(C)$, so there are two arm paths $\gamma_1$ and $\gamma_2$  from $Z$ and $\bar Z$ to $\tau^+(\overline{BC})$ and $\tau^+(\overline{AC})$ in $\tau^+(D)$ (see Fig. 3).
  Moreover,  on $({\cal H}(D, \bar{Z}, Z))^C$,  the complement  of ${\cal H}(D, \bar{Z}, Z)$,  either  no open $\bullet$-paths go to $\overline{BC}$
  from $Z$ or no closed $\circ$-paths go to $\overline{AC}$ from $\bar Z$.  The first case cannot hold since $\gamma_1$ exists. 
  On the above second case, 
   $\gamma_2$ goes to $\tau^+(\overline{AC})$, but cannot go further to $\overline{AC}$.
  If the above closed $\circ$-path $\gamma_2$ in $\tau^+(D)$ cannot go further to $\overline{AC}$, then the closed $\circ$-cluster of $\bar Z$  inside $\tau^+(D)$ has a boundary vertex $S$ at $\overline{AC}$ and all its boundary vertices stay inside $D\cup \tau^+(D)$ (see Fig. 3).  By the circuit lemma,  we know that the boundary inside $ D\cup \tau^+(D)$ consists of  an open $\bullet$-path 
  from a point in $\overline{A \bar{Z}}$ to a vertex of $\overline{AC}$ 
  and then to $ \tau^+(\overline{BC})$ (see Fig. 3).  In other words, this open $\bullet$-path is blocked, called  {\em block property}, $\gamma_2$ to reach  $\overline{AC}$. 
  We denote this event of the existence of  the above open $\bullet$-path  by 
  $${\cal D}^-( \bigtriangleup ABC, \tau^+(\bigtriangleup ABC), \bar{Z})= {\cal D}^-( D, \tau^+(D), \bar{Z}).$$
  More precisely, there is an open $\bullet$-path from $\overline{ A\bar Z}$ to the left side of $D$ and then to the right side of  $\tau^+(D)$.
  Thus, 
  $${\cal H}^C(D, \bar{Z}, Z)\cap {\cal H}(\tau^+(D), \bar{Z}, {Z})\subset {\cal H}(\tau^+(D), \bar{Z}, Z)\cap {\cal D}^-(D, \tau^+(D),   \bar{Z}).$$
  On the other hand, by the circuit lemma,  if ${\cal D}^-(D, \tau^+(D), \bar{Z}) $ occurs, then the closed $\circ$-cluster containing 
  $\gamma_2$ has to stay inside $\tau^+(D)$, so 
   the closed $\circ$-path $\gamma_2$ from $\bar Z$ cannot  reach to $\overline{AC}$.  This implies that ${\cal H}^C(D, \bar{Z}, Z)$ occurs.
   Thus, we also have
  $${\cal H}(\tau^+(D), \bar{Z}, {Z})\cap {\cal D}^-(D, \tau^+(D), \bar{Z})\subset {\cal H}^C(D, \bar{Z}, Z)\cap {\cal H}(\tau^+(D), \bar{Z}, {Z}).$$
  By the above two  inequalities, 
  $${\cal H}^C(D, \bar{Z}, Z)\cap {\cal H}(\tau^+(D), \bar{Z}, {Z})= {\cal H}(\tau^+(D), \bar{Z}, {Z})\cap {\cal D}^-(D, \tau^+(D), \bar{Z}).\eqno{(2.3.3)}$$
  
On the above event ${\cal D}^-(D, \tau^+(D), \bar{Z})$,  note that $\gamma_2$ $\circ$ adjacent to $\overline{AC}$,  so it has to $\circ$-adjacent to its open $\bullet$-circuit boundary. In other words,  the blocked open $\bullet$-path has a point $S\in \overline{AC}$ $\circ$-adjacent to $\gamma_2$ (see Fig. 3).
We construct a square $B(S, jl^{-1})$  with the center at $S$ and a side length $jl^{-1}$. We divide the following two cases: (1) 
$d(S, \tau^+(C))\leq \iota'$ and  (2) $d(S, \tau^+(C))> \iota'$  for some small
positive constant $\iota'$. For case (2), 
there are three arm arm below $\overline{AC}$ from 
$S$ to $\partial B(S, \iota')$, where two open $\bullet$-paths  and  a closed $\circ$-path, a part of $\gamma_2$,  are    from $S$ to $\partial B(S, \iota')$.
In addition, 
there are the other two arm paths from $\bar Z$ and $Z$ to  the boundary of 
$Z+[-\iota', \iota']^2$  that is  the disjoint square from $S+[\iota', \iota']^2$.  By taking $\iota'$ small but positive, these three arm paths and the two arm paths are disjoint. Thus, by using (2.1.7) for $k=3$ and $k=2$,
$${\bf P}_l( {\cal H}(\tau^+(D), \bar{Z}, Z)\cap {\cal D}^-(D, \tau^+(D), \bar{Z}),\mbox{ case (2)})\leq c_1 l^{-1} l^{-1}\leq c_2 l^{-2}.\eqno{(2.3.4)} $$

Now we focus on case (1).
We suppose that $d(S, \tau^+(C))=jl^{-1}$ for $j l^{-1}\leq \iota'$. Thus, there are three arm paths below $\overline{AC}$ from 
$S$ to $\partial B(S, jl^{-1})$, where two open $\bullet$-paths  and  a closed $\circ$-path, a part of $\gamma_2$,  are  from $S$ to $\partial B(S,jl^{-1})$.
In addition, there are two arm paths from $\partial B(S, jl^{-1})$ to $\partial(S+[-\iota', \iota']^2)$ in $D=\bigtriangleup ABC$ (see Fig. 3). 
Note that the two arm paths are in the two rays with the angle $\pi/3 < \pi$.
By taking $\iota' =O(1)>0$ small
and  using (2.1.7) and (2.1.8), if we sum all $j$, then we can show that  there exists $c_1=c_1(\iota)$ and $c_2=c_2(\iota')$ such that
\begin{eqnarray*}
&&\!\!\!\sum_{j=1}^{\iota' l}{\bf P}_l(\exists\small \mbox{ three arm paths from $S$ to $\partial B(S, jl^{-1})$ and  two arm paths in $D$ from }B(S, jl^{-1}) \mbox{ to }\partial (S+[-\iota', \iota']^2))\\
&&\leq c_1  \sum_{j=1}^{\iota'l} j^{-2}( j l^{-1})^{1+d_3}  \leq c_2l^{-1-d_3} l^{d_3} \leq c_1l^{-1}.\hskip 8cm {(2.3.5)}
\end{eqnarray*}
By using (2.1.6) for the above two arm paths  in $Z+[-\iota', \iota']^2$, together with (2.3.5),  note  that two squares  are disjoint, so 
there exists a positive constant $c_1=c_1(\iota)$ independent of $l$ such that  
$${\bf P}_l( {\cal H}(\tau^+(D), \bar{Z}, Z)\cap {\cal D}^-(D, \tau^+(D), \bar{Z}),\mbox{ case (1)} )\leq c_1l^{-1}l^{-2}\leq c_1 l^{-2}.\eqno{(2.3.6)} $$
Together with (2.3.4) and (2.3.6),
$${\bf P}_l( {\cal H}(\tau^+(D), \bar{Z}, Z)\cap {\cal D}^-(D, \tau^+(D), \bar{Z}))\leq c_1 l^{-2}.\eqno{(2.3.7)} $$

\begin{figure}
\begin{center}
\begin{picture}(70,200)(40,70)

\setlength{\unitlength}{0.0125in}%
\begin{tikzpicture}
\thicklines
\begin{scope}[>={Stealth[black]},
              every edge/.style={draw=blue,very thick}]
 \path [-] (5,1) edge [right=0](5,1);
  \path [-] (6.7,3.1) edge [right=40](9.1,7.3);

\end{scope}
\begin{scope}[>={Stealth[black]},
              every edge/.style={draw=red,very thick}]
 \path [-] (5,1) edge [right=0](5,1);
 \path [-] (6,3.1) edge [right=40](1.5,7.5);
 
 \path [-] (3.8, 5.3) edge [right=40](7.4,11.5);
  \path [-] (7.4,11.5) edge [right=40](9.9,7.7);
 \path [-] (7.3,3.1) edge [right=40](9.9,7.7);

\end{scope}

\put(-0,100){\line(2,3){210}}
\put(30,100){\line(2,3){210}}

\put(420,100){\line(-2,3){210}}
\put(450,100){\line(-2,3){210}}

\put(0,100){\line(1,0){450}}

\put(255,100){\circle*{8}}
\put(277,100){\circle*{8}}
\put(233,100){\circle*{8}}
\put(348,254){\circle*{8}}


\put(-10,90){$A$}
\put(15,90){$\tau^+(A)$}
\put(450,90){$\tau^+(B)$}
\put(410,90){$B$}
\put(250,410){$\tau^+(C)$}
\put(190,410){$C$}
\put(215,375){$C''$}
\put(352,257){$S$}
\put(250,83){$Z$}
\put(230,83){$\bar {Z}$}
\put(180,160){$\gamma_2$}
\put(275,150){$\gamma_1$}
\end{tikzpicture}
\end{picture}
\end{center}
\caption{\em  \small{On ${\cal H}(D, \bar{Z}, Z)$, there are two arm paths from $\bar Z$ and $Z$  to $\overline{AC}$ and to $\overline{BC}$, respectively. On $ \{{\cal H}(\tau^+(D), \bar{Z}, {Z})\}^C$,
$\gamma_1$ cannot go further to $\tau^+(\overline{BC})$.  By the circuit lemma, there is a closed $\circ$-path in $D\cup\tau^+(D)$ from $\overline{ ZB}$ meeting a vertex  $S\in \tau^+(\overline{BC})$  and 
then to $\overline{AC}$. This event is denoted by ${\cal D}^+(D, \tau^+(D), Z)$. }}
\end{figure}

Now we estimate ${\bf P}_l({\cal H}(D,\bar{Z}, Z)\cap  {\cal H}^C(\tau^+(D), \bar{Z}, {Z}))$ for the above $ \hat{Z}, \bar Z, Z  \in [A_\iota, B_\iota]\cap {\bf N}_{\l^{-1}}$.
On ${\cal H}(D, \bar{Z}, Z)$,  there are two arm paths $\gamma_1$ and $\gamma_2$  from $\bar Z$ and $Z$ in $D$ to $\overline{AC}$ and to $\overline{BC}$, respectively. On $ \{{\cal H}(\tau^+(D), \bar{Z}, {Z})\}^C$, $\gamma_2$ does not exist or 
$\gamma_1$ cannot go further to $\tau^+(\overline{BC})$ (see Fig. 3).   Thus, the only second case occurs.
If $\gamma_1$ cannot go further to $\tau^+(\overline{BC})$, then
the open connected cluster of ${Z}$  inside $D$ has a boundary vertex $S$ at $ \tau^+(\overline{BC)}$ and all its boundary stays inside $D\cup \tau^+(D)$ (see Fig. 4). By the circuit lemma,  the boundary inside $D\cup \tau^+(D)$ consists of  a closed $\circ$-path from  a point in $ \overline{Z B}$ to $S$  and then to   $\overline{AC}$ (see Fig. 4). We denote  this event  of the existence of the above closed $\circ$-path
by ${\cal D}^+(D, \tau^+(D), Z)$.   More precisely,  there is a closed $\circ$-path from $ \overline{Z B}$ to $S$ and then to $\overline{AC}$ in $ D\cup \tau^+(D)$. 
  Thus, 
  $${\cal H}(D,\bar{Z}, Z)\cap  {\cal H}^C(\tau^+(D), \bar{Z}, {Z})\subset {\cal H}(D, \bar{Z}, Z)\cap {\cal D}^+(D, \tau^+(D),   {Z}).\eqno{(2.3.8)}$$
  On the other hand, by the circuit lemma,  if ${\cal D}^+(D, \tau^+(D), {Z}) $ occurs, then the open $\bullet$-cluster containing $\gamma_1$  will stay inside $D$, so 
  the  open $\bullet$-path $\gamma_1$ from $ Z$ cannot  reach to $\tau^+(\overline{BC})$.  This implies that
  ${\cal H}^C(\tau^+(D), \bar{Z}, {Z})$ occurs.
  Thus, we also have
  $${\cal H}(D, \bar{Z}, Z)\cap {\cal D}^+(D, \tau^+(D), {Z})\subset {\cal H}(D,\bar{Z}, Z)\cap  {\cal H}^C(\tau^+(D), \bar{Z}, {Z}).$$
  By the above two inequalities, 
  $${\cal H}(D, \bar{Z}, Z)\cap {\cal H}^C(\tau^+(D), \bar{Z}, {Z})={\cal H}(D, \bar{Z}, Z)\cap {\cal D}^+(D, \tau^+(D),   {Z}).\eqno{(2.3.9)}$$
  On the above event ${\cal H}(D,\bar{Z}, Z)\cap {\cal D}^+(D, \tau^+(D), {Z})$, 
  by using the same estimate of  (2.3.7),
there exists a positive constant $c_1=c_1(\iota,\eta,  \iota')$ independent of $l$ such that  
$${\bf P}_l({\cal H}(D,\bar{Z}, Z)\cap  \{{\cal H}(\tau^+(D), \bar{Z}, {Z}) \}^C)={\bf P}_l({\cal H}(D, \bar{Z}, Z)\cap {\cal D}^+(D, \tau^+(D),   {Z}))\leq c_1l^{-2}.\eqno{(2.3.10)}$$

Thus, we have 
\begin{eqnarray*}
&&f_l(Z)-f_l(\bar{Z})-(f_l(\bar{Z})-f_l(\hat{Z})) ={\bf P}_l({\cal H}(D, \bar{Z}, Z))-{\bf P}_l({\cal H}(D, \hat{Z}, \bar{Z}))\\
&=&{\bf P}_l({\cal H}(D, \bar{Z}, Z))-{\bf P}_l({\cal H}(\tau^+(D), \bar{Z}, {Z}))\\
&=&{\bf P}_l({\cal H}(D, \bar{Z}, Z)\cap {\cal H}^C(\tau^+(D), \bar{Z}, {Z}))-{\bf P}_l({\cal H}^C(D, \bar{Z}, Z)\cap {\cal H}(\tau^+(D), \bar{Z}, {Z}))\\
&=&{\bf P}_l( {\cal H}(D, \bar{Z}, Z)\cap {\cal D}^+(D, \tau^+(D), {Z}))-{\bf P}_l({\cal H}(\tau^+(D), \bar{Z}, Z)\cap {\cal D}^-(D, \tau^+(D),   \bar{Z})).
\hskip 1 cm (2.3.11)
\end{eqnarray*}
By using the estimates (2.3.7) and (2.3.10), 
 for  $\circ$-adjacent triple 
 $ \hat{Z}< \bar Z <Z  \in [A_\iota, B_\iota]\cap {\bf N}_{l^{-1}}$,  there is $c_1=c_1(\iota)$ such that
 $$|f_l(Z)-f_l(\bar{Z})-(f_l(\bar{Z})-f_l(\hat{Z}))| \leq c_1 l^{-2}.\eqno{(2.3.12)}$$
 
 
 We also consider $\circ$-adjacent lattice pairs $Z_1, Z_2$ and $Z_3, Z_4$ in $ [A_\iota, B_\iota] \cap {\bf N}_{l^{-1}}$  with $Z_2 < Z_3$ and $d(Z_2, Z_3)= \epsilon$ for  
 for some small positive  $\epsilon$. 
 Let $W_1, W_1, \cdots, W_n$ be the $\circ$-adjacent lattice points in $ [A_\iota, B_\iota] \cap {\bf N}_{l^{-1}}$ with
 $W_1=Z_1 < W_2=Z_2< Z_3<  \cdots < W_{n-1}= Z_3< W_n=Z_4$ for $n=\epsilon l$.   Here we assume that $\epsilon l$   is an integer,
 otherwise we can use $\lfloor \epsilon l\rfloor$.
 By (2.3.12)
 $$|f_l(Z_4)-f_l(Z_3)-(f_l(Z_2)-f_l(Z_1)) |\leq \sum_{i=1}^{n-2} |f_l(W_{i+2})-f_l(W_{i+1})-( f_l(W_{i+1})- f_l(W_i))|\leq  c_1 \epsilon.\eqno{(2.3.13)}$$


 \subsection{The third order  difference  of quotient for $f_l(Z)$.}
In  this section,  we will investigate  the third order difference   of quotient  $f_l(U)$ on $\bigtriangleup ABC=D$.
We consider the $\circ$-adjacent   lattice points $Z_0, Z_1, Z_2, Z_3$ on $[A_\iota, B_\iota]$ with $Z_0<  Z_1< Z_2< Z_3$. By (2.3.11),
 \begin{eqnarray*}
 &&[(f_l(Z_3)- f_l(Z_2)) - (f_l(Z_2)- f_l(Z_2))]- [ (f_l(Z_2)- f_l(Z_1) - (f_l(Z_1)- f_l(Z_0))]\\
 &=&[{\bf P}_l( {\cal H}(D, Z_2, Z_3)\cap {\cal D}^+(D, \tau^+(D), Z_3))- {\bf P}_l( {\cal H}(\tau^+(D), Z_2, Z_3)\cap {\cal D}^-(D, \tau^+(D), Z_2))]
 \\
 &&-[{\bf P}_l( {\cal H}(D, Z_1, Z_2)\cap {\cal D}^+(D, \tau^+(D), Z_2))- {\bf P}_l( {\cal H}(\tau^+(D), Z_1, Z_2)\cap {\cal D}^-(D, \tau^+(D), Z_1))].\hskip .1cm (2.4.1)
  \end{eqnarray*}
  By (2.3.12), there exists $c_1 >0$ such that
  \begin{eqnarray*}
 &&|[ (f_l(Z_3)- f_l(Z_2)) - (f_l(Z_2)- f_l(Z_1))]- [ (f_l(Z_2)- f_l(Z_1) - (f_l(Z_1)- f_l(Z_0))]|\\
 &\leq &|{\bf P}_l( {\cal H}(D, Z_2, Z_3)\cap {\cal D}^+(D, \tau^+(D), Z_3))- {\bf P}_l( {\cal H}(\tau^-(D), Z_2, Z_3)\cap {\cal D}^-(D, \tau^+(D), Z_2))|\\
 &&+|{\bf P}_l( {\cal H}(D, Z_1, Z_2)\cap {\cal D}^+(D, \tau^+(D), Z_2))- {\bf P}_l( {\cal H}(\tau^+(D), Z_1, Z_2)\cap {\cal D}^-(D, \tau^+(D), Z_1))|\\
 &\leq &c_1 l^{-1 }.\hskip 13.5 cm (2.4.2)
  \end{eqnarray*}
  Thus, we have an estimate of the third  order difference of quotient
  $$\left| {f_l(Z_3)- f_l(Z_2) \over l^{-1}} - {f_l(Z_2)- f_l(Z_1)\over l^{-1}}- \left [ {f_l(Z_2)- f_l(Z_1) \over l^{-1}}- {f_l(Z_1)- f_l(Z_0)\over l^{-1}}\right]\right |\leq c_1 l^{-2}.\eqno{(2.4.3)}$$
  The upper bound in the right side of (2.4.3) is not good enough to show the theorem.  We need to improve a more tighter bound to be $o(l^{-2-\alpha})$ for some $\alpha >0$.  
  
  We  consider  two pairs of $\circ$-adjacent triples $ Z_1 < Z_2 < Z_3$ and $Z_4< Z_5 < Z_6$ on $[A_\iota, B_\iota]$ with $d(Z_3, Z_4) =nl^{-1}$ and
 $Z_3 < Z_4$ for small $nl^{-1}>0$.   By (2.4.1),
 \begin{eqnarray*}
 &&[ (f_l(Z_6)- f_l(Z_5)) - (f_l(Z_5)- f_l(Z_4))]- [ (f_l(Z_3)- f_l(Z_2)) - (f_l(Z_2)- f_l(Z_1))]\\
 &=&[{\bf P}_l( {\cal H}(D, Z_5, Z_6)\cap {\cal D}^+(D, \tau^+(D), Z_6))- {\bf P}_l( {\cal H}(D, Z_2, Z_3)\cap {\cal D}^+(D, \tau^+(D), Z_3))]
 \\
 &&-[{\bf P}_l( {\cal H}(\tau^+(D), Z_5, Z_6)\cap {\cal D}^-(D, \tau^+(D), Z_5))- {\bf P}_l( {\cal H}(\tau^+(D), Z_1, Z_2)\cap {\cal D}^-(D, \tau^+(D), Z_1))].\hskip .1cm (2.4.4)
 \end{eqnarray*}
 
 We  show the following lemma.\\

  {\bf Lemma  2.4.1.} {\em For  $Z_i$ with $i=1,2,3,4,5,6$ defined above and for any small $nl^{-1}>0$, there exist constants $\alpha>0$ and $c_1$ independent of $\epsilon$ and $l$  such that
  $$|{\bf P}_l( {\cal H}(D, Z_5, Z_6)\cap {\cal D}^+(D, \tau^+(D), Z_6))- {\bf P}_l( {\cal H}(D, Z_2, Z_3)\cap {\cal D}^+(D, \tau^+(D), Z_3))|\leq c_1(nl^{-1})^\alpha l^{-2}$$
  and 
  $$|{\bf P}_l( {\cal H}(\tau^+(D), Z_5, Z_6)\cap {\cal D}^-(D, \tau^+(D), Z_5))- {\bf P}_l( {\cal H}(\tau^+(D), Z_1, Z_2)\cap {\cal D}^-(D, \tau^+(D), Z_1))|\leq c_1 (nl^{-1}) ^\alpha l^{-2}.$$}

A small $nl^{-1}$ means that $nl^{-1}$ is a small fraction such as $1/100$. In fact, we only need $n=1$ and $nl^{-1} < \epsilon$ for any $\epsilon >0$ in Lemma 2.4.1 to show the Theorem. The proof of Lemmas 2.4.1 is  involved and long, so we show it in the Appendix.\\

  By Lemma 2.4.1, if $nl^{-1} \leq \epsilon$, then
   $$ \left|  {{f_l(Z_6)- f_l(Z_5)\over l^{-1}} -{f_l(Z_4)- f_l(Z_3)\over l^{-1}}\over l^{-1}}-  {{f_l(Z_3)- f_l(Z_2)\over l^{-1}} -{f_l(Z_2)- f_l(Z_1)\over l^{-1}}\over l^{-1}}\right |\leq c_1 \epsilon^{\alpha}.\eqno{(2.4.5)}$$
In particular, 
if  $n=1$, then  $Z_1, Z_2, Z_3, Z_4, Z_5, Z_6$ are $\circ$-adjacent lattice points on $[A_\iota, B_\iota]$. We have the third
 difference order of $f_l$.
 $$ \left|  {{f_l(Z_3)- f_l(Z_2)\over l^{-1}} -{f_l(Z_2)- f_l(Z_1)\over l^{-1}}\over l^{-1}}-  {{f_l(Z_2)- f_l(Z_1)\over l^{-1}} -{f_l(Z_1)- f_l(Z_0)\over l^{-1}}\over l^{-1}}\right |= o(l^{-\alpha}). \eqno{(2.4.6)}$$

\section{ A diagonal invariant shifting for separated probabilities and discrete contour  integrals.}
In section 2,  we defined the separated probabilities on $\overline{AB}$. Now we will define them on all the lattice points on $\bigtriangleup ABC$.
We consider $D=\bigtriangleup ABC$  and $A_\iota$ and $B_\iota$  with the mesh $l^{-1}$ defined before.   
We denote   two unit vectors by
$$\vec{e}_1=\langle 1/2,\sqrt{3}/2\rangle\mbox{ and } \vec{e}_2=\langle 1/2,-\sqrt{3}/2\rangle.$$
We consider  unit equilateral triangle $D_1=\bigtriangleup A_1B_1 C_1 $ by moving
$D=\bigtriangleup ABC$ one unit $l^{-1}$ in  the right and up  along $\vec{e}_1$ direction (see Fig. 5).
We also consider  unit equilateral triangles $D_1'=\bigtriangleup A_1'B_1' C_1' $ by moving  $\bigtriangleup ABC$ one unit $l^{-1}$ in the left  and up   along $\vec{e}_2$ direction (see Fig. 5). 
Moreover, we move $D_1$ one unit $l^{-1}$ in the left and up along $\vec{e}_2$ direction to obtain $D_2=\bigtriangleup A_2B_2C_2$. Similarly,  we move $\bigtriangleup A'_1 B_1'C'_1$ one unit $l^{-1}$ in the right and up along $\vec{e}_1$ direction to obtain 
$D_2'=\bigtriangleup A_2'B_2' C_2'$ (see Fig. 5). Thus
$$D_2=\bigtriangleup A_2B_2 C_2=D_2'=\bigtriangleup A_2'B_2'C_2'.\eqno{(3.1)}$$
We do the same construction from $D_2$,  as we did for $D_1$, $D_1'$,  $D_2$, $D_2'$ from $D$,  to obtain 
$$D_3=\bigtriangleup A_3B_3 C_3, D'_3=\bigtriangleup A_3'B_3'C_3', 
D_4=\bigtriangleup A_4B_4 C_4=D_4'=\bigtriangleup A_4'B_4'C_4'.\eqno{}$$
More precisely, , $D_2, D_3, D_3', D_4$ are obtained by moving $D, D_1, D_1', D_2$ up  $2 \sqrt{3} l^{-1}$ vertically.
We continue this vertical moving to construct $D_4, D_5, D_5', D_6, \cdots,$ $D_{2i}, D_{2i+1}, D_{2i+1}', D_{2i+2}$.
For convenience,  we denote by $D_0=D_0'= \bigtriangleup ABC$. In other words, $D_0$ is moved $l^{-1}$ in the right and up direction for $D_1$, $D_1$ is 
moved $l^{-1}$ in the left and up direction for $D_2$, $\cdots$, $D_{2i}$ is moved $l^{-1}$ in the right and up direction for $D_{2i+1}$ and $D_{2i+1}$ is moved $l^{-1}$ in the left and up direction for $D_{2i+2}$.  
The vertices between  horizontal lines $y=i$ and $y=i+1$ are called the $i$-{\em layer}. The set sequence 
 $\{D_i\}$ grows vertically in  a zag-zig way and  on alternative  layers.
Moreover,
$D_0'$ is moved $l^{-1}$ in the left and up direction for $D_1'$, $D_1'$ is moved $l^{-1}$ in the right and up direction for $D_2'$, $\cdots$, $D_{2i}'$ is moved $l^{-1}$ in the
 left and up direction for $D_{2i+1}'$ and $D_{2i+1}'$ is moved $l^{-1}$ in  the right and up direction for $D_{2i+2}'$.  $\{D_i\}$ and $\{D_i'\}$ grow
 vertically in a mirror symmetry about the vertical line $y=1/2$. 
 \begin{figure}
\begin{center}
\begin{picture}(100,260)(40,75)

\setlength{\unitlength}{0.0125in}%
\begin{tikzpicture}
\thicklines
\begin{scope}[>={Stealth[black]},
              every edge/.style={draw=blue,very thick}]
\path [-] (3.4,3.2) edge [right=40](7.1,8.9);
\end{scope}
\begin{scope}[>={Stealth[black]},
              every edge/.style={draw=green,very thick}]
\path [-] (4.4,4.6) edge [right=40](8.1,10.4);
\end{scope}

\begin{scope}[>={Stealth[black]},
              every edge/.style={draw=red,very thick}]
\path [-] (5,1) edge [right=0](5,1);
 
 \path [-] (1.5, 3.2) edge [right=40](-0.3,6.3);
 

\end{scope}
\begin{scope}[>={Stealth[black]},
  every edge/.style={draw=pink,very thick}]
 
  \path [-] (2.4, 4.6) edge [right=40](0.6,7.6);
 

\end{scope}

\put(-30,100){\line(2,3){240}}
\put(90,100){\line(2,3){30}}
\put(150,100){\line(2,3){30}}

\put(390,100){\line(-2,3){240}}
\put(420,145){\line(-2,3){240}}
\put(210,100){\line(-2,3){30}}
\put(150,100){\line(-2,3){30}}

\put(-30,100){\line(1,0){420}}
\put(-60,145){\line(1,0){480}}
\put(-30,190){\line(1,0){420}}
\put(-60,145){\line(2,3){240}}

\put(90,100){\circle*{4}}
\put(150,100){\circle*{4}}
\put(210,100){\circle*{4}}



\put(-40,90){$A=A_0$}
\put(420,150){$B_1$}
\put(-10,150){$A_1$}
\put(360,150){$B_1'$}
\put(380,90){$B=B_0$}
\put(170,410){$C=C_0$}
\put(400,190){$B_2=B_2'$}
\put(-70,200){$A_2=A_2'$}
\put(-50,150){$A_1'$}
\put(215,460){$C_1$}
\put(190,500){$C_2=C_2'$}
\put(130,460){$C_1'$}


\put(85,90){\tiny{${Z}_{1}$}}
\put(143,90){\tiny{${Z}_{2}$}}
\put(200,90){\tiny{${Z}_{3}$}}
\put(120,170){$\gamma_2'$}
\put(280,310){$\gamma_1'$}
\put(80,130){$\gamma_2$}
\put(155,130){$\gamma_1$}

\put(175,135){\tiny{${a}_{2}$}}
\put(115,135){\tiny{${a}_{1}$}}
\end{tikzpicture}
\end{picture}
\end{center}
\caption{ \em \small{ $\bigtriangleup AB C$ is moved $l^{-1}$ in the right and up along $\vec{e}_1$ to obtain $\bigtriangleup A_1 B_1 C_1$.
$\bigtriangleup A_1B_1 C_1$ is moved $\l^{-1}$ in the left and up along $\vec{e}_2$  to obtain $\bigtriangleup A_2 B_2C_2$. 
$\bigtriangleup AB C$ is moved  $l^{-1}$ in the left and up along $\vec{e}_2$ to obtain $\bigtriangleup A_1' B_1'C_1'$. $\bigtriangleup A_1'B_1' C_1'$ is moved  $l^{-1}$ in the right and up along $\vec{e}_1$ to obtain $\bigtriangleup A_2' B_2'C_2'=\bigtriangleup A_2 B_2C_2$. 
 After shifting  $\bigtriangleup A_1B_1 C_1$ down to  $\bigtriangleup ABC$,
the green path $\gamma_1'$ (the open $\bullet$-path in ${\cal H}(D_1,a_1, a_2)$) and the pink path $\gamma_1'$ (the closed $\circ$-path in ${\cal H}(a_1, a_2)$)  become to be the blue path $\gamma_1$
( the open $\bullet$-path in ${\cal H}(D,Z_1, Z_2)$) and the red path $\gamma_2$ (the closed $\circ$-path in ${\cal H}(D, Z_1, Z_2)$).}}
\end{figure}

 For any $a\in \overline{A_iB_i}$, let ${\cal E}(D_i, a)$ be the event that there is open path separating $a$ from $\overline{A_iC_i}$ in $D_i$.
 Similarly, for any $a\in \overline{A_i'B_i'}$, let ${\cal E}(D_i',a)$ be the event that there is open path separating $a$ from $\overline{A_i'C_i'}$ in $D_i'$.  For each $a\in \overline{A_iB_i}$, we denote by
 $$f_l(a)= {\bf P}_l( {\cal E}(D_i, a)).$$
 For each 
 $a\in \overline{A_i'B_i'}$, we denote by
 $$f_l'(a)= {\bf P}_l( {\cal E}(D_i', a)).$$
 By (3.1),
 $$f_l(a)=f_l'(a) \mbox{ if }a\in {A_i B_i} \mbox{ for even }i.\eqno{(3.2)}$$
 For any two $\circ$-adjacent vertices  $a_1< a_2$ in $\overline{A_iB_i}$, let ${\cal H}(D_i, a_1, a_2)$  be the event that there are two arm paths: one open $\bullet$-path from 
$a_2$ to $\overline {{ B_iC_i}}$ and one closed $\circ$-path from $a_1$ to $ \overline{A_iC_i}$ in $D_i$ (see Fig. 5).
Similarly, for any two $\circ$-adjacent vertices  $a_1< a_2$ in $\overline{A_i'B_i'}$, let ${\cal H}(D_i', a_1, a_2)$  be the event that there are two arm paths: one open $\bullet$-path from 
$a_2$ to $\overline {{ B_i'C_i'}}$ and one closed $\circ$-path from $a_1$ to $ \overline{A_i'C_i'}$ in $D_i'$. 
In fact, the events ${\cal H}(D, a_1, a_2)$ and ${\cal H}(D_i, a_1, a_2)$, and  ${\cal H}(D, a_1, a_2)$ and ${\cal H}(D_i', a_1, a_2)$ are the same events by shifting
$\bigtriangleup ABC$ to $\bigtriangleup A_{i} B_{i} C_{i}$, and to $\bigtriangleup A_i'B_i'C_i'$, respectively.
 For $a_1, a_2$ defined above, if $a_1, a_2\in \overline{A_i B_i}\cap \overline{A_i'B_i'}$ are $\circ$-adjacent for an even $i$,
 we  denote by  
$$g_l( a_1)=h_l(a_1)={\bf P}_l( {\cal H}(D_i, a_1, a_2))={\bf P}_l( {\cal H}(D_i', a_1, a_2)).$$
 If $a_1, a_2\in \overline{A_i B_i}\cap \overline{A_i'B_i'}$ are $\circ$-adjacent for  an odd $i$, we denote by
 $$ g_l( a_1)={\bf P}_l( {\cal H}(D_i, a_1, a_2))\mbox{ and } h_l( a_1)={\bf P}_l( {\cal H}(D_1', a_1, a_2)).$$
Thus, by (3.1), $g_l$ and $h_l$ are well defined on all the lattice points of $\bigtriangleup ABC$. If $a$ is a lattice point in $\bigtriangleup ABC$
with $a\in \overline{A_iB_i}$  for even $i$, then
 $$ g_l(a)=h_l(a).\eqno{(3.3)}$$
 If $a$ and $b$ for $ a=(a_1, a_2), b=(b_1, b_2)$ are lattice points  in $\bigtriangleup ABC$ with $b_2= a_2+ 2\sqrt{3}l^{-1}, a_1=b_1$,  then
 it follows from our construction that
 $$g_l(a)=g_l(b)\mbox{ and } h_l(a)=h_l(b).\eqno{(3.4)}$$
We select  $\circ$-adjacent lattice  points $ Z_1< Z_2< Z_3 \in \overline{ AB}$.
We  also select $\circ$-adjacent lattice  points $ a_1< a_2\in\overline{A_{1}B_{1}}$ such that   $\bigtriangleup Z_1a_1 Z_2$, $\bigtriangleup Z_2 a_2Z_3$  are two equilateral triangles with side length $l^{-1}$ (see Fig. 5). By using symmetry,  we show the following lemma.\\

{\bf Lemma 3.1.} {\em   For each $i \geq 0$ and for  $Z_1, Z_2, Z_3, a_1, a_2$ defined above,  
 $$g_l( a_1)={\bf P}_l( {\cal H}(D_{1},  a_1, a_2))=   g_l(Z_1)={\bf P}_l( {\cal H}(D,  Z_1, Z_2)),\eqno{(3.5)}$$
 $$h_l(a_1)={\bf P}_l( {\cal H}(D_{1}',  a_1, a_2))=   h_l( Z_2)={\bf P}_l( {\cal H}(D,  Z_2, Z_3)),\eqno{(3.6)}$$
\begin{eqnarray*}
 &&g_l(Z_2)-g_l(a_1)\\
 &= &{\bf P}_l( {\cal H}(D, Z_2, Z_3),{\cal D}^+(D, \tau^+(D), Z_3))-{\bf P}_l( {\cal H}(\tau^+(D),  Z_2, Z_3),{\cal D}^-(D, \tau^+(D), Z_2)), \hskip 1.5cm (3.7)
 \end{eqnarray*}
\begin{eqnarray*}
 &&h_l( a_1)-h_l(Z_1)\\
 &= &
  {\bf P}_l( {\cal H}(\tau^-(D),  Z_1, Z_2),{\cal D}^+(\tau^-(D), D, Z_2))-{\bf P}_l( {\cal H}(D,  Z_1, Z_2),{\cal D}^-(\tau^{-1}(D), D, Z_1)).\hskip 1.5 cm {(3.8)}
  \end{eqnarray*}}
  
 {\bf Proof.}   
 We show the equations of (3.5) in Lemma 3.1.
 We shift $D_{1}$ down along $\vec{e}_1$ to  be $D$ (see Fig. 5).
 Thus, $a_1$ and $ a_2$ are shifted to be $Z_1$ and $Z_2$.
  By translation invariance,
  $$g_l(a_1)={\bf P}_l( {\cal H}(D_1,  a_1, a_2))={\bf P}_l( {\cal H}(D,  Z_1, Z_2))=g_l(Z_1).\eqno{(3.9)}$$
 Thus, the  equations  of Lemma 3.1 in (3.5) follows.
 The equations of Lemma 3.1 in (3.6) can be showed by the same way if we shift $D_{1}'$ down along $\vec{e}_2$ to be $D$.
 
 Now we focus on the equation in (3.7). We shift $D_{1}$ down along $\vec{e}_2$ to  be $\tau^+(D)$ (see Fig. 6).
 Thus, $a_1$ and $ a_2$ are shifted to be $Z_2$ and $Z_3$.
 After shifting, the open $\bullet$-path and closed $\circ$-path in ${\cal H}(D_{1},  a_1, a_2)$ will   be  two new paths $\gamma_1'$ and $\gamma_2'$ from 
 $Z_3$ to $\tau^+(\overline{BC })$ and from $Z_2$ to $\tau^+( \overline{ AC})$, respectively.
 We denote this event  with these two new paths by $\hat {\cal H}(\tau^+(D),  Z_2, Z_3)$.
We have
 \begin{eqnarray*}
 &&g_l(Z_2)-g_l(a_1)=  {\bf P}_l( {\cal H}(D,  Z_2, Z_3))-{\bf P}_l( \hat {\cal H}(\tau^+(D),  Z_2, Z_3))\\
 &=& {\bf P}_l( {\cal H}(D,  Z_2, Z_3)\cap (\hat{\cal H}(\tau^+(D),  Z_2, Z_3))^C)-{\bf P}_l( \hat {\cal H}(\tau^+(D),  Z_2, Z_3)\cap ({\cal H}(D,  Z_2, Z_3))^C).\hskip .5cm {(3.10)}
 \end{eqnarray*}
 On ${\cal H}(D_{i},  Z_2, Z_3)\cap (\hat{\cal H}(\tau^+(D),  Z_2, Z_3))^C$, the open $\bullet$-path $\gamma_1$ in ${\cal H}(D,  Z_2, Z_3)$ from $Z_3$ to $\overline{BC}$ in $D$
 cannot go further to $\tau^+(\overline{BC})$ (see Fig. 6). By the same proof in (2.3.9),
 $${\bf P}_l( {\cal H}(D,  Z_2, Z_3)\cap (\hat{\cal H}(\tau^+(D),  Z_2, Z_3))^C)= {\bf P}_l( {\cal H}(D,  Z_2, Z_3),{\cal D}^+(D, \tau^+(D), Z_3)).\eqno{(3.11)}$$
 On $ \hat {\cal H}(\tau^+(D),  Z_2, Z_3)\cap ({\cal H}(D,  Z_2, Z_3))^C$, the closed $\circ$-path $\gamma_2'$ from $Z_2$ in $\hat {\cal H}(\tau^+(D),  Z_2, Z_3)$ to $\tau^+( \overline{ AC})$ cannot reach to $ \overline{ AC}$ (see Fig. 6). By the same proof of (2.3.3),
 $${\bf P}_l( \hat {\cal H}(\tau^+(D),  Z_2, Z_3)\cap ({\cal H}(D,  Z_2, Z_3))^C)= {\bf P}_l( {\cal H}(\tau^+(D),  Z_2, Z_3),{\cal D}^-(D, \tau^+(D), Z_2)).\eqno{(3.12)}$$
 The  equation in Lemma 3.1 of (3.7)  follows from (3.10)-(3.12). The equation in (3.8) can be showed by the same way if we shift $D_{1}'$ along $\vec{e}_1$ to be $\tau^{-}(D)$.
 $\blacksquare$\\
 
In Lemma 3.1, we investigate  the diagonal shifting for $g_l$ and $h_l$ in the  first layer. We also need to investigate the diagonal shifting for  $g_l$ and $h_l$ in the second layer.  We select  $\circ$-adjacent lattice  points $ Z_1< Z_2< Z_3 \in \overline{A_1B_1}$.
We  also select $\circ$-adjacent lattice  points $ a_1< a_2\in\overline{A_{2}B_{2}}$ such that   $\bigtriangleup Z_1a_1 Z_2$, $\bigtriangleup Z_2 a_2Z_3$ 
 are two equilateral triangles with side length $l^{-1}$.  After the zag-zig change from first layer to the second layer,  
we can use   the same proof of Lemma 3.1 by exchanging the  situations of $g_l$ and $h_l$ to show
 $$g_l( a_1)={\bf P}_l( {\cal H}(D_{2},  a_1, a_2))=   g_l(Z_2)={\bf P}_l( {\cal H}(D_1,  Z_2, Z_3)),\eqno{(3.13)}$$
 $$h_l(a_1)={\bf P}_l( {\cal H}(D_{2},  a_1, a_2))=   h_l(Z_1)={\bf P}_l( {\cal H}(D_1',  Z_1, Z_2)),\eqno{(3.14)}$$
\begin{eqnarray*}
 &&h_l(Z_2)-h_l(a_1)\\
 &= &{\bf P}_l( {\cal H}(D_1', Z_2, Z_3),{\cal D}^+(D_1', \tau^+(D_1'), Z_3))-{\bf P}_l( {\cal H}(\tau^+(D_1'),  Z_2, Z_3),{\cal D}^-(D_1', \tau^+(D_1'), Z_2)), \hskip 0.5cm (3.15)
 \end{eqnarray*}
\begin{eqnarray*}
 &&g_l( a_1)-g_l(Z_1)\\
 &= &
  {\bf P}_l( {\cal H}(\tau^-(D_1),  Z_1, Z_2),{\cal D}^+(\tau^-(D_1), D_1, Z_2))-{\bf P}_l( {\cal H}(D_1,  Z_1, Z_2),{\cal D}^-(\tau^{-1}(D_1), D_1, Z_1)).\hskip .5 cm {(3.16)}
  \end{eqnarray*}

Now we work on discrete line integrals.   
We select $A^\iota$ and $B^\iota$ in $\overline{ A_{\iota l} B_{\iota l}}$ for a small $\iota >0$, where
$A^\iota$ and $B^\iota$ are two lattice points by moving $A_\iota$ and $B_\iota$ up  $\iota$ along $\vec{e}_1$ and $\vec{e}_2$, respectively. 
\begin{figure}
\begin{center}
\begin{picture}(250,220)(40,70)
\setlength{\unitlength}{0.0125in}%
\begin{tikzpicture}
\thicklines
\begin{scope}[>={Stealth[black]},
              every edge/.style={draw=blue,very thick}]
\path [-] (5.3,3.2) edge [right=40](8.1,7.4);
\path [-] (1.5,3.2) edge [right=40](1.5,8,5);
\path [-] (-.5,6) edge [right=40](2.8,11);
\path [-] (2.8,11) edge [right=40](7,12);
\end{scope}
\begin{scope}[>={Stealth[black]},
              every edge/.style={draw=red,very thick}]
 \path [-] (5,1) edge [right=0](5,1);
 
 \path [-] (3.4, 3.2) edge [right=40](1.4,6);
 \path [-] (8,3.2) edge [right=40](10.5,6.85);
  \path [-] (10.4, 6.8) edge [right=40](8.8,9.2);
  \path [-] (8.8,9.2) edge [right=40](6,5.4);
   \path [-] (6,5.4) edge [right=40](3,8.5);
 

\end{scope}
\begin{scope}[>={Stealth[black]},
              every edge/.style={draw=green,very thick}]
\path [-] (-0.53,3.15) edge [right=40](6.2,13.2);
\path [-] (6.2,13.2) edge [right=40](13, 3.1);
\path [-] (-0.53,3.15) edge [right=40](13,3.1);

\end{scope}

\put(-30,100){\line(2,3){240}}
\put(90,100){\line(2,3){30}}
\put(150,100){\line(2,3){30}}
\put(210,100){\line(2,3){30}}

\put(390,100){\line(-2,3){210}}
\put(420,145){\line(-2,3){210}}
\put(270,100){\line(-2,3){30}}
\put(210,100){\line(-2,3){30}}
\put(150,100){\line(-2,3){30}}

\put(-30,100){\line(1,0){420}}
\put(0,145){\line(1,0){420}}

\put(90,100){\circle*{4}}
\put(150,100){\circle*{4}}
\put(210,100){\circle*{4}}
\put(270,100){\circle*{4}}



\put(-40,90){$A=A_0$}
\put(420,150){$B_1$}
\put(-10,150){$A_1$}
\put(380,90){$B=B_0$}
\put(170,410){$C$}
\put(-50,150){$\alpha_1$}
\put(215,460){$C_1$}


\put(85,90){\tiny{${Z}_{1}$}}
\put(143,90){\tiny{${Z}_{2}$}}
\put(200,90){\tiny{${Z}_{3}$}}
\put(260,90){\tiny{${Z}_{4}$}}
\put(120,170){$\gamma_2'$}
\put(280,230){$\gamma_1$}

\put(175,135){\tiny{${a}_{2}$}}
\put(115,135){\tiny{${a}_{1}$}}
\put(235,135){\tiny{${a}_{3}$}}
\end{tikzpicture}
\end{picture}
\end{center}
\caption{ \em \small{$D_1= \bigtriangleup A_1B_1C_1$ is shifted down along $\vec{e}_2$ to  be $\tau^+(\bigtriangleup ABC)$,
the green triangle.
 Thus, $a_1$ and $ a_2$ are shifted to be $Z_2$ and $Z_3$.
$\gamma_1$ (blue path), in event ${\cal H}(\bigtriangleup ABC, Z_2, Z_3)$ cannot reach to $\tau^+(\overline{BC})$, so it is  blocked by a closed $\circ$-path (red path). $\gamma_2'$ (red path)  in event ${\cal H}(\bigtriangleup A_1B_1C_1, a_1, a_2)$  is moved to be a closed $\circ$-path
from $Z_2$ to $\tau^+(\overline{ AC})$. $\gamma_2'$ cannot reach to $\overline{AC}$, so it is blocked by an open $\bullet$-path.}}
\end{figure}
 We denote by  $T= T_\iota$ the isosceles trapezium consisting of four vertices  
$A_\iota, B_\iota, B^\iota, A^\iota$ for a small $\iota$.
For any $\bullet$-path $L=\{Z_0=(x_0, y_0), \cdots, Z_k=(x_k, y_k), \cdots, Z_n=(x_n, y_n)\}$ from $Z_0$ to $Z_n$ in $\bigtriangleup ABC$, if $F_l$ is a discrete function defined on $T$, we define its discrete line integral on $L$  by
$$ \int_L F_l(Z) dx= \sum_{k=0}^{n-1} F_l(Z_k)(x_{k+1}- x_k)\mbox{ and }  \int_L F_l(Z) dy= \sum_{k=0}^{n-1} F_l(Z_k)(y_{k+1}- y_k).\eqno{(3.17)}$$ 
For any $\bullet$-adjacent circuit, called a discrete {\em contour}, it can be viewed a path enclosed with $l^{-1}$-parallelograms. We show the following lemma.\\

{\bf Lemma 3.2.} {\em For  a discrete contour $C\subset T$ enclosed by at most the $Ml$ many edges of  $l^{-1}$-parallelograms for a constant $M$, then there exist $c_1=c_1(M)$ and  $\alpha>0$  such that
$$\left | \oint_C( \sqrt{3} l g_l(Z) dx+ l h_l(Z) dy)\right |\leq  c_1 l^{-\alpha}.$$}

{\bf Proof.} 
We consider a discrete contour $C\subset T$ enclosed by at most the $Ml^{}$ many bonds of  $l^{-1}$-parallelograms. We first focus on an $l^{-1}$-parallelogram  $P$ in the first  layer. More precisely, $P$ has   its vertices $Z_1, \alpha_1, Z_2, a_1$ in $T$ with $\alpha_1\in \overline{AB}$,
$Z_1< Z_2\in \overline{A_1B_1}$ and $a_1\in \overline{A_2B_2}$. In addition, we denote by $Z_3 > Z_2$, $\alpha_2 > \alpha_1$ and $a_2 > a_1$ such that $Z_3$ $\circ$-adjacent to $Z_2$, $\alpha_2$ $\circ$-adjacent to $\alpha_1$ and
$a_2 $ $\circ$-adjacent to $a_1$ in $\overline{A_1B_1}$, $\overline{A B}$ and $\overline{A_2B_2}$, respectively.  By the definition in (3.17),
$$\oint_P l g_l(Z) dx=1/2 (g_l(\alpha_1)-g_l(Z_2)- g_l(a_1)+ g_l(Z_1)).\eqno{(3.18)}$$
By (3.18) and (3.5) in Lemma 3.1,
$$\oint_P l g_l(Z) dx=1/2 (- g_l(a_1)+ g_l(Z_1)).\eqno{(3.19)}$$
By (3.16) and (3.19) that
\begin{eqnarray*}
 &&\oint_P l g_l(Z) dx=-1/2(g_l( a_1)-g_l(Z_1))\\
 &= &
  {-1\over 2}[{\bf P}_l( {\cal H}(\tau^-(D_1),  Z_1, Z_2),{\cal D}^+(\tau^-(D_1), D_1, Z_2))-{\bf P}_l( {\cal H}(D_1,  Z_1, Z_2),{\cal D}^-(\tau^{-1}(D_1), D_1, Z_1))].\hskip .1 cm {(3.20)}
  \end{eqnarray*}

Now we work on the contour integral of $P$ for $h_l$ in the first layer. For the above $Z_1, \alpha_1, Z_2, a_1$,
$$\oint_P lh_l(Z) dy=\sqrt{3}/2( h_l(\alpha_1)+h_l(Z_2)- h_l(a_1)- h_l(Z_1))=\sqrt{3}/2( h_l(\alpha_1)-h_l(Z_1)+h_l(Z_2)- h_l(a_1)).\eqno{(3.21)}$$
It follows from (3.6) in Lemma 3.1 that
$$\oint_P lh_l(Z) dy=\sqrt{3}/2( h_l(Z_2)- h_l(a_1)).\eqno{(3.22)}$$
By (3.15) and (3.22),
\begin{eqnarray*}
 &&\oint_P lh_l(Z) dy=\sqrt{3}/2(h_l(Z_2)-h_l(a_1))\\
 &= &\sqrt{3}/2[{\bf P}_l( {\cal H}(D_1', Z_2, Z_3),{\cal D}^+(D_1', \tau^+(D_1'), Z_3))-{\bf P}_l( {\cal H}(\tau^+(D_1'),  Z_2, Z_3),{\cal D}^-(D_1', \tau^+(D_1'), Z_2))].\hskip .1cm (3.23)
 \end{eqnarray*}
Together with (3.20) and (3.23),
\begin{eqnarray*}
&&\oint_P \sqrt{3}lg_l(Z) dx+lh_l(Z) dy\\
&&={\sqrt{3}\over 2} \{-[{\bf P}_l( {\cal H}(\tau^-(D_1),  Z_1, Z_2),{\cal D}^+(\tau^-(D_1), D_1, Z_2))-{\bf P}_l( {\cal H}(D_1,  Z_1, Z_2),{\cal D}^-(\tau^{-1}(D_1), D_1, Z_1))]\\
&&+[{\bf P}_l( {\cal H}(D_1', Z_2, Z_3),{\cal D}^+(D_1', \tau^+(D_1'), Z_3))-{\bf P}_l( {\cal H}(\tau^+(D_1'),  Z_2, Z_3),{\cal D}^-(D_1', \tau^+(D_1'), Z_2))]\} \\
&&={\sqrt{3}\over 2} \{[-{\bf P}_l( {\cal H}(\tau^-(D_1),  Z_1, Z_2),{\cal D}^+(\tau^-(D_1), D_1, Z_2))+{\bf P}_l( {\cal H}(D_1', Z_2, Z_3),{\cal D}^+(D_1', \tau^+(D_1'), Z_3))]\\
&&+[{\bf P}_l( {\cal H}(D_1,  Z_1, Z_2),{\cal D}^-(\tau^{-1}(D_1), D_1, Z_1))-{\bf P}_l( {\cal H}(\tau^+(D_1'),  Z_2, Z_3),{\cal D}^-(D_1', \tau^+(D_1'), Z_2))]\}
\hskip .5cm (3.24)\\
\end{eqnarray*}
Note that $\tau^{-1}(D_1)=D_1'$ and $\tau^+(D_1')=D_1$, so by (3.24),
\begin{eqnarray*}
&&\oint_P \sqrt{3}lg_l(Z) dx+lh_l(Z) dy\\
&&={\sqrt{3}\over 2} \{[-{\bf P}_l( {\cal H}(D_1',  Z_1, Z_2),{\cal D}^+(D_1', D_1, Z_2))+{\bf P}_l( {\cal H}(D_1', Z_2, Z_3),{\cal D}^+(D_1', D_1, Z_3))]\\
&&+[{\bf P}_l( {\cal H}(D_1,  Z_1, Z_2),{\cal D}^-(D_1', D_1, Z_1))-{\bf P}_l( {\cal H}(D_1,  Z_2, Z_3),{\cal D}^-(D_1', D_1, Z_2))]\}
\end{eqnarray*}
By Lemma 2.4.1, for the above  $Z_1, \alpha_1, Z_2, a_1$ in the first  layer,
there exist $c_1$ and  $\alpha>0$  such that
$$\left | \oint_P( \sqrt{3} l g_l(Z) dx+ l h_l(Z) dy)\right |\leq  c_1 l^{-2-\alpha}.\eqno{(3.25)}$$

We focus on an $l^{-1}$-parallelogram  $P$ in the second layer. More precisely, $P$ has   its vertices $Z_1, \alpha_1, Z_2, a_1$ in $T$ with $\alpha_1\in \overline{A_1B_1}$,
$Z_1< Z_2\in \overline{A_2B_2}$ and $a_1\in \overline{A_3B_3}$. In addition, we denote by $Z_3 > Z_2$, $\alpha_2 > \alpha_1$ and $a_2 > a_1$ such that $Z_3$ $\circ$-adjacent to $Z_2$, $\alpha_2$ $\circ$-adjacent to $\alpha_1$ and
$a_2 $ $\circ$-adjacent to $a_1$ in $\overline{A_2B_2}$, $\overline{A_1 B_1}$ and $\overline{A_3B_3}$, respectively.
By the definition in (3.17),
$$\oint_P l g_l(Z) dx=1/2 (g_l(\alpha_1)-g_l(Z_2)- g_l(a_1)+ g_l(Z_1)).\eqno{(3.26)}$$
Note that, by (3.4), the values of $g_l$ in the first and the third layers are the same, so by (3.5) in Lemma 3.1,
$$\oint_P l g_l(Z) dx=1/2 (g_l(\alpha_1)-g_l(Z_2)).\eqno{(3.27)}$$
By (3.27), if we replace $Z_1$ by $\alpha_1$ and $a_1$ by $Z_2$ in (3.16), then
\begin{eqnarray*}
 &&\oint_P l g_l(Z) dx=1/2(g_l(\alpha_1)-g_l(Z_2))=-1/2(g_l(Z_2)-g_l(\alpha_1))\\
 &= &
  {-1\over 2}[{\bf P}_l( {\cal H}(\tau^-(D_1),  \alpha_1, \alpha_2),{\cal D}^+(\tau^-(D_1), D_1, \alpha_2))-{\bf P}_l( {\cal H}(D_1,  \alpha_1, \alpha_2),{\cal D}^-(\tau^{-1}(D_1), D_1, \alpha_1))].\hskip .1 cm {(3.28)}
  \end{eqnarray*}

Now we work on the contour integral of $P$ for $h_l$ in the second  layer. For the above $Z_1, \alpha_1, Z_2, a_1$,
$$\oint_P lh_l(Z) dy=\sqrt{3}/2( h_l(\alpha_1)+h_l(Z_2)- h_l(a_1)- h_l(Z_1))=\sqrt{3}/2( h_l(\alpha_1)-h_l(Z_1)+h_l(Z_2)- h_l(a_1)).\eqno{(3.29)}$$
Note that, by (3.4), the values of  $h_l$ in the first and the third layers are the same,   so  by (3.6) in Lemma 3.1,
$$\oint_P lh_l(Z) dy=\sqrt{3}/2( h_l(\alpha_1)-h_l(Z_1) ).\eqno{(3.30)}$$
If we replace  $Z_2$ by $\alpha_1$ and $a_1$ by $Z_1$ in (3.15), then
\begin{eqnarray*}
 &&h_l( \alpha_1)-h_l(Z_1)\\
 &= &
  {\bf P}_l( {\cal H}(D_1',  \alpha_1, \alpha_2),{\cal D}^+(D_1', \tau^+(D_1'), \alpha_2))-{\bf P}_l( {\cal H}(\tau^+(D_1'),  \alpha_1, \alpha_2),{\cal D}^-(D_1', \tau^+(D_1'), \alpha_1)).\hskip .5 cm {(3.31)}
  \end{eqnarray*}
Thus, 
\begin{eqnarray*}
 &&\oint_P lh_l(Z) dy=\sqrt{3}/2[h_l(\alpha_1)+h_l(Z_1))]\\
 &= &\sqrt{3}/2[{\bf P}_l( {\cal H}(D_1',  \alpha_1, \alpha_2),{\cal D}^+(D_1', \tau^+(D_1'), \alpha_2))-{\bf P}_l( {\cal H}(\tau^+(D_1'),  \alpha_1, \alpha_2),{\cal D}^-(D_1', \tau^+(D_1'), \alpha_1))], \hskip .1cm (3.32)
 \end{eqnarray*}
Together with (3.28) and (3.32),
\begin{eqnarray*}
&&\oint_P \sqrt{3}lg_l(Z) dx+lh_l(Z) dy\\
&&={\sqrt{3}\over 2} \{-[{\bf P}_l( {\cal H}(\tau^-(D_1),  \alpha_1, \alpha_2),{\cal D}^+(\tau^-(D_1), D_1, \alpha_2))-{\bf P}_l( {\cal H}(D_1,  \alpha_1, \alpha_2),{\cal D}^-(\tau^{-1}(D_1), D_1, \alpha_1))]\\
&&+[{\bf P}_l( {\cal H}(D_1',  \alpha_1, \alpha_2),{\cal D}^+(D_1', \tau^+(D_1'), \alpha_2))-{\bf P}_l( {\cal H}(\tau^+(D_1'),  \alpha_1, \alpha_2),{\cal D}^-(D_1', \tau^+(D_1'), \alpha_1))]
\}. \hskip .5cm (3.33)\\
\end{eqnarray*}
Note that $\tau^{-1}(D_1)=D_1'$ and $\tau^+(D_1')=D_1$, so by (3.33) for the $P$ in the second layer,
\begin{eqnarray*}
&&\oint_P \sqrt{3}lg_l(Z) dx+lh_l(Z) dy\\
&&={\sqrt{3}\over 2} \{-[{\bf P}_l( {\cal H}(D_1',  \alpha_1, \alpha_2),{\cal D}^+(D_1',D_1, \alpha_2))-{\bf P}_l( {\cal H}(D_1,  \alpha_1, \alpha_2),{\cal D}^-(D_1', D_1, \alpha_1))]\\
 &&+[{\bf P}_l( {\cal H}(D_1',  \alpha_1, \alpha_2),{\cal D}^+(D_1', D_1, \alpha_2))-{\bf P}_l( {\cal H}(D_1,  \alpha_1, \alpha_2),{\cal D}^-(D_1', D_1, \alpha_1))]\}=0. \hskip 1.5cm (3.34)\\
\end{eqnarray*}
By (3.4), (3.33) and (3.34) together with translation invariance, for any $l^{-1}$-parallelogram $P\subset T$, there exist $c_1$ and  $\alpha>0$  such that
$$\left | \oint_P( \sqrt{3} l g_l(Z) dx+ l h_l(Z) dy)\right |\leq  c_1 l^{-2-\alpha}.\eqno{(3.35)}$$

Now we focus on a discrete contour. Suppose that the boundary of the contour $C$ contains  less than $Ml^{-1}$ edges of $l^{-1}$-parallelograms .  Let $I$ be the $l^{-1}$-parallelograms enclosed by $C$.  If the boundary of the contour $C$ is less than $Ml^{-1}$ edges of $l^{-1}$-parallelograms,  it is easy to verify that the number of $I$  is at most $c_2 l^{2}$  for a constant $c_2=c_2(M)$.  Thus, by (3.35),
$$\left | \oint_C( \sqrt{3} l g_l(Z) dx+ l h_l(Z) dy)\right |\leq \left |\sum_{P\in I} \oint_P( \sqrt{3} l g_l(Z) dx+ l h_l(Z) dy)\right |\leq c_2 l^{2} c_1l^{-2-\alpha}\leq  c_2c_1 l^{-\alpha}.\eqno{(3.36)}$$
Lemma 3.2 follows (3.35). $\blacksquare$\\

\section{ The convergence  on the separating probabilities.}
The functions discussed above are only defined on the lattice points in $T$. We need to extend them to all the points on $T$ by linear interpolation.
 We consider a discrete  function $\phi_\delta(Z)$ defined on the  triangular lattice $\bf T_\delta$.
For each triangle $\omega$ of $\bf T_\delta$ with three vertices $\alpha_i=(a_i, b_i)$ for $i=1,2,3$,  we construct the plane by using the three points 
$(a_i, b_i, \phi(\alpha_i))$ for $i=1,2,3$.  The piece of the plane on the domain $\omega$ is called  a triangular {\em plaquette}. All these triangular  plaquettes consist of a surface extending discrete function values of $\phi_\delta$ with domain on  ${\bf N}_{l^{-1}}$ to $\bf R_2$. The extension is called a {\em linear extension} and the extension function on ${\bf R}_2$  is  denoted by $\bar \phi_\delta$.  By the construction, $\bar \phi_\delta (Z)$ is continuous. Note that our ${\bf N}_{l^{-1}}$ is a subset of the triangular lattice. We do the linear extensions for
$\bar f_{l}$,  $\bar G_{l}$ and $\bar H_{l}$ from  the lattice points of $T$ to  all the points of $T$, denoted by $\bar T$.
For each $Z_1\in \overline{A_iB_i}\cap T$, let $Z_2, Z_3$ also in $\overline{A_i B_i}\cap T$ such that
$Z_1< Z_2 < Z_3$ are $\circ$-adjacent. We denote by
$$S_l(Z_1)= {G_l(Z_2)-G_l(Z_1)\over l^{-1}}.$$
We also do the linear extension for $\bar S_l$.
Note that $f_l\leq 1$  on ${\bf N}_{l^{-1}}\cap \bigtriangleup ABC$, so $\bar f_l$ is uniformly bounded on all the lattice points of $\bigtriangleup ABC$. 
It follows from (2.3.1) that $G_l$ and $H_l$ are also uniformly bounded on ${\bf N}_{l^{-1}}\cap T$.
so $\bar G_l$ and $\bar H_l$ are uniformly bounded on  $\bar T$ by Lemma 4.1. 
By (2.3.12), for $Z\in T$,
$$|S_l(Z)| \mbox{ is uniformly bounded on all the lattice points in }T.\eqno{(4.1)}$$
By the linear extension for $f_l$, $G_l$, $H_l$,  and $S_l$, it is straightforward  to show that $\bar f_l$, $\bar G_l$, $\bar H_l$ and $\bar S_l$  are uniformly bounded on $\bar T$. 
We summarize the above observation as the following lemma.\\

{\bf Lemma 4.1.} {\em $\bar f_l$, $\bar G_l$, $\bar H_l$  and $\bar S_l$ are uniformly bounded on $\bar T$. }\\

By (2.3.1), (2.3.12), (2.4.5), and their linear extensions,  for any $Z$ in an equilateral $\bigtriangleup a_1 a_2 a_3$ with side length 
$l^{-1}$,  we have, for $i=1,2,3$,
\begin{eqnarray*}
&&|\bar f_l(Z)- f_l(a_i)|\leq O(l^{-1}),|\bar G_l(Z)- G_l(a_i)|\leq O(l^{-1}),\\
&&|\bar H_l(Z)- H_l(a_i)|\leq O(l^{-1}),|\bar S_l(Z)- S_l(a_i)|\leq O(l^{-1}).\hskip 6cm (4.2)
\end{eqnarray*}

 We will show that  these functions are equicontinuous as the following lemma.\\

{\bf Lemma 4.2}. {\em $\bar f_l$ is uniformly equicontinuous  $\bigtriangleup ABC$, and $\bar G_l$,  $\bar H_l$  and $\bar S_l$ are uniformly equicontinuous on $ \bar T$.}\\

{\bf Proof.}  It is well known (see Chapter 7, Claim 22 in  Bollobas, B. and Riordan, O. (2006)) that $\bar f_l$ is uniformly equicontinuous  $\bigtriangleup ABC$.
Now we show the other functions are also uniformly equicontinuous. Without loss of generality, we only show that $\bar S_l$ is uniformly equicontinuous on $ \bar T$, since  $\bar G_l$ and $\bar H_l$ are easier to handle by the same way.
 Given $\epsilon >0$ small, we must show that there is a $\delta$ such that for all $X, Y\in \bar T$ with $d(X, Y) < \delta$, we have
 $$|\bar S_l(X)-\bar S_l(Y)|< \epsilon.\eqno{(4.3)}$$
 If $\epsilon=  c_1l^{-\alpha}_0$ for the $\alpha$ and $c_1$ defined in Lemma 2.4.1, then by using the linear extension of $\bar S_l$ and Lemma 4.1,  we can take $\delta(\epsilon)$  such that for all $l \geq l_0$,
(4.3) holds.  

 For each $l <  l_0$, we assume that $n l^{-1} = l_0^{-1}$ for some $n$.
 If  $X,Y\in \bar T$ are not lattice points of ${\bf N}_{l^{-1}}$ between $\overline {A_i B_i}$ and $\overline{A_{i+2} B_{i+2}}$,  
 then we select $\circ$-adjacent lattice points $X_1, X_2$ and $Y_1, Y_2$  on  $\overline {A_i B_i}$
such that $X_1 <  X_2$ and $Y_1< Y_2$, and $X$ and $Y$ stay inside the equilateral  triangles with the left and the right vertices $X_1, X_2$ and $Y_1, Y_2$, respectively.  It follows from (4.2) that
$$|\bar S_l(X)- S_l(X_1)|\leq c_1 l^{-1}\leq c_1 l^{-\alpha} \leq \epsilon\mbox{ and } |S_l(Y_1)-\bar S_l(Y)|\leq c_1l^{-1} \leq c_1 l^{-\alpha}\leq \epsilon.\eqno{(4.4)}$$
In addition,  if we take with $d(X, Y)\leq \delta$ with $\delta = (\epsilon /c_1)^{1/\alpha}/2$, 
then  $d(X_1, Y_1) \leq l_0^{-1}=nl^{-1}$.
By Lemma 2.4.1,
$$|S_l(X_1)- S_l(Y_1)|\leq c_1 (n l^{-1})^{-\alpha} \leq c_1 l_0^{-\alpha} \leq  \epsilon.\eqno{(4.5)}$$
By (4.4)and (4.5),
$$|\bar S_l(Y)- \bar S_l(X)|\leq |\bar S(Y)- S(Y_1)|+ |\bar S(X)-S(X_1)|+ |S(Y_1)-S(X_1)|\leq 3\epsilon.\eqno{(4.6)}$$

Now we assume that $X=(x_1, x_2)$ and $Y=(y_1, y_2)$ with $x_2\neq y_2$ such that $X$ and $Y$ are not between $\overline {A_i B_i}$ and $\overline{A_{i+2} B_{i+2}}$, but $X$ is between $\overline {A_i B_i}$ and $\overline{A_{i+2} B_{i+2}}$. 
We suppose that $y_2 > x_1$ without loss of generality.
We suppose that $Z=(y_1, x_2)$ and $X$ are not in the same $l^{-1}$-parallelogram.
We go from $X$ horizontally to $Z$ for $Y$ and $Z$ not in a $l^{-1}$-parallelogram. 
There are three different $l^{-1}$-parallelograms containing $X$, $Y$ and $Z$. We pick the left vertices $X_1$, $Y_1$ and $ Z_1$ from these three  $l^{-1}$-parallelograms.
By Lemma 3.1 and Lemma 2.4.1, note that there is a zig-zag path from $Y_1$ to $Z_1$, so
$$|S_l(Y_1)- S_l(Z_1)|\leq c_1 l^{-\alpha}\leq \epsilon.\eqno{(4.7)}$$
On the other hand, by (4.2),
$$| S(X_1)-\bar S(X)|\leq \epsilon\mbox{ and } | S(Y_1)-\bar S(Y)|\leq \epsilon \mbox{ and }|\bar S(Z)-S(Z_1)|\leq \epsilon.\eqno{(4.8)}$$
Thus, by (4.6), (4.7) and (4.8), if we take with $d(X, Y)\leq \delta$ with $\delta = (\epsilon /c_1)^{1/\alpha}/2$, 
$$|\bar S_l(X)-\bar S_l(Y)|\leq |\bar S(X)-\bar S(Z)|+ |\bar S(Z)-S(Z_1)| + |S(Z_1)- S(Y_1)|+ S(Y_1)-\bar S(Y)|\leq 4\epsilon.\eqno{(4.9)}$$
If $X$ and $Z$ are in the same $l^{-1}$-parallelogram, we do not need to use $Z$, but just use $X$ and $Y$ in the estimates in (4.8) and (4.9)
to show (4.10). This showed that $\bar S_l$ is uniformly equicontinuous.
 Lemma 4.2 follows.
$\blacksquare$\\

 By Lemmas 4.1-4.2, and the Arzela-Ascoli theorem,   there exists a subsequence $\{l(m)\}$ of $\{l\}$ such that
 $\bar f_{l(m)}$ converges uniformly to $ f$  on $\bar T$. In addition,   there exists another subsequence $\{l(m(s))\}$ of $\{l(m)\}$ such that
 $\bar G_{l(m(s))}$ converges uniformly to $G$ on $\bar T$. In addition, there exists another subsequence $\{l(m(s(t)))\}$ of $\{l(m(s))\}$ such that
 $\bar H_{l(m(s(t)))}$ converges uniformly to $H$ on $\bar T$. Moreover, there exists another subsequence $\{l(m(s(t(q))))\}$ of $\{l(m(s(t)))\}$ such that
 $\bar S_{l(m(s(t(q))))}$ converges uniformly to $S$ on $\bar T$.
 By (3.3), we know that $G=H$ on $\bar T$. We summarize the above observations by the following proposition.\\
 
 {\bf Proposition 1.} {\em There exists a subsequence $\{l_k\}$ from $\{l\}$ such that
 $\bar f_{l_k} $ converges uniformly to a continuous function $f$ on $\bigtriangleup ABC$. In addition, $\bar G_{l_k}$ and $\bar H_{l_k}$ converge uniformly to a continuous function $G$ on $\bar T$. Moreover, $\bar S_{l_k}$ converges uniformly a continuous $S$ on $\bar T$.}\\

 By Proposition 1, we show the following lemma.\\
 
 {\bf Lemma 4.3} {\em For $Z=(x, y)\in \bar T$, 
 $${\partial f (Z)\over \partial x}=G(Z) \mbox{ and } {\partial^2 f(Z) \over \partial x^2}={\partial G(Z)\over \partial x}=S(Z).$$}
 
 {\bf Proof.}  
 We consider lattice points $Z=(z_1, z_2)\in T$ and $(z_1+h, z_2)\in T$. By Proposition 1, we take $h$ small such that for each $X$ with  $d(X, Z)\leq h$,
 $$G(X)- G(Z)=o(h).\eqno{(4.10)}$$
 We take $l_k$ large and divide interval $[z_1, z_1+h]$ into $ l_k^{-1}$ many smaller equal intervals with  a length $hl_k^{-1}$.
 We denote  these intervals  by $\{[Z_i, Z_{i+1}]\}$ for $i=1, \cdots, l_k-1$ for $Z_1=(z_1, z_2)$ and $Z_{l_k}=(z_1+h, z_2)$ and for lattice points $\{ Z_i\}$ in $T$.  By Proposition 1,  we take $l_k$ small such that
 $$| f(z_1+h, z_2)- f_{l_k}(z_1+h, z_2)|=o( h^2) \mbox{ and } | f(z_1, z_2)- f_{l_k}(z_1, z_2)|=o( h^2)\eqno{(4.11)}$$
 and  for each $i$,
 $$| G_{l_k} (Z_i)- G(Z_i)|=o( h).\eqno{(4.12)}$$
 By (4.10)-(4.12), 
 \begin{eqnarray*}
 &&{f(z_1+h, z_2)- f(z_1, z_2)\over h}
={ f_{l_k}(z_1+h, z_2)-  f_{l_k}(z_1, z_2)\over h}+o(h)\\
&=&{ \sum_{i=1}^{l_k}} l_k^{-1}  { f_{l_k}(Z_{i+1})-  f_{l_k}(Z_i) \over  hl^{-1}_k}+o(h)={ \sum_{i=1}^{l_k}}{ l_k^{-1}   G_{l_k}(Z_i)}+o(h) \\
&=  &   \sum_{i=1}^{l_k} l_k^{-1}  (G(Z_i)+o(h))+o(h)= \sum_{i=1}^{l_k} l_k^{-1}  (G(Z)+o(h))+o(h)=G(Z)+o(h) .\hskip 2.5cm {(4.13)}
\end{eqnarray*}
By (4.13), 
$$\lim_{h\rightarrow 0} {f(z_1+h, z_2)- f(z_1, z_2)\over h}= G(Z).\eqno{(4.14)}$$

If $Z=(z_1, z_2)$ and $(z_1+h, z_2)$ are not lattice points in $\bar T$,  then by (4.2), we take $Z'=(z_1', z_2')$ and $(z_1'+h, z_2)$ are lattice points with
$d(Z, Z') \leq l^{-1}_k$ for large $l_k$ such that
$$|\bar f_{l_k}(z_1+h, z_2)- f_{l_k} (z_1'+h, z_2')|=o(h^2) \mbox{ and } |\bar f_{l_k}(z_1, z_2)- f_{l_k} (z_1', z_2')|=o(h^2).\eqno{(4.15)}$$
By using the same estimate in (4.13) 
 $${\bar f(z_1+h, z_2)- \bar f(z_1, z_2)\over h}
={ f_{l_k}(z_1'+h, z_2')-  f_{l_k}(z_1', z_2')\over h}+o(h)=G(Z)+o(h) .\eqno{(4.16)}$$
 (4.14) still holds if $Z=(z_1, z_2)$ and $(z_1+h, z_2)$ are not lattice points.
 Thus, the first equation in Lemma 4.4 follows.  
 If we replace $f$ by $G$ and $G$ by $S$ in the proofs in (4.13)-(4.16), then we can show the second and third equations in Lemma 4.4.
Lemma 4.3 follows. $\blacksquare$\\

 By Proposition 1, we also show the following lemma.\\

{\bf Lemma 4.4.} {\em   For $Z=(x, y)\in \bar T$,
$ {\partial H (Z)\over \partial y}=0.$}\\

{\bf Proof.} For a smaller $h$ and $Z=(z_1, z_2)\in \bar T$, 
we take $l_k^{-1}$ small such that $ (z_1, z_2)\in w=\bigtriangleup \alpha_1\alpha_2\alpha_3$ and $(z_1, z_2+h)\in w'=\bigtriangleup \alpha_1'\alpha_2'\alpha_3' $ for two equilaterals with side length $l^{-1}_k$ such that $\alpha_i$ and $\alpha_i'$ for $i=1,2,3$ are lattice points with $\alpha_1$ and $\alpha_1'$ on a vertical line $L$.
By Proposition 1 and (4.2), we take $l_k$ large such that
$$|H(z_1, z_2)-\bar H_{l_k}(z_1, z_2)|\leq h^2\mbox{ and } |H(z_1, z_2+h)-\bar H_{l_k}(z_1, z_2+h)|\leq h^2\eqno{(4.17)}$$
and 
$$| H_{l_k}(\alpha_1)- \bar H_{l_k}(z_1, z_2)|\leq h^2\mbox{ and } | H_{l_k}(\alpha_1')- \bar H_{l_k}(z_1, z_2+h)|\leq h^2.\eqno{(4.18)}$$
For the small $l_k^{-1}$, by Lemma 4.2 and Proposition 1 and (3.4), there are  lattice points $\beta_1, \beta_2, \cdots, \beta_i, \cdots, \beta_m$ on $L$ between $\alpha_1$ and $\alpha_1'$ with
$$d(\beta_i, \beta_{i+1})=2\sqrt{3} l^{-1}_k \mbox{ and } d(\alpha, \beta_1)\leq l^{-1}_k \mbox{ and } d(\alpha' ,\beta_m)\leq l^{-1}_{k}\eqno{(4.19)}$$
such that  for each $i$, 
$$H_{l_k}(\beta_i)- H_{l_k}(\beta_{i+1})=0, |H_{l_k}(\alpha_1)- H_{l_k}(\beta_1)|\leq h^2\mbox{ and } |H_{l_k}(\alpha_1')- H_{l_k}(\beta_m)|\leq h^2.\eqno{(4.20)}$$
Thus, by (4.17))-(4.20),
\begin{eqnarray*}
&&{ H(z_1, z_2+h)- H(z_1, z_2)\over h}={\bar H_{l_k}(z_1, z_2+h)-\bar  H_{l_k}(z_1, z_2)\over h}+o(h)\\
&= &{\bar H_{l_k}(z_1, z_2+h)- H_{l_k}(\alpha_1') +  H_{l_k}( \alpha_1')- H_{l_k}(\beta_m)\over h} + {\sum_{i=1}^m  H_{i_k}(\beta_{i+1})- H_{l_k}(\beta_i)\over h}\\
&&+ { H_{l_k}(\beta_1)-  H_{l_k}( \alpha_1) \over h} +  { H_{l_k}( \alpha_1) -\bar H_{l_k}(z_1, z_2)\over h} +o(h) =o(h).\hskip 5cm (4.21)
\end{eqnarray*}
Thus,
$$\lim_{h\rightarrow 0} {H(z_1, z_2+h)- H(z_1, z_2)\over h}=0.\eqno{(4.22)}$$
Lemma 4.4 follow. $\blacksquare$\\




{\bf Proposition 2.} For $Z\in \bar T$,
$${\partial^2 f(Z)\over \partial  x^2}=0.$$

{\bf Proof.}  Since $G$ and $ H$ are  uniformly  continuous on $\bar T$,  for any smooth contour $C$ on $\bar T$, if we take $l_k$ large, there exists a discrete contour $C'$
included by at most $O(l_k)$ many edges of $l_k^{-1}$-parallelograms  in ${\bf N}_{l^{-1}_k}$ such that
 $$\oint_C( \sqrt{3} \bar G_{l_k}(Z) dx+ \bar H_{l_k} (Z) dy)= \oint_{C'}( \sqrt{3} G_{l_k}(Z) dx+ H_{l_k} (Z) dy)+o(l^{-1}_k).\eqno{(4.24)}$$
By Lemma 3.2 and Proposition 1,
$$\lim_{ k\rightarrow \infty} \oint_{C'}( \sqrt{3}  G_{l_k}(Z) dx+  H_{l_k} (Z) dy)=0\eqno{(4.25)}.$$
Note that $\bar G_{l_k}$ and $\bar H_{l_k}$ converges uniformly to $G $ and $H$ on $T$, so, by (4.24) and (4.25),   if we  take $k\rightarrow \infty$ in both sides of (4.24),  then
$$\oint_C( \sqrt{3} G(Z) dx+ H (Z)) dy=0.\eqno{(4.26)}$$
It follows from  (4.26) that we can set 
$$t(x, y)= \int^{(x, y)} _{(a, b)}( \sqrt{3} G(Z) dx+ H (Z)dy) \mbox{ for any $(a, b)$ and $(x, y) \in T$}.\eqno{(4.27)}$$
  Note that $\partial G/\partial x$ and $\partial H/\partial y$  are continuous on $T$, so
$${\partial^2 t(Z) \over \partial y \partial x}=\sqrt{3} {\partial G\over \partial x}={\partial^2 t(Z) \over \partial x \partial y}={\partial H\over \partial y}.\eqno{(4.28)}$$
By (4.28), Lemma 4.3 and Lemma 4.4,
$${\partial^2 f(Z)\over \partial  x^2}={\partial G\over \partial x}={1\over \sqrt{3}}{\partial H\over \partial y}=0.\eqno{(4.29)}$$
Proposition 2 follows from (4.29). $\blacksquare$\\

\section{ Proof of Theorem.} 
If we solve the equation in Proposition 2, then there exist fixed constant $a$  and $b$ such that
$$f(U)=au+b\mbox{ for }U=(u, 0)\in [A_\iota, B_\iota]\mbox{ for  any }\iota >0. \eqno{(5.1)}$$ 
Note that $f$ is continuous on $\overline{AB}$ with $f(0)=0$ and $f(1)=1$ by Proposition 1 and (2.2.8),  so 
$$\lim_{u\rightarrow 0} f(U)=f(A)=0\mbox{ and } f(B)=f_l(B)=1.\eqno{(5.2)}$$
By (5.2),
$$\lim_{u\rightarrow 0}au+b=0.\eqno{(5.3)}$$
By (5.3),  $b=0$. By using the same argument of (5.3),  $a=1$ if $u\rightarrow 1$ in (5.3).
Thus, 
$$f(U)=u \mbox{ for }U=(u, 0)\in \overline{AB}.\eqno{(5.4)}$$
  Note  that $\pi_{l_k} (D, \overline{AX}, \overline {BC} )=f_{l_k}(X)$ for a lattice point $X=(x, 0)\in \overline{AB}$ by (2.2.1), so by (5.4) and Proposition 1, for $X\in \overline{AB}$,  
$$\lim_{k\rightarrow \infty  } \pi_{l_k^{-1}} (D, \overline{AX}, \overline {BC} )=X.\eqno{(5.5)}$$
Now we  assume that  Theorem  1 does not hold. There exists a subsequence $\{n\}$ and $X\in \overline{AB}$  such that
$$\lim_{n\rightarrow \infty } \pi_{n^{-1}} (D, \overline{AX}, \overline {BC} )\neq X.\eqno{(5.6)}$$
By the same proof of (5.5),  we take a subsequence of $l_{k}$ from $\{n\}$ such that 
$$ \lim_{m\rightarrow \infty}f_{l_k}(X)=\lim_{m\rightarrow \infty  } \pi_{l_k^{-1}} (D, \overline{AX}, \overline {BC} )=X.\eqno{(5.7)}$$
The contradiction in (5.6) and (5.7) tells us that, for any sequence $\{l\}$,
$$\lim_{l \rightarrow \infty} \pi_{l^{-1}} (D, \overline{AX}, \overline {BC} )= X.\eqno{(5.8)}$$
The Theorem   follows from (5.8). $\blacksquare$

\section {Appendix: Proof of Lemmas 2.4.1.}

We only show the first inequality of Lemma 2.4.1:
 $$
 |{\bf P}_l( {\cal H}(D, Z_2, Z_3)\cap {\cal D}^+(D, \tau^+(D), {Z_3}))-{\bf P}_l( {\cal H}(D, Z_1,  Z_2)\cap {\cal D}^+(D, \tau^+(D), Z_2))|\leq c_1 (nl^{-1})^\alpha l^{-2}, \eqno{ (A.1)}$$
 for $Z_1, Z_2$ and $Z_3$ defined in Lemma 2.4.1 and for some positive constants $c_1$ and $\alpha$.
The same proof can be adapted to show the second inequality of Lemma 2.4.1. 
It is easier to understand the proof when $n=1$, so we first prove this case. 
Without loss of generality, we assume that $Z_2={1\over 2}(A+B)\in {\bf N}_{l^{-1}}$ is the center of $\overline{AB}$.
 We still divide $\overline{AB}$ into
$l-1$ many equal sub-intervals with length $l^{-1}$. 
 Let $\bigtriangleup  \alpha\beta\gamma$ be a smaller  equilateral  triangle similar to $\bigtriangleup  ABC$ with the center $Z_2$ of $\overline{\alpha\beta}\subset \overline{AB}$ 
and with
$$d(\alpha,\beta)= l^{-\eta }\mbox{ for some }0<  \eta <1,\eqno{(A.2)}$$
where $\eta$ can be chosen for any number less than one, so we may choose $\eta=1/2$.
In addition, for  a small $0< \delta < \eta$, let $\bigtriangleup  \alpha'\beta'\gamma'$ be an  equilateral  triangle  similar to $\bigtriangleup ABC$ also with the center $Z_2$ of $\overline{\alpha'\beta'}\subset \overline{AB}$  (see Fig. 7)
and
$$d({\alpha',\beta'})=  l^{-\delta }\mbox{ for a small }0< \delta < \eta < 1.\eqno{(A.3)}$$
On ${\cal H}(D, Z_1, Z_2)\cap {\cal D}^+(D, \tau^+(D), Z_2)$, there are two arm paths:  an open $\bullet$-path  $\gamma_1$  and a closed  $\circ$-path $\gamma_2$ from $Z_2$ and $Z_1$ in $\bigtriangleup ABC $ to $\overline{AC}$ and $\overline{BC}$, respectively. There might be many such paths $\{\gamma_1\}$ and $\{\gamma_2\}$. We select $\gamma_1$ to be the {\em innermost} open path: the area $\Gamma_1'$ enclosed by $\gamma_1$ and $\overline{Z_2B}$ and  a part of $\overline{BC}$ is the smallest. Similarly, we select $\gamma_2$ to be the innermost path: the area $\Gamma_2'$  enclosed by $\gamma_2$, $\overline{Z_1 A}$, and  a  part of $\overline{AC}$ is the smallest (see Fig. 7).
\begin{figure}
\begin{center}
\begin{picture}(100,130)(35,10)
\setlength{\unitlength}{0.0125in}%
\begin{tikzpicture}
\thicklines
\begin{scope}[>={Stealth[black]},
              every edge/.style={draw=blue,very thick}]
 \path [-] (7.9, -2.5) edge [bend left=40](9.9,3.3);
 \path [-] (9.5, -0.5) edge [bend left=10](12,1.2);
 \path [-] (6.5, -0.8) edge [bend left=10](6.5,-2.5);
 \path [-] (6.6, -1.6) edge [bend left=10](7.2,-1.6);
\path [-] (7.7, -1.6) edge [bend left=10](8.7,-1.6);
\end{scope}
\begin{scope}[>={Stealth[black]},
              every edge/.style={draw=red,very thick}]
              
   \path [-] (10.2, 3.15) edge [bend right=50](9.5,-0.5);
    \path [-] (7.7, -2.5) edge [bend right=20](3.7,1.2);
     \path [-] (9.5, -0.5) edge [bend right=20](13.5,-2.5);
      \path [-] (5.8, 3.3) edge [bend right=40](9.8,3.5);
       \path [-] (9.8,3.5) edge [bend right=10](10.15,3.16);
\path [-] (6.6, -2) edge [bend left=10](7.45,-2);
\path [-] (7.8, 0.4) edge [bend left=10](8.9,-2.5);
\path [-] (7.56, -1.2) edge [bend left=10](8.6,-1.2);
\end{scope}

\put(0,-80){\line(1,1){250}}
\put(500,-80){\line(-1,1){250}}
\put(495,-80){\line(-1,1){248}}
\put(50,-80){\line(1,1){200}}
\put(450,-80){\line(-1,1){200}}
\put(100,-80){\line(1,1){150}}
\put(400,-80){\line(-1,1){150}}
\put(150,-80){\line(1,1){100}}
\put(350,-80){\line(-1,1){100}}
\put(0,-80){\line(1,0){500}}

\put(240,-37){\circle*{3}}
\put(300,-15){\circle*{3}}
\put(310,-20){$v$}
\put(320,105){$S$}
\put(275,-50){\circle*{3}}

\put(227,-50){\circle*{3}}
\put(207,-62){\circle*{3}}
\put(243,12){\circle*{3}}
\put(208,-23){\circle*{3}}

\put(265,-20){$\bar{r}_1$}
\put(240,-30){${r}_1$}
\put(265,-20){$\bar{r}_1$}
\put(260,-65){$\Gamma_1$}
\put(215,-75){$\Gamma_2$}
\put(220,-40){${r}_2$}
\put(195,-60){$\bar{r}_2$}
\put(500,-85){$\tau^+(B)$}
\put(485,-90){$B$}
\put(250,170){$\tau^+(C)$}
\put(235,165){$C$}
\put(350,-55){$\lambda_1$}
\put(390,-20){$\Lambda_1$}
\put(440,-90){$b'$}
\put(270,120){$\Lambda_2$}
\put(250,120){$c'$}
\put(200,100){$\lambda_2$}
\put(240,70){$\gamma'$}

\put(-10,-90){$A$}
\put(400,-90){$\beta'$}
\put(250,20){$\gamma$}

\put(250,-95){$Z_2$}
\put(235,-95){${Z}_1$}

\put(45,-90){$a'$}
\put(85,-90){$\alpha'$}
\put(150,-90){$\alpha$}
\put(340,-90){$\beta$}
\put(192,-25){$\nu_2$}
\put(230,10){$\nu_1$}
\put(162,-0){$\beta_2$}
\put(240,35){$\beta_1$}


\end{tikzpicture}
\end{picture}
\end{center}
\caption{\em \small {We construct the innermost paths $r_1$ and $r_2$ from $Z_2$ and $Z_1$ to
$\overline {\alpha \gamma}\cup\overline{\beta\gamma}$ at $\nu_1$ and $\nu_2$.
The areas included by $r_1$ and $\bar{r}_1$, and by $r_2$ and $\bar{r}_2$ are $\Gamma_1$ and $\Gamma_2$.  For each point in $r_1$, there is a closed $\circ$-path to $\bar{r}_1$, and for each point in $\bar{r}_1$, there is 
an open $\bullet$-path to $r_1$. For each point in $\bar{r}_2$, there is a closed $\circ$-path to ${r}_2$. For each point in ${r}_2$, there is 
an open $\bullet$-path to $\bar{r}_2$. $\gamma_1$ is selected to be the innermost path from $Z_2$ to $\overline{BC}$, but not 
to $\tau^+(\overline{BC})$. $\gamma_1$ can be viewed  from $Z_2$ along $r_1$ and then along a sub-path  $\beta_1\subset \gamma_1$ from $r_1$  to $S$. There are three arm paths below $\overline{BC}$ (the half-space) from
$S\in\overline{BC}$: $\lambda_1$ and $\lambda_2$  and $\gamma_1$. $\lambda_1$ and $\lambda_2$ are also selected to be the innermost path. 
For each point $v$ in  $\lambda_1$, there is an open $\bullet$-path from the point to $\overline{BC}$.  $\gamma_2$ is selected to be the innermost path from $Z_1$ to $\overline{AC}$. $\gamma_2$ can be viewed from $Z_1$ along $r_2$ and then along $\beta_2\subset \gamma_2$   from $r_2$ to $\overline{AC}$.}}
\end{figure}

$\gamma_1$ and $\gamma_2$ have to  cross  out $\overline{\beta \gamma}\cup\overline{\alpha \gamma}$.
Thus, there are two  arm paths: closed $\circ$-path $r_2$ and open $\bullet$-path $r_1$ from  $Z_1$ and from $Z_2$ to  $\overline{\alpha\gamma}\cup \overline{\beta\gamma}$, respectively. In addition,  there are two paths: a closed $\circ$-path $\beta_2$, an open $\bullet$-path $\beta_1$, a part of $\gamma_2$ and a part of $\gamma_1$,  from  $r_2$,  and from $r_1$, but without using $r_2$ and $r_1$, to  $\overline{AC}$ and $\overline{BC}$, respectively (see Fig. 7).  
If there are two arm paths from $Z_1$ and $Z_2$ to $\overline{\alpha \gamma}\cup \overline{\beta\gamma}$ , then we select the innermost two arm paths, still denoted by $r_1$ and $r_2$
in $\bigtriangleup \alpha\beta\gamma$ intersecting $\overline{\alpha \gamma}\cup \overline{\beta\gamma}$ at $\nu_1$ and $\nu_2$, respectively (see Fig. 7).  
In other words, if $r_1$ is the innermost one, then the area including $r_1$, a part of  $\overline{\alpha \gamma}\cup \overline{\beta\gamma}$, and  $\overline{Z_2\beta}$ is  the smallest (see Fig. 7). Similarly, $r_2$ can be selected as the innermost one in the same way.
Since $r_1$ and $r_2$ are the innermost paths, $\nu_1$ and $\nu_2$ are unique for each configuration. 
Thus, 
\begin{eqnarray*}
 &&{\bf P}_l({\cal H}(D, Z_1, Z_2)\cap {\cal D}^+(D, \tau^+(D), {Z}_2))\\
 &=&\sum_{h_1, h_2}{\bf P}_l({\cal H}(D, Z_1, Z_2)\cap {\cal D}^+(D, \tau^+(D), Z_2),  \nu_1=h_1, \nu_2=h_2),\hskip 3.5cm (A.4)
  \end{eqnarray*}
where the above sum is taken over all $h_1$ and $h_2$ with $h_2, h_1\in \overline{\alpha \gamma}\cup \overline{\beta\gamma}$.
On the other hand, if $r_1$ and $r_2$ are the innermost paths, then by Proposition 2.3 in Kesten (1982),  there are two paths: one closed $\circ$-path  $\bar{r}_1$ and one open $\bullet$-path $\bar{r}_2$ in $\bigtriangleup \alpha\beta\gamma$
from $\nu_1$ to $\overline {Z_2\beta}$ and from $\nu_2$ to $\overline{\alpha Z_1}$, respectively (see Fig. 7).
$\bar{r}_1$ and $\bar{r}_2$ are selected to be the innermost paths such that the areas enclosed by $r_1$, $\bar r_1$ and a part of $\overline {Z_2\beta}$ and by $r_2$, $\bar r_2$ and a part of $\overline {\alpha Z_1}$ are the smallest, respectively.  
Let $\Gamma_1$  be the vertices enclosed by $r_1$ and $\bar{r}_1$ and a part of $\overline{Z_2\beta}$, and 
let $\Gamma_2$  be the vertices enclosed by $r_2$ and $\bar{r}_2$ and a part of $\overline{\alpha{Z}_1}$ (see Fig. 7).
The vertex sets of $\Gamma_1$ and $\Gamma_2$ are smallest when $\bar{r}_1$ and $\bar{r}_2$ are selected to be the innermost paths.
 With these constructions,
\begin{eqnarray*}
 &&{\bf P}_l({\cal H}(D, Z_1, Z_2)\cap {\cal D}^+(D, \tau^+(D), {Z}_2))\\
 &=&\sum_{G_1, G_2}{\bf P}_l({\cal H}(D, Z_1, Z_2)\cap {\cal D}^+(D, \tau^+(D), {Z}_2),  \Gamma_1=G_1, \Gamma_2=G_2),\hskip 4cm (A.5)
  \end{eqnarray*}
where the above sum is taken over all possible vertex sets $G_1$ and $G_2$ in $\bigtriangleup\alpha\beta\gamma$. 
It follows from Proposition 2.3 in Kesten (1982) that if $\Gamma_1=G_1$ and  $\Gamma_2=G_2$ for some fixed vertex sets $G_1$ and $G_2$,
$$\{ \Gamma_1=G_1,\Gamma_2=G_2\} \mbox{ only depends on the configurations  in the vertices of }G_1 \cup G_2.\eqno{(A.6)}$$
It also follows from Proposition 2.3 in Kesten (1982) that for each $u\in r_1$ and each $v\in r_2$ (see Fig. 7),
$$\exists \mbox{ an open $\bullet$-path from $v$ to $\bar{r}_2$ and a closed $\circ$-path  from $u$ to $\bar{r}_1$, respectively}.\eqno{(A.7)}$$

On ${\cal H}(D, Z_1, Z_2)\cap {\cal D}^+(D, \tau^+(D), {Z}_2)$,  there are three arm paths from some point on $\tau^+(\overline{BC})$. One open path is $\gamma_1$ from $Z_2$ meeting $\overline{ BC}$  and $\circ$-adjacent  $S$, and the others are closed $\circ$-paths 
$\lambda_1$ and $\lambda_2$ from $S$ to $\overline{AC}$ and to $\overline{Z_2 B}$, respectively.  Since $\gamma_1$ is the innermost path, $S$ is unique and closed.
With these two closed $\circ$-paths, we go along $\lambda_2$ from $\overline{AC}$ to meet
$S$ and go along $\lambda_1$ from $S$ to $\overline{ZB}$. We assume that $S\in \lambda_1$ (see Fig. 7).
We  select these two closed $\circ$-paths, still denoted by  $\lambda_1$ and $\lambda_2$, to be the innermost paths. 
In other words, if $\lambda_1$ is the innermost one, then the area including  $\lambda_1$,  a part of $\tau^+(\overline{BC})$, and  a part of $\overline{Z\tau^+(B)}$ is  the smallest (see Fig. 7). Let $\Lambda_1$ be the vertices in the above area. Similarly, $\lambda_2$ is selected as the innermost one such that the area included  by $\lambda_2$,  a part of $\tau^+(\overline{BC})$,  and  a part of $\overline{AC}$ is the smallest (see Fig. 6).  Let $\Lambda_2$ be the vertices in the area. We have
\begin{eqnarray*}
 &&{\bf P}_l({\cal H}(D, Z_1, Z_2)\cap {\cal D}^+(D, \tau^+(D), {Z}_2))\\
 &=&\!\!\!\!\!\sum_{G_1, G_2, H_1, H_2}\!\!\!\!{\bf P}_l({\cal H}(D, Z_1, Z_2\cap {\cal D}^+(D, \tau^+(D), {Z}_2),  \Gamma_1=G_1, \Gamma_2=G_2, \Lambda_1=H_1, \Lambda_2=H_2),\hskip .5cm (A.8)
  \end{eqnarray*}
where the above sum is taken over all possible vertex sets $G_1$, $G_2$, $H_1$, and $H_2$.
It follows from Proposition 2.3 in Kesten (1982) that if  there are three arm paths from $S$ defined above with $\Lambda_1=H_1$ and  $\Lambda_2=H_2$ for some fixed vertex sets $H_1$ and $H_2$, then
$$\{ \Lambda_1=H_1, \Lambda_2=H_2\} \mbox{ only depend on configurations  in the vertices of }H_1\cup H_2.\eqno{(A.9)}$$
Similar to (A.6), 
 for each $u\in \lambda_1$ and each $v\in \lambda_2$  (see Fig. 7),
\begin{eqnarray*}
&&\exists \mbox{ an open $\bullet$-path from $v$ to $\tau^+(\overline{BC})\cup \overline{\tau^+(B)Z_2}$};\\
&&\exists \mbox{  an open   $\bullet$-path  from $u$ to $\tau^+(\overline{BC})\cup\overline{AC}$}.\hskip 7.5cm {(A.10)}
\end{eqnarray*}
Now we show that the three arm paths from $S$ will stay outside $\bigtriangleup \alpha'\beta'\gamma'$ as the following lemma.\\

{\bf Lemma A.1.} {\em For the $\delta>0$ defined in (A.3), there exist $c_i >0$ for $i=1,2$ such that}
$${\bf P}_l ({\cal H}(D, Z_1, Z_2)\cap {\cal D}^+(D, \tau^+(D), {Z}_2), \lambda_1\cup \lambda_2\cap \bigtriangleup \alpha' \beta' \gamma'\neq \emptyset )\leq c_1l^{-2-c_2\delta}.$$

{\bf Proof.}  We construct an equilateral  triangle $\bigtriangleup a'b'c'$  with the center $Z_2$ on $\overline{a'b'}$ (see Fig. 7) similar to $\bigtriangleup ABC$ such that $Z_2$ is the center of $\overline{a'b'}$ and 
$$d(a', b')= \kappa d(A, B)\mbox{ for a small } 0< \kappa < 1.\eqno{}$$
On ${\cal H}(D, Z_1, Z_2)\cap {\cal D}^+(D, \tau^+(D), {Z}_2)$,
there are two arm paths from $Z_1$ and $Z_2$ to the boundary of $\bigtriangleup a'b'c'$ outside of $\Lambda_1\cup\Lambda_2$ denoted by  event ${\cal G}_0$. We consider the line $L_1$ containing $\overline{c'b'}$ and the $(B-b')$-parallelogram  with the center at $S$ denoted
by $S(B-b')$.
On ${\cal H}(D, Z_1, Z_2)\cap {\cal D}^+(D, \tau^+(D), {Z}_2)$,
 three arm paths  from $S$ below $\overline{BC}$ but above $L$ meet the boundary  of $S(B-b')$,
  denoted by event ${\cal G}_1$ (see Fig. 7).
If $\lambda_1\cap \bigtriangleup \alpha' \beta' \gamma'\neq \emptyset$, then $\lambda_1$ first meets $L_1$ and then meets
$v\in \bigtriangleup \alpha'\beta'\gamma'$ (see Fig. 7).
Since $\lambda_1$ is innermost,   there is an open $\bullet$-path in $\Lambda_1$  from $v$ to $\overline{BC}$ disjoint to the three arm paths in ${\cal G}_1$ for some $v$, 
 denoted by event ${\cal G}_2$.  
 Thus, 
  the  open $\bullet$-path from $v$,  the two arm paths from $Z_1$ and $Z_2$, and the three arm paths from $S$ to $L_2$ are disjoint (see Fig. 7). In other words,  three events ${\cal G}_0$, ${\cal G}_1$, and ${\cal G}_2$ occur disjointly (see section 2.3 in Grimmett (1999) for the definition of the disjoint events). 
 By Reimer's inequality, (2.1.3), (2.1.6),  and (2.1.7), there exist $c_i$ for $i=1,2$,  independent of $l$ and $\epsilon$ such that
\begin{eqnarray*}
&&{\bf P}_l ({\cal H}(D, Z_1, Z_2)\cap {\cal D}^+(D, \tau^+(D), {Z}_2), \lambda_1\cap \bigtriangleup \alpha' \beta' \gamma'\neq \emptyset )\\
&\leq & {\bf P}_l ({\cal G}_0) {\bf P}_l ({\cal G}_1){\bf P}_l ({\cal G}_2)\leq c_1  l^{-1} l^{-1}  l^{-c_2\delta}\leq c_1l^{-2-c_2\delta}.\hskip 6.5cm{(A.11)}
\end{eqnarray*}
The same proof also implies that
$${\bf P}_l ({\cal H}(D, Z_1, Z_2)\cap {\cal D}^+(D, \tau^+(D), {Z}_2),  \lambda_2\cap \bigtriangleup \alpha' \beta' \gamma'\neq \emptyset )\leq  c_1l^{-2-c_2\delta}.\eqno{(A.12)}$$
Thus, Lemma A.1 follows from  (A.11) and (A.12). $\blacksquare$\\

   We denote by ${\cal G}_3$ there are two arm paths from $\bigtriangleup abc$ to  $\overline{a'c'}\cup\overline{b'c'}$, where 
$\bigtriangleup abc$ is an equilateral  triangle such that $Z_2$ is the center of $\overline{ab}$ and
$d(a, b)=2d(\alpha, \beta)$. 
Let 
$${\cal D}_1={\cal D}_1(\delta)={\cal G}_3\cap \{\exists \,\, \lambda: \lambda_1 \cap \bigtriangleup \alpha' \beta' \gamma'= \emptyset\} \cap \{\exists\,\, \lambda_2: \lambda_2 \cap \bigtriangleup \alpha' \beta' \gamma'= \emptyset\}.$$
Note that ${\cal G}_3$ and ${\cal G}_1$ are disjoint, so by (2.1.7) for $k=2$ and $k=3$, there exists $c_1$ such that
$${\bf P}_l({\cal D}_1)\leq {\bf P}_l({\cal G}_1\cap {\cal G}_3)\leq {\bf P}_l({\cal G}_1){\bf P}_l( {\cal G}_3)\leq c_1l^{-\eta} l^{-1}=c_1l^{-1-\eta}.\eqno{(A.13)}$$
By  Lemma A.1, 
 \begin{eqnarray*}
 &&{\bf P}_l ({\cal H}(D, Z_1, Z_2)\cap {\cal D}^+(D, \tau^+(D), {Z}_2),  {\cal D}_1^C)\\
 &\leq &{\bf P}_l ({\cal H}(D, Z_1, Z_2)\cap {\cal D}^+(D, \tau^+(D), {Z}_2), \lambda_1\cup \lambda_2\cap \bigtriangleup \alpha' \beta' \gamma'\neq \emptyset )\leq  c_1l^{-2-c_2\delta}.\hskip 2.5cm {(A.14)}
 \end{eqnarray*}
 
 Now we show that $\beta_1$ and $\beta_2$ could not  be very close. Let ${\cal D}_2={\cal D}_2(\epsilon)$ be the event (see Fig. 8) that
$$ d(\gamma_1, \beta_2)\geq l^{-\eta-\epsilon M} \mbox{ and } d(\gamma_2, \beta_1) \geq  l^{-\eta-\epsilon M},\eqno{}$$
where  $M$ is a positive constant independent of $\epsilon$ and $l$ with $M^{-1}\geq \epsilon$. 
\\

{\bf Lemma A.2.} {\em If  $0<\epsilon< \eta$ defined above,  there exist  constants  $M$  and $c_i=c_i( M)$ for $i=1,2$, independent of $\epsilon$ and $l$, such that
}$${\bf P}_l ( {\cal H}(D, Z_1, Z_2)\cap {\cal D}^+(D, \tau^+(D), {Z}_2), {\cal D}_1, {\cal D}_2^C)\leq c_1l^{-(2+c_2\epsilon) }.$$

\begin{figure}
\begin{center}
\begin{picture}(140,220)(40,-20)
\setlength{\unitlength}{0.0125in}%
\begin{tikzpicture}
\thicklines
\begin{scope}[>={Stealth[black]},
              every edge/.style={draw=blue,very thick}]
 \path [-] (8, -2.5) edge [bend right=40](8.5,-0.5);
 \path [-] (8.5, -0.5) edge [bend left=100](10.1,3.3);

 \path [-] (2.5, 0) edge [bend left=10](7.8,0);
\end{scope}
\begin{scope}[>={Stealth[black]},
              every edge/.style={draw=red,very thick}]
              
   \path [-] (10.2, 3.15) edge [bend right=20](15.5,-2.5);
   \path [-] (7.7, -2.5) edge [bend left=40](7.8,-0.6);
   \path [-] (7.8, -0.6) edge [bend right=60](4.5,2);
      \path [-] (6.8, 4.3) edge [bend right=10](9.8,3.5);
\path [-] (8, 0) edge [bend right=10](12.1,0);
\end{scope}

\put(0,-80){\line(1,1){250}}
\put(500,-80){\line(-1,1){250}}
\put(50,-80){\line(1,1){200}}
\put(450,-80){\line(-1,1){200}}
\put(100,-80){\line(1,1){150}}
\put(400,-80){\line(-1,1){150}}
\put(150,-80){\line(1,1){100}}
\put(350,-80){\line(-1,1){100}}
\put(0,-80){\line(1,0){500}}
\put(0,-40){\line(1,0){500}}
\put(240, -10){\framebox(20, 20)[br]{T}}


\put(270,-30){${\gamma_1}$}
\put(220,-40){$\gamma_2$}
\put(490,-90){$\tau^+(B)$}
\put(250,180){$\tau^+(C)$}
\put(410,-20){$\lambda_1$}
\put(440,-90){$b'$}
\put(250,120){$c'$}
\put(250,70){$c$}
\put(220,120){$\lambda_2$}

\put(-10,-90){$A$}
\put(400,-90){$b$}
\put(250,20){$\gamma$}
\put(-10,-30){$y=2l^{-\eta-\epsilon}$}

\put(250,-95){$Z_2$}
\put(235,-95){$Z_1$}

\put(85,-90){$$}
\put(150,-90){$\alpha$}
\put(100,-90){$a$}
\put(50,-90){$a'$}
\put(340,-90){$\beta$}

\put(370,60){$\tau^+(\overline{BC})$}

\end{tikzpicture}
\end{picture}
\end{center}
\caption{\em \small{On ${\cal D}_2^C$,  $\gamma_1$ and $\gamma_2$ will meet in a square $T$ with its center inside $\bigtriangleup \alpha \beta\gamma$ and above $y=2l^{-\eta-\epsilon}$. By (A.6), there are six arm paths from
$T$ to $\overline{ac}\cup\overline{bc}$. }}
\end{figure}

{\bf Proof.}  We estimate the probability of  $d(\gamma_1, \beta_2) < l^{-\eta-\epsilon M}$, denote by ${\cal D}_2'^C$.
We first estimate the probability of event that $d(\beta_1, Z_2) < l^{-\eta-\epsilon}$.
On $d(\beta_1, Z_2) < l^{-\eta-\epsilon}$, there exists an open $\bullet$-path, a part of $\beta_1$,  from the boundary of $Z_2+ [-l^{-\eta-\epsilon}, l^{-\eta-\epsilon}]^2$ to
$\overline{\alpha \gamma} \cup \overline{\beta\gamma}$ outside $\Gamma_1$ and $\Gamma_2$, denote by event ${\cal G}_4$. 
Thus, by (2.1.3), there exist $c_1$ and $c_2=c_2$ such that
$$\!{\bf P}_l({\cal G}_4)\leq  {\bf P}_l(\exists \mbox{ an open $\bullet$-path from the boundary }Z_2+ [-l^{-\eta-\epsilon}, 
l^{-\eta-\epsilon}]^2 \mbox{ to } \overline{\alpha \gamma} \cup \overline{\beta\gamma})\leq c_1l^{-c_2\epsilon}.\eqno {(A.15)}$$
Note that the above open path in ${\cal G}_4$ is outside of $\Gamma_1$ and $\Gamma_2$ and away 
$\lambda_1$ and $\lambda_2$ on ${\cal D}_1$, 
 so by (A. 13), (A.15) and  Reimer's inequality,  there exist $c_i$ for $i=1,2$ such that
\begin{eqnarray*}
&&{\bf P}_l ({\cal H}(D, Z_1, Z_2)\cap {\cal D}^+(D, \tau^+(D), {Z}_2), {\cal D}_1, d(\beta_1, Z_2) < l^{-\eta-\epsilon})\\
&\leq &  {\bf P}_l ( {\cal G}_4, {\cal H}(\bigtriangleup \alpha \beta\gamma, Z_1, Z_2),  {\cal D}_1)\leq   {\bf P}_l ( {\cal G}_4){\bf P}_l ({\cal H}(\bigtriangleup \alpha \beta\gamma, Z_1, Z_2)) {\bf P}_l ( {\cal D}_1)\\
&\leq  &c_1l^{-c_2\epsilon} l^{-\eta} l^{-1+\eta}l^{-1}\leq c_1l^{-2-c_2\epsilon}.\hskip 9.5cm {(A.16)}
\end{eqnarray*}
Similarly, we have
$${\bf P}_l ({\cal H}(D, Z_1, Z_2)\cap {\cal D}^+(D, \tau^+(D), {Z}_2), {\cal D}_1, d(\beta_2, Z_2) < l^{-\eta-\epsilon})\leq c_1l^{-2-c_2\epsilon}.\eqno {(A. 17)}$$
Here we work on the innermost paths $\beta_1$ and $\beta_2$ in (A.16) and (A.17). For any open $\bullet$-path $\dot \beta_1$ from $r_2$ to $S$ and any closed 
$\circ$-path  $\dot\beta_2$ from $r_2$ to $\overline{AC}$,  the same proofs of (A.16) and (A.17) also work.
\begin{figure}
\begin{center}
\begin{picture}(100,220)(35,10)
\setlength{\unitlength}{0.0125in}%
\begin{tikzpicture}
\thicklines
\begin{scope}[>={Stealth[black]},
              every edge/.style={draw=blue,very thick}]
 \path [-] (7.9, -2.5) edge [bend left=20](15.8,-2);
 \path [-] (9.5, -1.2) edge [bend right=10](5,5);
\end{scope}
\begin{scope}[>={Stealth[black]},
              every edge/.style={draw=red,very thick}]
              
    \path [-] (7.7, -2.5) edge [bend right=20](3.7,1.2);
\path [-] (5.5, 0.2) edge [bend right=40](14,0.5);
\end{scope}

\put(0,-80){\line(1,1){250}}
\put(500,-80){\line(-1,1){250}}
\put(0,-80){\line(1,0){500}}
\put(0,-80){\framebox(70,70)[br]{\mbox{$$}}}
\put(70,-80){\framebox(70,70)[br]{\mbox{$$}}}
\put(140,-80){\framebox(70,70)[br]{\mbox{$$}}}
\put(210,-80){\framebox(70,70)[br]{\mbox{$$}}}
\put(280,-80){\framebox(70,70)[br]{\mbox{$$}}}
\put(350,-80){\framebox(70,70)[br]{\mbox{$$}}}
\put(230,-80){\framebox(40,40)[br]{\mbox{$$}}}
\put(420,-80){\framebox(70,70)[br]{\mbox{$$}}}

\put(215,-60){${r}_2$}
\put(240,-30){${\beta}_2$}
\put(260,-65){$S$}
\put(300,-30){${w}_1$}
\put(300,-65){${w}_2$}
\put(390,-55){$r_1$}
\put(-10,-85){$\alpha$}
\put(500,-85){$\beta$}
\put(250,170){$\gamma$}
\put(0,-5){$y=2l^{-\eta-\epsilon}$}
\put(340,-80){$T$}

\put(250,-95){$Z_2$}
\put(235,-95){${Z}_1$}


\end{tikzpicture}
\end{picture}
\end{center}
\caption{\em \small {The rectangle $[\alpha, \beta]\times [0, 2l^{\eta-\epsilon}] $ is divided by squares. $\beta_2$ will meet a square $T$ at $w_1$. There are three arm paths from $w_1$ to $\overline{ ac} \cup \overline{bc}$. $r_1$ will also meet the square $T$ at $w_2$. There is an open $\bullet$ path from $w_2$ to $\overline{ ac} \cup \overline{bc}$. In addition, there are two arm paths from $Z_1$ and $Z_2$ to the boundary of $S=Z_2+ [-l^{-\eta-\epsilon}, l^{-\eta-\epsilon}]^2$. }}
\end{figure}

We now estimate ${\cal D}_2'^C$.   
We construct   equal squares $\{T\}\subset {\bf Z}_{l^{-\eta -M\epsilon}}^2$ such that the centers of these squares in $ \bigtriangleup \alpha \beta\gamma$.
Note that each $T$ is a square with side length $l^{-\eta -\epsilon M }$,
so there are at most $l^{2M\epsilon}$ many such squares.  On ${\cal D}_2'^C$, there is a $T$ such that
$$T\cap \beta_2 \cap \gamma_1\neq \emptyset.\eqno{(A.18)}$$
We first suppose that  the center of each $T$ is  above $y=2l^{-\eta-\epsilon}$.
Since $\gamma_2$ and $\gamma_1$ are the innermost paths, on (A.18),  by Proposition 2.3 in Kesten (1982), 
there are a closed $\circ$-path in $\Gamma_1'$ from $T$ to $\overline{BC}\cup \overline{ Z_2 C}$
and an open $\bullet$-path in $\Gamma_2'$  from $T$ to $\overline{AC}\cup \overline{AZ_1}$ besides $\gamma_2$ and $\gamma_1$ (see Fig. 8). 
Thus, there  is a square $T$ with its center  inside $\bigtriangleup \alpha\beta\gamma\cap \{y \geq  2l^{-\eta-\epsilon}\}$ such that there are six arm paths from  $T$ to the boundary of  $\overline{ ac}\cup \overline{ bc} \cup  \{y\geq l^{-\eta-\epsilon}\}$ (see Fig. 8). We denote the above event by ${\cal G}_5$.  We also denote by ${\cal G}_6$
  there exist two arm paths from $Z_1$ and $Z_2$  to the upper, the left and the right boundaries  of $Z_2+ [-l^{-\eta-\epsilon}, l^{-\eta-\epsilon}]^2$.
Note that  $ {\cal D}_1$,  ${\cal G}_5$ and ${\cal G}_6$ are independent since they occur in different vertex sets, so by Reimer's inequality,
\begin{eqnarray*}
&&{\bf P}_l ({\cal H}(D, Z_1, Z_2)\cap {\cal D}^+(D, \tau^+(D), {Z}_2), {\cal D}_1, {\cal D}_2'^C, (A. 18) \mbox{ with the center of } T \mbox{ above } y= 2l^{-\eta-\epsilon})\\
&\leq &  {\bf P}_l ({\cal D}_1,{\cal G}_5,{\cal G}_6)\leq {\bf P}_l ({\cal D}_1) {\bf P}_l ({\cal G}_5) {\bf P}_l ({\cal G}_6).\hskip 8.5cm {(A. 19)}
\end{eqnarray*}
By (2.1.5), note that there are $l^{2\epsilon M}$ many choices to select $T$, so  there exists  $c_1$   such that 
$${\bf P}_l ({\cal G}_5) \leq c_1 l^{2M\epsilon }( l^{-\eta-\epsilon M }  l^{\eta+\epsilon} )^{2+d_2}\leq c_1 l^{-\epsilon(d_2(M-1)-2) } ,\eqno{(A.20)}$$
where $d_2$ is the constant defined in (2.1.5).
By (2.1.6) and (A.13), there exist $c_i$ for $i=2,3$ such that
$${\bf P}_l ({\cal G}_6)\leq c_2 l^{-1+\eta+\epsilon},\,\,\,\, {\bf P}_l ({\cal D}_1)\leq c_3 l^{-\eta-1}.\eqno{(A. 21)}$$
By (A.19)--(A.21), if $M$ is taken to be $d_2 (M-1)-2 > 1$, then there exist $c_6$ and $c_7$ independent of  $\epsilon$ and $l$ such that
\begin{eqnarray*}
&&{\bf P}_l ( {\cal H}(D, \bar{Z}, Z)\cap {\cal D}^+(D, \tau^+(D), {Z}), {\cal D}_1, {\cal D}_2'^C, (A. 18) \mbox{ with the center of } T \mbox{ above } y= 2l^{-\eta-\epsilon})\\
&\leq &c_1 l^{-\epsilon(d_2(M-1)-2)} c_3 l^{-1+\eta+\epsilon}  c_4l^{-\eta-1}\leq  c_6l^{-(2+c_7\epsilon )}.\hskip 6.5 cm {(A.22)}
\end{eqnarray*}

Now we focus on the center of the $T$ in (A.18) below the horizontal line $y=2l^{\eta-\epsilon}$.  We divide the rectangle  $[\alpha, \beta]\times [0, 2l^{-\eta-\epsilon}]$ into 
equal squares 
$$[\alpha, \alpha+  2l^{-\eta-\epsilon}]\times [0, 2l^{-\eta-\epsilon}],\cdots, [ \alpha+ 2 j l^{-\eta-\epsilon}, \alpha+ 2 (j+1) l^{-\eta-\epsilon}]\times [0, 2l^{-\eta-\epsilon}]$$ 
such that $\beta\in  [ \alpha+ 2 j l^{-\eta-\epsilon}, \alpha+  2(j+1) l^{-\eta-\epsilon}]$ (see Fig. 9).
There are at most $O(l^\epsilon)$ many such squares such that $T$ is one of them. 
Thus, $\beta_2$ passes a lattice points $w_1\in T$. Since $\beta_2$ is innermost, there are three arm paths  from $w_1$ to 
$\overline{ac}\cup \overline{bc}$, denoted by ${\cal G}_7$.
In addition,  $r_1$ also passes a lattice point $w_2\in T\setminus \{Z_2+ [-l^{-\eta-\epsilon}, l^{-\eta-\epsilon}]\}$  and then from $w_2$ to 
$\overline{ac}\cup \overline{bc}$ below the above three arm paths without using points of $Z_2+ [-l^{-\eta-\epsilon}, l^{-\eta-\epsilon}]$ again (see Fig. 9).
Thus, there are three arm path from $w_1$ to $\overline{ac}\cup \overline{bc}$  above
the $x$-axis and one disjoint open $\bullet$-path from $w_2$ to $\overline{ac}\cup \overline{bc}$.
In addition, on $ d(\beta_2, Z_2)\geq l^{-\eta-\epsilon}$, the above three paths in ${\cal G}_7$ and two arm paths in ${\cal G}_6$ are also disjoint.
On the other hand, $ {\cal D}_1$ only depends on the configurations  outside $\bigtriangleup abc$.
 By (A.16) and  Reimer's inequality
\begin{eqnarray*}
&&{\bf P}_l ( {\cal H}(D, \bar{Z}, Z)\cap {\cal D}^+(D, \tau^+(D), {Z}), {\cal D}_1,  {\cal D}_2'^C, (A. 18) \mbox{ with  the center of } T \mbox{ below } y= 2l^{-\eta-\epsilon})\\
&\leq &{\bf P}_l ( {\cal H}(D, \bar{Z}, Z)\cap {\cal D}^+(D, \tau^+(D), {Z}), {\cal D}_1,d(\beta_2, Z_2)\geq l^{-\eta-\epsilon},  (A. 18) \mbox{ with the center of } T \mbox{ below } y= 2l^{-\eta-\epsilon})\\
&&+c_1l^{-2-c_2\epsilon }\\
&\leq & {\bf P}_l ({\cal G}_6,  {\cal D}_1, {\cal G}_7)\leq {\bf P}_l ( {\cal G}_6) {\bf P}_l ({\cal D}_1) {\bf P}_l ({\cal G}_7)+ c_1l^{-2-c_2\epsilon}.\hskip 6.5cm {(A.23)}
\end{eqnarray*}
By (2.1.7) for  $k=3$  to estimate the three arm paths from $w_1$, by (A.13) and (A.15)
with translation invariance to estimate   one path from $w_2$,  we use (A.23) together with the estimates in (A.21),
 $${\bf P}_l ( {\cal D}_1) {\bf P}_l ({\cal G}_6) {\bf P}_l ({\cal G}_7)\leq 
c_1l^{-1-\eta}  l^{-1+\eta+\epsilon} l^{\epsilon} (l^{-\eta-\epsilon}/l^{-\eta})^2 l^{-c_2\epsilon}\leq c_1 l^{-2-c_2\epsilon}.\eqno{(A. 24)}$$
By (A.22), (A.23) and (A.24), 
$${\bf P}_l ( {\cal H}(D, \bar{Z}, Z)\cap {\cal D}^+(D, \tau^+(D), {Z}), {\cal D}_1, {\cal D}_2'^C) \leq c_1 l^{-2-c_2\epsilon}.\eqno{(A. 25)}$$
Similarly, if ${\cal D}_2''^C$ be the event that $d(\gamma_2,\beta_1) < l^{-\eta-\epsilon M}$, then by the same proof of (A.25), 
$${\bf P}_l ( {\cal H}(D, \bar{Z}, Z)\cap {\cal D}^+(D, \tau^+(D), {Z}), {\cal D}_1, {\cal D}_2''^C) \leq  c_1 l^{-2-c_2\epsilon}.\eqno{(A. 26)}$$
Lemma A.2 follows from (A.25) and (A. 26).  $\blacksquare$\\

On $ {\cal H}(D, Z_1, Z_2)\cap {\cal D}^+(D, \tau^+(D), {Z}_2)$, there exist a closed  $\circ$-path and an open $\bullet$-path  from $r_2$   to $\overline{AC}$,  and from $r_1$  to $S$ outside $\Gamma_1\cup\Gamma_2\cup \Lambda_1\cup\Lambda_2$,
respectively, as we define before (see Fig. 7).  We denote them by $\{\dot \beta_1\}$ and $\{\dot\beta_2\}$, where $\beta_1$ and $\beta_2$ are one of them.  Let ${\cal D}_3$ be the event that  there are an open $\bullet$-path $\dot{\beta}_1$ and  a closed $\circ$-path $\dot{\beta}_2$ form  $\{\dot \beta_1\}$ and $\{\dot\beta_2\}$
such that 
$$d( \dot{\beta}_1, \bar{r}_1) > \sqrt{3}/2l^{-1}\mbox{ and }   d( \dot{\beta}_2, \bar{r}_2) > \sqrt{3}/2l^{-1}.\eqno{(A. 27)}$$
Event ${\cal D}_3$ implies that  $\dot{\beta}_1$ has  to be away from $ \bar{r}_1$, and   $\dot{\beta}_2$  has to be away from  $ \bar{r}_2$.
 We show the following lemma.\\
\\
\begin{figure}
\begin{center}
\begin{picture}(80,180)(-20,-20)
\setlength{\unitlength}{0.0125in}%
\begin{tikzpicture}
\thicklines
\begin{scope}[>={Stealth[black]},
              every edge/.style={draw=blue,very thick}]
\path [-] (7.9, -2.5) edge [bend left=40](9.1,-0.5);
 \path [-] (6.5, -0.8) edge [bend left=10](6.5,-2.5);
 \path [-] (9.2, -1) edge [bend right=20](10,3.3);
\path [-] (9.1, -1) edge [bend left=40](8.6,-0.7);
\path [-] (8.3, -1) edge [bend left=40](8.5,0.5);
\path [-] (9.2, -1) edge [bend right=40](8.5,0.5);

\end{scope}
\begin{scope}[>={Stealth[black]},
              every edge/.style={draw=red,very thick}]
    \path [-] (10.2, 3.15) edge [bend right=20](13.5,-2.5);
    \path [-] (7.7, -2.5) edge [bend right=20](6.55,-0.78);
      \path [-] (6.8, 4.3) edge [bend right=40](9.8,3.5);
\path [-] (9.1, -0.5) edge [bend left=10](8.9,-2.5);
\path [-] (7.1, -1.2) edge [bend right=20](5.5,3);
\path [-] (9.2, -1) edge [bend right=70](6.2,2);

\end{scope}

\put(0,-80){\line(1,1){250}}
\put(500,-80){\line(-1,1){250}}
\put(150,-80){\line(1,1){100}}
\put(350,-80){\line(-1,1){100}}
\put(0,-80){\line(1,0){500}}
\put(0,-45){\line(1,0){500}}
\put(280, -40){\framebox(20, 20)[br]{T}}
\put(250,-80){\circle*{3}}
\put(243,-80){\circle*{3}}

\put(291,-32){\circle*{3}}
\put(293,-35){\scriptsize$\mu$}
\put(286,-16){\circle*{3}}
\put(280,-12){\scriptsize$\nu_1$}
\put(318,103){$S$}

\put(275,-50){$\bar{r}_1$}
\put(260,-40){${r}_1$}
\put(260,-65){$\Gamma_1$}
\put(215,-75){$\Gamma_2$}
\put(220,-50){${r}_2$}
\put(210,-60){$\bar{r}_2$}
\put(490,-90){$\tau^+(B)$}
\put(250,180){$\tau^+(C)$}
\put(320,-25){$\Gamma_1'$}
\put(390,-20){$\Lambda_1$}
\put(270,120){$\Lambda_2$}
\put(-10,-90){$A$}
\put(-10,-40){$y=l^{-\eta-\epsilon}$}
\put(250,20){$\gamma$}
\put(250,-95){\footnotesize $Z_2$}
\put(240,-95){\footnotesize $Z_1$}

\put(150,-90){$\alpha$}
\put(340,-90){$\beta$}
\put(200,20){${\beta}_2$}
\put(300,35){${\beta}_1$}
\put(270,45){${\zeta'}$}

\put(370,60){$\tau^+(\overline{BC})$}

\end{tikzpicture}
\end{picture}
\end{center}
\caption{\em {The graph shows that some $\dot{\beta}_1$ has to reach to $\bar{r}_1$  away in a distance more than $\sqrt{3}l^{-1}/2$, called event  ${\cal D}_3'$. On ${\cal D}'^C_3$, there is $\mu\in \beta_1$  with $d(\mu, \bar r_1)=\sqrt{3} l^{-1}/2$. If $\mu$ is near $\nu_1$,  there is  an extra $\circ$-path
$\zeta'$ from the box $\nu_1+[-l^{-\eta-2\epsilon M}, l^{-\eta-2\epsilon M}]^2$ with a length at least $l^{-\eta-\epsilon M}$ before reaching to other $\circ$-paths. 
If $\mu$ is not near $\nu_1$ and $Z$, there  are six arm paths from $\mu$ to the boundary of  $\mu+[-l^{-\eta-2M\epsilon }, l^{-\eta-2M\epsilon }]^2$.}}
\end{figure}

{\bf Lemma A.3.} {\em If $\eta$ and $ \epsilon$ are defined above for a small $\epsilon$, then there exist $c_1$ and $c_2$ independent of $l$ and $\epsilon$ such that
$${\bf P}_l ( {\cal H}(D, Z_1, Z_2)\cap {\cal D}^+(D, \tau^+(D), {Z}_2), {\cal D}_1, {\cal D}_2, {\cal D}^C_3)\leq c_1  l^{-(c_2\epsilon+2) }.$$}

{\bf Proof.}   
We only show that $d( \dot{\beta}_1, \bar{r}_1) \leq \sqrt{3}l^{-1}/2$ for all $\dot \beta_1$ in Lemma A.3, denoted by event ${\cal D}_3'^C$   with a probability less than $ c_1l^{-(c_2\epsilon+2) }$.  The other case in ${\cal D}_3^C$  can be proved by the same estimate. 
 On ${\cal D}'^C_3$,  we know that  there is $\mu\in \beta_1$ for  the innermost open $\bullet$-path $\beta_1$ from $\{\dot{ \beta}_1\}$ such that
$$d( \mu, \bar{r}_1) \leq   \sqrt{3} l^{-1}/2\mbox{ and all the other open $\bullet$-paths in $\{\dot{\beta}_1\}$ have to use $\mu$}.\eqno{}$$
Thus,  there is a lattice point of $\bar r_1$ with $\circ$-adjacent to $\mu$.
If there are many such $\{\mu\}$, we select the first $\mu$ along $\gamma_1$ from $Z_2$.
 We consider the closed $\circ$-cluster
adjacent $\mu$ outside $\Gamma_1'$.  The cluster has to $\circ$-adjacent $\overline{AB}\cup\overline{BC}$, otherwise, by  the circuit lemma,
there is an open $\bullet$-circuit, the boundary of the $\circ$-cluster,  containing $\mu$, so we can find an open $\bullet$-path $\dot \beta_2$ without using $\mu$.
The contradiction tells that $\circ$-cluster contains a closed $\circ$-path $\zeta$ $\circ$-adjacent to $\mu$ outside of  $\Gamma_1'$ to  the boundary of $\overline{AB}\cup\overline{BC}$ (see Fig. 10).  
  
 We divide our examination of the situations of $\mu$ into the following  two cases: (a) $d( \nu_1, \mu) \leq l^{-\eta-2M\epsilon }$, 
   (b) (a) does not hold, where $\eta$ is defined in (A.2) and $M$ is defined in (A.14).
On case (a) and on ${\cal D}_1\cap {\cal D}_2\cap {\cal D}_3'^C$,  we know that $d(\mu, \lambda_1\cup\lambda_2) \geq l^{-\epsilon}$, 
$ d (\mu, \gamma_2\cap \bigtriangleup  \alpha'\beta'\gamma')\geq l^{-\eta -\epsilon M}$. Thus,  for a small $\epsilon>0$, on ${\cal D}_1\cap {\cal D}_2\cap {\cal D}_3'^C$,
a sub-path $\zeta'\subset \zeta$ from $\nu_1+[-l^{-\eta-2M\epsilon }, l^{-\eta-2M\epsilon }]^2$ to $\nu_1+[-l^{-\eta-M\epsilon }, l^{-\eta-M\epsilon }]^2$  will be away  from $\gamma_2$,  $\lambda_1\cup\lambda_2$,  and outside  $\Gamma_1'$ (see Fig. 10), denoted by event ${\cal G}_8(\nu_1)$. 
Note that if $\Gamma_1$, $\gamma_2$,  $\Gamma_1'$, $\Lambda$, $\Lambda_2$ and $\nu_1=h$ are fixed, so  ${\cal G}_8(h)$ is independent of  the configurations of these sets
by  Proposition 2.3 of Kesten (1982).  On the other hand, for a fixed $h$, 
$${\bf P}_l( {\cal G}_8(h))\leq {\bf P}_l(\exists \mbox{ an open $\bullet$-path  from } h+[-l^{-\eta-2M\epsilon }, l^{-\eta-2M\epsilon }]^2\mbox{ to } h+[-l^{-\eta-M\epsilon }, l^{-\eta-M\epsilon }]^2).\eqno{(A. 28)}$$
By translation invariance and using (2.1.3) for the right side of (A.28), there exist $c_1$ and $c_2$ such that
$${\bf P}_l({\cal G}_8(h))\leq c_1 (l^{-\eta-2M\epsilon}/ l^{-\eta-M\epsilon})^{d_1}= c_1 l^{-d_1M\epsilon} \leq c_1l^{-c_2\epsilon}.\eqno{(A.29)}$$
By  the estimates in (A.11), (A.13) and (A.29),  if we select $\epsilon$ to be small,  
  then there exists $c_i$ for $i=1,2,3,4$ such that
\begin{eqnarray*}
&&{\bf P}_l ({\cal H}(D, Z_1, Z_2)\cap {\cal D}^+(D, \tau^+(D), {Z}_2),{\cal D}_1, {\cal D}_2, {\cal D}'^C_3,\mbox{ case (a)})\\
&\leq &\sum_{G_1, G_1', G_2, H_1, H_2, h}{\bf P}_l ({\cal H}(D, Z_1, Z_2)\cap {\cal D}^+(D, \tau^+(D), {Z}_2),{\cal D}_1, {\cal D}_2, {\cal D}'^C_3, \Gamma_1=G_1, \Gamma_1'=G_1',\\ 
&&\hskip 8cm \Gamma_2=G_2, \Lambda_1=H_1, \Lambda_2=H_2, \nu_1=h)\\
&\leq &\sum_{G_1, G_1' G_2, H_1, H_2, h}{\bf P}_l ({\cal G}_8(h)){\bf P}_l ({\cal H}(D, Z_1, Z_2)\cap {\cal D}^+(D, \tau^+(D), {Z}_2),{\cal D}_1, {\cal D}_2,  \Gamma_1=G_1, \Gamma_1'=G_1'\\
&&\hskip 9cm \Gamma_2=G_2, \Lambda_1=H_1, \Lambda_2=H_2, \nu_1=h)\\
&\leq & c_1l^{-c_2\epsilon}{\bf P}_l({\cal G}_0, {\cal G}_1)\leq  c_1l^{-c_2\epsilon}{\bf P}_l({\cal G}_0){\bf P}_l({\cal G}_1) \leq  c_1  l^{-c_2\epsilon}c_3 l^{-2}\leq c_4l^{-(2-c_2\epsilon)},\hskip 3 cm {(A.30)}
\end{eqnarray*}
where ${\cal G}_0$ and ${\cal G}_1$ are the events defined in Lemmas A.1, and the  sum in (A.30) takes over all possible sets $G_1$, $G_1'$, $G_2$ $H_1$, $H_2$ and  $h$. 
By the same proof in case (a) above and using the innermost of $\beta_1$,
if $d(u_1, \nu_1)\leq l^{-\eta-M\epsilon}$ for $u_1= \beta_1\cap r_1$, then ${\cal G}_8$ occurs. Similarly, we can work on the same situation for $\beta_2$.
Thus,  there exist $c_1$ and $c_2$ such that for $i=1,2$
$${\bf P}_l({\cal H}(D, Z_1, Z_2)\cap {\cal D}^+(D, \tau^+(D), {Z}_2),  d(u_i, \nu_i)\leq l^{-\eta-M\epsilon})\leq c_1 l^{-2-c_2\epsilon}.\eqno{( A.31)}$$

Now we focus on case (b).  Let $T$ be a square with the center at $\mu$ and with side length $l^{-\eta-2\epsilon M}$.
On case (b), note that $\nu$ is away $\nu_1$ at least $l^{-\eta-M\epsilon}$ and away $Z$  at least  $l^{-\eta+\epsilon}$. 
By using $\beta_1$, $\zeta'$, $\bar r_1$ and an open $\bullet$-path next $\mu$ inside $\Gamma_1$ to $r_1$ (see Fig. 10), there exist
 six arm paths from $\mu+[-3, 3]^2$  to the boundary of $T$.
 More precisely, $\beta_1$ is from the  boundary of $T$ passing $\mu$ and  back to the boundary of $T$, 
 $\zeta'$  $\circ$-adjacent to $ \mu$ to
 the boundary of $T$,  $\bar r_1$ from the boundary of $T$ passing  a point $\circ$-adjacent  to $\mu$ and then back to the  boundary of $T$, and finally an open $\bullet$-path,
 $\circ$-adjacent to a point of $\bar r_1$ and the point  also $\circ$-adjacent $\mu$, to the boundary of $T$.
 We simply call there are six arm paths from $\mu$.
Let ${\cal Q}_6(h)$ be the event that there is a fixed  $h\in \bigtriangleup \alpha\beta\gamma$ with six arm paths from $h$ to the square boundary, where the square has a  side length $l^{-\eta-2\epsilon M}$ and  the center at $h$.
On ${\cal D}_2$,  $d(h,Z_2) \geq l^{-\eta- 2M\epsilon}$. 
By  (2.1.5) and translation invariance, if $h$ is fixed, 
$${\bf P}_l({\cal Q}_6(h))\leq  {\bf P}_l({\cal Q}_6(l^{-1}, l^{-\eta-2\epsilon M})) \leq c_1( l^{-1} /l^{-\eta-2M\epsilon})^{(2+d_2)}=c_1l^{-2 -2d_2 + 2\eta+ 2d_2\eta +2 M\epsilon(2+d_2)}.\eqno{}$$
On the other hand, there are at most  $O(l^{2-2\eta})$ many lattice points in $\bigtriangleup \alpha\beta\gamma$.  
$${\bf P}_l( \exists \,\,h \mbox{ such that } {\cal Q}_6(h))\leq c_1  l^{-2 -2d_2 +2\eta+ 2d_2\eta+ 2M\epsilon(2+d_2)+ 2-2\eta}= c_1 l^{-2d_2(1-\eta)+
 2M\epsilon(2+d_2)}. \eqno{(A. 32)}$$
 If we take $\epsilon$ small, then  there exist $c_1$  and $c_2$ for $i=1,2$ independent of $\epsilon$ such that
$${\bf P}_l( \exists \,\,h \mbox{ such that } {\cal Q}_6(h)\mbox{ occurs})\leq c_1l^{-c_2 }.\eqno{(A. 33)}$$
On $\{{\cal H}(D, Z_1, Z_2)\cap {\cal D}^+(D, \tau^+(D), {Z}_2),{\cal D}_1, {\cal D}_2, {\cal D}'^C_3,\mbox{ case (b)}\}$, 
${\cal G}_6$ occurs and  ${\cal Q}_6(h)$ occurs for some $h$ with $d(h, Z_2) \geq 
l^{-\eta- 2M\epsilon}$. Moreover,  ${\cal D}_1$ occurs.
Note that these above events occur in the different domains, so by the estimates in (A.13) and  Lemma A.2 for ${\cal G}_6$, 
if we take $\epsilon$ small, then  there exist $c_1$, $c_2$ and  $c_3$ such that 
\begin{eqnarray*}
&&{\bf P}_l ({\cal H}(D, Z_1, Z_2)\cap {\cal D}^+(D, \tau^+(D), {Z}_2),{\cal D}_1, {\cal D}_2, {\cal D}'^C_3,\mbox{ case (b)})\\
&\leq &  {\bf P}_l({\cal G}_6) {\bf P}_l( \exists \,\,h \mbox{ such that } {\cal Q}_6(h)\mbox{ occurs }){\bf P}_l( {\cal D}_1) \leq c_1 l^{-1+\eta+\epsilon} l^{-c_2 }l^{-1-\eta}\leq c_1 l^{-2-c_3}.\hskip 1cm {(A.34)}
\end{eqnarray*}
Lemma A.3 follows from (A.30) and (A.34).   $\blacksquare$\\

$\beta_1$ will meet $r_1$ at $u_1$ and $\beta_2$ will meet $r_2$ at $u_2$.
We denote by ${\cal D}_4$ the event that $d(u_1, \nu_1) \geq l^{-\eta-\epsilon M}$ and $d(u_2, \nu_2) \geq l^{-\eta-\epsilon M}$
for $\eta$ and $M$ defined (A.13). We showed in (A.31) that
 $${\bf P}_l ({\cal H}(D, Z_1, Z_2)\cap {\cal D}^+(D, \tau^+(D), {Z}_2),{\cal D}_1, {\cal D}_2, {\cal D}_3, {\cal D}_4^C)
\leq c_1 l^{-2-c_2\epsilon}.\eqno {(A.35)}$$
Let $\dot \beta_1$ be an open $\bullet$-path from $\overline{BC}$ to $r_1$  at $u_1'$ and  let $\dot\beta_2$ be  a closed $\circ$-path from
$\overline{BC}$ to $r_2$ at $u_2'$.  We denote by  ${\cal D}_5^C$  the event that 
$d(u_1', Z_2) <  l^{-\eta-\epsilon }$ or $d(u_2', Z_2) < l^{-\eta-\epsilon }$ for some $\dot \beta_1$ or  $\dot \beta_2$.
By the same proofs of (A.16), as we mentioned in the sentences below (A.17),
 $${\bf P}_l ({\cal H}(D, Z_1, Z_2)\cap {\cal D}^+(D, \tau^+(D), {Z}_2),{\cal D}_1, {\cal D}_2, {\cal D}_3, {\cal D}_4, {\cal D}_5^C)\leq c_1 l^{-2-c_2\epsilon}.\eqno {(A.36)}$$
 Now we are ready to show (A.1) in Lemma 2.4.1 for $n=1$.\\

{\bf Proof of Lemma 2.4.1 for $n=1$.}\\

By (A.14), Lemmas A.2 and A.3,  (A.35) and (A.36),  there exist $c_i>0$ for $i=1,2$ such that for $\epsilon$ and $M$ defined in (A.14),
\begin{eqnarray*}
&&{\bf P}_l ({\cal H}(D, Z_1, Z_2)\cap {\cal D}^+(D, \tau^+(D), {Z}_2))\\
&\leq &{\bf P}_l ({\cal H}(D, Z_1, Z_2)\cap {\cal D}^+(D, \tau^+(D), {Z}_2),{\cal D}_1, {\cal D}_2,{\cal D}_3, {\cal D}_4, {\cal D}_5)+c_1 l^{-(2+c_2\epsilon)}.\hskip 3cm {(A.37)}
\end{eqnarray*}
By (A.5), 
\begin{eqnarray*}
&&{\bf P}_l ({\cal H}(D, Z_1, Z_2)\cap {\cal D}^+(D, \tau^+(D), {Z}_2),{\cal D}_1, {\cal D}_2,{\cal D}_3, {\cal D}_4, {\cal D}_5)\\
&=& \sum_{  G_1, G_2} {\bf P}_l ({\cal H}(D, Z_1, Z_2)\cap {\cal D}^+(D, \tau^+(D), {Z}_2), {\cal D}_1, {\cal D}_2,{\cal D}_3, {\cal D}_4, {\cal D}_5, \Gamma_1=G_1,\Gamma_2=G_2),\hskip 0.6cm {(A.38)}
\end{eqnarray*}
where the sum takes over all possible vertex sets $G_1\subset \bigtriangleup \alpha\beta\gamma$ and $G_2\subset \bigtriangleup \alpha\beta\gamma$. 
Let $\sigma^\circ( G_1\cup G_2)$ be the lattice points $\circ$-adjacent to $G_1\cup G_2$  above the $x$-axis together with the lattice points of $G_1$ and $G_2$ (see Fig. 11).
 On $\{{\cal H}(D, Z_1, Z_2)\cap {\cal D}^+(D, \tau^+(D), {Z}_2), {\cal D}_1, {\cal D}_2,{\cal D}_3, {\cal D}_4, {\cal D}_5, \Gamma_1=G_1,\Gamma_2=G_2\}$, there are 
 two paths: an  open $\bullet$-path $\dot \beta_1$ from $S$ $\circ$-adjacent $u_1\in \sigma^\circ( G_1\cup G_2)$  and a closed $\circ$-path $\dot \beta_2$ from $\overline {AC}$ also $\circ$-adjacent to $u_2\in \sigma^\circ( G_1\cup G_2)$, where both paths do not use the vertices of $\sigma^\circ( G_1\cup G_2)$. Furthermore, $u_1$ is $\circ$-adjacent $r_1$,  and  $u_2$ is $\circ$-adjacent to $r_2$  (see Fig. 11).
We also require that $\dot \beta_1$ and $\dot \beta_2$ satisfy the requirements in ${\cal D}_i$ for $i=4,5$ (see Fig. 11). More precisely,  for $\epsilon$ and $M$ defined in (A.14),
 $$d(u_1, \nu_1)\geq l^{-\eta-\epsilon M}, d(u_1, Z_2) \geq l^{-\eta-\epsilon }, d(u_2, \nu_2)\geq l^{-\eta-\epsilon M}, d(u_2, Z_2) \geq l^{-\eta-\epsilon}. \eqno{(A.39)}$$
Let 
${\cal E}({G_1, G_2})$ be  the event of existence of  $\dot \beta_1$ and  $\dot \beta_2$. 
  Thus, for fixed $G_1$ and $G_2$,
\begin{eqnarray*}
&&\{ {\cal H}(D, Z_1, Z_2)\cap {\cal D}^+(D, \tau^+(D), {Z}_2),{\cal D}_1, {\cal D}_2,{\cal D}_3,  {\cal D}_4, {\cal D}_5, \Gamma_1=G_1,\Gamma_2=G_2\}\\
&\subset & \{{\cal E}({G_1, G_2}),  {\cal D}_1, \{\Gamma_1=G_1,\Gamma_2=G_2\}\}.\hskip 8.5cm (A.40)
\end{eqnarray*}
If $\Gamma_1=G_1,\Gamma_2=G_2$ for  fixed $G_1$ and $G_2$, then $r_1=g_1$ and $r_2=g_2$ for fixed  $g_1$ and $g_2$, and $\nu_1=h_1$ and $\nu_2=h_2$ for fixed $h_1$ and $h_2$.

\begin{figure}
\begin{center}
\begin{picture}(80,180)(-20,-10)
\setlength{\unitlength}{0.0125in}%
\begin{tikzpicture}
\thicklines
\begin{scope}[>={Stealth[black]},
              every edge/.style={draw=blue,very thick}]
\path [-] (7.9, -2.5) edge [bend left=40](9.1,-0.5);
 \path [-] (6.5, -0.8) edge [bend left=10](6.5,-2.5);
 \path [-] (7.9, -1.5) edge [bend left=40](10,3.4);
\end{scope}
\begin{scope}[>={Stealth[black]},
              every edge/.style={draw=red,very thick}]
    \path [-] (10.2, 3.15) edge [bend right=20](13.5,-2.5);
    \path [-] (7.7, -2.5) edge [bend right=20](6.55,-0.78);
      \path [-] (6.8, 4.3) edge [bend right=40](9.8,3.5);
\path [-] (9.1, -0.5) edge [bend left=10](8.9,-2.5);
\path [-] (7.1, -1) edge [bend right=20](5.5,3);
\end{scope}
\begin{scope}[>={Stealth[black]},
              every edge/.style={draw=green,very thick}]
    \path [-] (7.8, -2.1) edge [bend right=20](6.4,-0.5);
\path [-] (9.2, -0.5) edge [bend left=10](9.1,-2.5);
\path [-] (7.8, -2.1) edge [bend left=40](9.2,-0.43);
 \path [-] (6.4, -0.5) edge [bend left=10](6.3,-2.5);
\end{scope}

\put(0,-80){\line(1,1){250}}
\put(500,-80){\line(-1,1){250}}
\put(150,-80){\line(1,1){100}}
\put(350,-80){\line(-1,1){100}}
\put(0,-80){\line(1,0){500}}
\put(256,-80){\circle*{3}}
\put(242,-80){\circle*{3}}
\put(250,-80){\circle*{3}}

\put(320,105){$S$}

\put(275,-50){$\bar{r}_1$}
\put(265,-35){${r}_1$}
\put(260,-65){$G_1$}
\put(215,-75){$G_2$}
\put(220,-50){${r}_2$}
\put(210,-60){$\bar{r}_2$}
\put(490,-90){$B$}
\put(250,180){$C$}
\put(350,-55){$\lambda_1$}
\put(390,-20){$\Lambda_1$}
\put(270,120){$\Lambda_2$}
\put(-10,-90){$A$}
\put(250,20){$\gamma$}
\put(255,-95){\footnotesize $Z_3$}
\put(245,-95){\footnotesize $Z_2$}
\put(235,-95){\scriptsize $Z_1$}

\put(150,-90){$\alpha$}
\put(340,-90){$\beta$}
\put(200,20){$\dot{\beta}_2$}
\put(255,35){$\dot{\beta}_1$}
\put(225,-25){$u_2$}
\put(245,-40){\small$u_1$}

\put(370,60){$\tau^+(\overline{BC})$}

\end{tikzpicture}
\end{picture}
\end{center}
\caption{\small{\em $\sigma^\circ ({G}_1\cup {G}_2)$ is  the vertices enclosed by the green curves above the $x$ axis, where the green curve is  $\circ$-adjacent to $G_1 \cup G_2$. $\dot{\beta}_1$ is an open $\bullet$-path from $S$ to $\sigma^\circ ({G}_1\cup {G}_2)$ and $\circ$-adjacent to $r_1$.
  $\dot{\beta_2}$ is a closed $\circ$-path from $\overline{AC}$ to $\sigma^\circ ({G}_1\cup {G}_2)$ and $\circ$-adjacent to $r_2$.
  We move $G_1\cup G_2$ one unit $l^{-1}$ to the right.
 $Z_1$ and $Z_2$ are moved to become $Z_2$ and $Z_3$.}}
\end{figure}

Now we  move vertex sets $G_1$ and $G_2$ one unit $l^{-1}$ to the right horizontally,
  but $\dot{\beta}_1$ and $\dot{\beta}_2$ and $\lambda_1$ and $\lambda_2$ remain the same (see Fig. 11).
Let  $\tau^+(\{\Gamma_1=G_1,\Gamma_2=G_2\})$ be the event with configurations of $\tau^+(\omega)$ for $\omega\in 
\{\Gamma_1=G_1,\Gamma_2=G_2\}$ (see Fig. 11). $r_1$ and $r_2$ are also moved to be $\tau^+(r_1)$ and $\tau^+(r_2)$.
If $r_1=g_1$ and $r_2=g_2$, then $g_1$ and $g_2$ are  moved to be $\tau^+(g_1)$ and $\tau^+(g_2)$.
After moving, each edge in $G_1\cup G_2$ is moved one unit $l^{-1}$ to the right, so by (A.6), $\tau^+(\{\Gamma_1=G_1,\Gamma_2=G_2\})$ only depends
on the configurations of $\sigma^\circ (G_1\cup G_2)$. 
On the other hand,  by the definitions of $ {\cal D}_1$, for  fixed $G_1$ and $G_2$, ${\cal E}(G_1, G_2)$ only depends on the configurations outside of $\sigma^\circ (G_1\cup G_2)$. Thus,  for fixed $G_1$ and $G_2$,
$$\!\!\{{\cal E}({G_1, G_2}),  {\cal D}_1\} \mbox{ and } \tau^{+}(\{\Gamma_1=G_1,\Gamma_2=G_2\})\mbox{  are independent.}\eqno{(A. 41)}$$
By (A.40), (A.41),  and translation invariance, 
\begin{eqnarray*}
&& \sum_{  G_1, G_2} {\bf P}_l ({\cal H}(D, Z_1, Z_2)\cap {\cal D}^+(D, \tau^+(D), {Z}_2), {\cal D}_1, {\cal D}_2,{\cal D}_3, {\cal D}_4, {\cal D}_5,  {\cal D}_5, \Gamma_1=G_1,\Gamma_2=G_2)\\
&\leq& \sum_{  G_1, G_2}{\bf P}_l( {\cal E}({G_1, G_2}),  {\cal D}_1){\bf P}_l( \Gamma_1=G_1,\Gamma_2=G_2)\\
&=& \sum_{  G_1, G_2}{\bf P}_l(  {\cal E}({G_1, G_2}),  {\cal D}_1){\bf P}_l( \tau^{+}(\Gamma_1=G_1,\Gamma_2=G_2))\\
&=& \sum_{  G_1, G_2}{\bf P}_l(  {\cal E}({G_1, G_2}),  {\cal D}_1,\tau^{+}( \Gamma_1=G_1,\Gamma_2=G_2)).\hskip 7cm (A. 42)
\end{eqnarray*}

Let ${\cal C}(\dot \beta_1)$  be the open $\bullet$-cluster containing $\dot\beta_1$ outside $\tau^+(G_1\cup G_2)$.
In addition, let ${\cal C}(\dot \beta_2)$ be the closed $\circ$-cluster containing $\dot\beta_2$ outside $\sigma^\circ(G_1\cup G_2)$.
On $ \{ {\cal E}({G_1, G_2}),   {\cal D}_1, \tau^{+}(\Gamma_1=G_1,\Gamma_2=G_2)\}$, if 
${\cal C}(\dot{\beta}_1)$  is $\bullet$-adjacent to $\tau^{+}(g_1)$ and ${\cal C}(\dot{\beta}_2)$ is $\circ$-adjacent to $\tau^{+}(g_2)$, denoted by ${\cal E}_1(G_1, G_2)$, respectively, then
$${\cal E}_1(G_1, G_2)\cap \tau^{+}( \Gamma_1=G_1,\Gamma_2=G_2)\cap {\cal D}_1\subset {\cal H}(D, Z_2, Z_3)\cap {\cal D}^+(D, \tau^+(D), Z_3).\eqno{(A. 43)}$$ 
 Note that $\{\tau^{+}( \Gamma_1=G_1,\Gamma_2=G_2)\}$ are disjoint in different $G_1$ or $G_2$, so by (A.42) and (A.43)
\begin{eqnarray*}
&& \sum_{  G_1, G_2} {\bf P}_l ({\cal H}(D, Z_1, Z_2)\cap {\cal D}^+(D, \tau^+(D), {Z}_2), {\cal D}_1, {\cal D}_2,{\cal D}_3, {\cal D}_4, {\cal D}_5, \Gamma_1=G_1,\Gamma_2=G_2)\\
&\leq & \sum_{  G_1, G_2}{\bf P}_l(   {\cal D}_1, {\cal E}({G_1, G_2}), ({\cal E}_1(G_1, G_2))^C, \tau^{+}( \Gamma_1=G_1,\Gamma_2=G_2))\\
&&+ \sum_{  G_1, G_2}{\bf P}_l(   {\cal D}_1,{\cal E}_1(G_1, G_2), \tau^{+}( \Gamma_1=G_1,\Gamma_2=G_2))\\
&\leq  & {\bf P}_l(   \bigcup_{G_1, G_2}{\cal D}_1,{\cal E}({G_1, G_2}),  ({\cal E}_1(G_1, G_2))^C, \tau^{+}( \Gamma_1=G_1,\Gamma_2=G_2))\\
&&+  {\bf P}_l(  \bigcup_{G_1, G_2} {\cal D}_1, {\cal E}_1(G_1, G_2), \tau^{+}( \Gamma_1=G_1,\Gamma_2=G_2))\\
&\leq &  {\bf P}_l(   \bigcup_{G_1, G_2}{\cal D}_1,{\cal E}({G_1, G_2}),  ({\cal E}_1(G_1, G_2))^C, \tau^{+}( \Gamma_1=G_1,\Gamma_2=G_2))\\
&&+ {\bf P}_l({\cal H}(D, Z_2, Z_3)\cap {\cal D}^+(D, \tau^+(D), Z_3)).\hskip 7.5cm (A. 44)
\end{eqnarray*}

\begin{figure}
\begin{center}
\begin{picture}(80,180)(-20,-20)
\setlength{\unitlength}{0.0125in}%
\begin{tikzpicture}
\thicklines
\begin{scope}[>={Stealth[black]},
              every edge/.style={draw=blue,very thick}]
\path [-] (7.9, -2.5) edge [bend left=40](9.1,-0.5);
 \path [-] (6.5, -0.8) edge [bend left=10](6.5,-2.5);
 \path [-] (8.4, -0.7) edge [bend left=40](10,3.4);
\end{scope}
\begin{scope}[>={Stealth[black]},
              every edge/.style={draw=red,very thick}]
    \path [-] (10.2, 3.15) edge [bend right=20](13.5,-2.5);
    \path [-] (7.7, -2.5) edge [bend right=20](6.55,-0.78);
      \path [-] (6.8, 4.3) edge [bend right=40](9.8,3.5);
\path [-] (9.1, -0.5) edge [bend left=10](8.9,-2.5);
\path [-] (7.1, -1) edge [bend right=20](5.5,3);
 \path [-] (8.5, -0.7) edge [bend left=40](8.7,-0.1);
 \path [-] (8.5, -0.7) edge [bend right=40](7.3,-1.6);
 \path [-] (8.5, -0.8) edge [bend left=40](9.1,-1.6);

\end{scope}

\put(0,-80){\line(1,1){250}}
\put(500,-80){\line(-1,1){250}}
\put(150,-80){\line(1,1){100}}
\put(350,-80){\line(-1,1){100}}
\put(0,-80){\line(1,0){500}}
\put(255, -35){\framebox(25, 25)[br]{T}}

\put(315,105){$S$}
\

\put(255,-50){\scriptsize$\tau^+({r}_1)$}
\put(490,-90){$\tau^+(B)$}
\put(250,180){$\tau^+(C)$}
\put(380,-55){$\lambda_1$}
\put(390,-20){$\Lambda_1$}

\put(270,120){$\Lambda_2$}

\put(240,100){$\lambda_2$}

\put(-10,-90){$\tau^+(A)$}
\put(250,20){$\gamma$}

\put(250,-95){\scriptsize$Z_3$}
\put(235,-95){\scriptsize$Z_2$}
\put(268,-22){\scriptsize$u_1$}

\put(150,-90){$\alpha$}
\put(340,-90){$\beta$}
\put(190,40){$\dot{\beta}_2$}
\put(275,60){$\dot{\beta}_1$}
\put(235,-30){\scriptsize$\zeta$}
\put(284,-12){\scriptsize$\nu_1$}

\put(370,60){$\tau^+(\overline{BC})$}

\end{tikzpicture}
\end{picture}
\end{center}
\caption{\em \small{The graph shows that  there is no open $\bullet$-path connecting $\dot{\beta}_1$ with $\tau^-(r_1)$.  Thus, there  are six arm paths
from $u_1$ to $\partial T$ for a square $T$   with an edge size $l^{-1-\epsilon}$.  }}
\end{figure}

If $({\cal E}_1(G_1, G_2))^C$ occurs, then  ${\cal C}(\dot{\beta}_1)$   is not $\bullet$-adjacent  to  $\tau^{+}(g_1)$  or ${\cal C}(\dot{\beta}_2)$ is not $\circ$-adjacent  to $\tau^{+}(g_2)$.   We suppose  that  ${\cal C}(\dot{\beta}_1)$   is not $\bullet$-adjacent  to  $\tau^{+}(g_1)$, denoted by event ${\cal G}_9$.
On ${\cal E}({G_1, G_2})\cap {\cal G}_9$,  ${\cal C}(\dot\beta_1)$ is not $\bullet$-adjacent to $r_1$,  by the circuit lemma,  the boundary of ${\cal C}(\dot{\beta}_1)$ contains a closed  $\circ$-path  $\zeta$ from $\overline{\alpha\gamma}\cup \overline{\beta\gamma}$ passing  $u_1$ and back to $\overline{\alpha\gamma}\cup \overline{\beta\gamma}$ (see Fig. 12). 
Note that $u_1$ is $\circ$-adjacent  to $r_1$,  so there are six arm paths from $u_1$: $\zeta$  passing $u_1$, $\dot \beta_1$ $\circ$-adjacent to  $u_1$, $\tau^+(r_1)$ $\circ$-adjacent to $u_1$, and a closed $\circ$-path in $\tau^+(G_1)$ with less than $3l^{-1}$
away from $u_1$ (see Fig. 12).
Let $T$ be the square with the center at $u_1$ and with edge length 
$l^{-1-\epsilon M}$ and let ${\cal Q}_6(u_1)$ be the event  that there are six arm paths from $u_1$ to the boundary of $T$.
We also denote by ${\cal G}_{10}$ the event that there are two arm paths from $Z_2$ and $Z_3$ to the upper or the left or the right boundaries of $Z_3+[-l^{-\eta-\epsilon}/2, l^{-\eta-\epsilon}/2]^2$.  Note that on ${\cal E}({G_1, G_2})$ with (A.39), $d(u_1, \tau^+(\nu_1)) \geq l^{-\eta-\epsilon M}-1$ and $d(u_1, Z_3)\geq l^{-\eta-\epsilon}-1$,   so on $\{ {\cal D}_1,{\cal E}({G_1, G_2}), {\cal G}_9, \tau^{+}( \Gamma_1=G_1,\Gamma_2=G_2)\}$, there are six arm path from $u_1$ to the boundary of $T$ for $u_1\in \bigtriangleup  \alpha\beta\gamma \setminus (Z_3+[-l^{-\eta-\epsilon}, l^{-\eta-\epsilon}]^2)=V$. For each  pair $G_1$ and $G_2$, 
$$ {\cal D}_1\cap{\cal E}({G_1, G_2}) \cap  {\cal G}_9\cap  \tau^{+}( \Gamma_1=G_1,\Gamma_2=G_2)\subset
{\cal G}_{10} \cap \{ \exists u_1\in V \mbox{ such that } {\cal Q}_6(u_1)\mbox{ occurs}\}\cap  {\cal  D}_1.\eqno{(A.45)}$$
Thus, 
\begin{eqnarray*}
&& {\bf P}_l( \bigcup_{G_1, G_2}{\cal D}_1,{\cal E}({G_1, G_2}), {\cal G}_9,  \tau^{+}( \Gamma_1=G_1,\Gamma_2=G_2))\\
& \leq & {\bf P}_l({\cal G}_{10} \cap \{ \exists u_1\in V  \mbox{ such that } {\cal Q}_6(u_1)\mbox{ occurs}\}\cap  {\cal  D}_1).\hskip 7cm (A.46)
\end{eqnarray*}
   On the other hand, 
${\cal G}_{10}$ and $\{ \exists u_1\in V  \mbox{ such that } {\cal Q}_6 (u_1)\mbox{ occurs}\}$ and $  {\cal  D}_1$ are disjoint events, so  by  the estimate in (A.33)
for the six arm paths and the estimate in (A.20) for ${\cal G}_{10}$, if $\epsilon$ is small,  then
there exist $c_i$ for $i=1,2, 4$ independent of $\epsilon$ such that
\begin{eqnarray*}
&& {\bf P}_l( \bigcup_{G_1, G_2}{\cal D}_1,{\cal E}({G_1, G_2}), {\cal G}_9,  \tau^{+}( \Gamma_1=G_1,\Gamma_2=G_2))\\
&\leq & {\bf P}_l({\cal G}_{10} \cap \{ \exists u_1\in V \mbox{ such that } {\cal Q}_6(u_1)\mbox{ occurs}\}\cap  {\cal  D}_1)\\
&\leq &  {\bf P}_l({\cal G}_{10}){\bf P}_l( \exists u_1\in V  \mbox{ such that } {\cal Q}_6(u_1)\mbox{ occurs}){\bf P}_l({\cal  D}_1)\leq c_1  l^{-1+\eta+\epsilon} l^{-c_4} l^{-\eta-1} 
 \leq c_1l^{-2-c_2}.\hskip 0.5cm {(A.47)}
 \end{eqnarray*}
We can obtain the same upper bound of (A.47) if $\dot{\beta}_2$ is not $\circ$-adjacent  to $\tau^+(g_2)$.
Thus, there exist $c_i$ for $i=1,2$ such that
$${\bf P}_l(   \bigcup_{G_1, G_2}{\cal D}_1,{\cal E}({G_1, G_2}), ({\cal E}_1(G_1, G_2))^C, \tau^{+}( \Gamma_1=G_1,\Gamma_2=G_2))\\
\leq c_1l^{-2-c_2}.\eqno{(A.48)}$$
Together with (A.37), (A.38), (A.42) and (A.48),  there are $c_i$ for $i=1,2,$ and $\epsilon >0$ such that
$${\bf P}_l ({\cal H}(D, Z_1, Z_2)\cap {\cal D}^+(D, \tau^+(D), {Z}_2))\leq {\bf P}_l ({\cal H}(D, Z_2, Z_3)\cap {\cal D}^+(D, \tau^+(D), {Z}_3))+ c_1 l^{-2-c_2\epsilon}.\eqno{ (A.49)}$$

On ${\cal H}(D, Z_2, Z_3)\cap {\cal D}^+(D, \tau^+(D), {Z}_3)$,  we can also define the same events ${\cal D}_i$ for $i=1,2,3,4,5$ with the  same $\lambda_1$, $\lambda_2$,
$r_1$, $\bar{r}_1$, $r_2$, and $\bar{r}_2$ on $\bigtriangleup  ABC$.
By the same estimates in Lemmas A.1--A.3,  
but with shifting  the corresponding event $\tau^- (\Gamma_1=G_1,\Gamma_2=G_2)$,  and by the same proof of (A. 45), we have
$${\bf P}_l ({\cal H}(D, Z_2, Z_3)\cap {\cal D}^+(D, \tau^+(D), {Z}_3)\leq {\bf P}_l ({\cal H}(D, Z_1, Z_2)\cap {\cal D}^+(D, \tau^+(D), {Z}_2))+ c_1 l^{-2-c_2\epsilon}.\eqno{ (A.50)}$$
By (A.49) and (A.50),   (A.1) follows for $n=1$. \\

{\bf Proof of Lemma 2.4.1 for any $n$.}
Now we focus on any $n$ with $nl^{-1}$ small. The proof is the same for $n=1$, but we just need to rescale the structures such as 
$\bigtriangleup  \alpha\beta\gamma$,  $\bigtriangleup  \alpha'\beta'\gamma'$, $\bigtriangleup  abc $ and $\bigtriangleup  a'b'c'$ from unit $l^{-1}$   as we did for $n=1$ to $nl^{-1}$.
We still divide $\overline{AB}$ into
$l-1$ many equal sub-intervals with length $l^{-1}$.  We also assume that $Z_2$ is the center of $\overline{AB}$ without loss of generality.
 Let $\bigtriangleup  \alpha\beta\gamma$ be a smaller  equilateral  triangle similar to $\bigtriangleup  ABC$ with the center $Z_2$ of $\overline{\alpha\beta}\subset \overline{AB}$ such that
$$d(\alpha,\beta)= (nl^{-1})^{\eta} \mbox{ for }\eta=1/2.\eqno{(A.51)}$$
In addition,  let $\bigtriangleup  \alpha'\beta'\gamma'$ be an  equilateral  triangle  similar to $\bigtriangleup ABC$ also with the center $Z_2$ of
 $\overline{\alpha'\beta'}\subset \overline{AB}$  (see Fig. 13)
such that
$$d({\alpha',\beta'})=  (nl^{-1})^{\delta}\mbox{ for small }0< \delta < 1/2.\eqno{(A.52)}$$
Similarly, we construct $\bigtriangleup  a'b'c'$ and  $\bigtriangleup  abc$ equilateral  triangles similar to $\bigtriangleup  ABC$ with the centers $Z_2$  of $\overline{a'b'}$ and  of $\overline{ab}$  such that 
$$d(a', b')= \kappa d(A, B) \mbox{ for a small }\kappa \mbox{ and } d(a,b)=2d(\alpha, \beta) \eqno{(A.53)}$$

On ${\cal H}(D, Z_1, Z_2)\cap {\cal D}^+(D, \tau^+(D), Z_2)$,
we construct the innermost two arm paths  $\gamma_1$ and $\gamma_2$,  and  the innermost $r_1$ and $r_2$
in $\bigtriangleup \alpha\beta\gamma$ intersecting $\overline{\alpha \gamma}\cup \overline{\beta\gamma}$ at $\nu_1$ and $\nu_2$
as we did for $n=1$.
We also construct the innermost $\bar r_1$ and $\bar r_2$,  and $\Gamma_1$ and $\Gamma_2$, and the innermost $\beta_1$ and  $\beta_2$ from
$r_1$ and $r_2$ to $S$ and to $\overline{AC}$ as we did for $n=1$. With these new constructions, if we replace $\eta$ in (A.51) by
$(nl^{-1})^\eta$ to define ${\cal D}_i$ for $i=1,2,3,4,5$ as we did for $n=1$, by the same proof, we can show that there exist  $c_1$ and $c_2$ such that
\begin{eqnarray*}
&&{\bf P}_l ({\cal H}(D, Z_1, Z_2)\cap {\cal D}^+(D, \tau^+(D), {Z}_2))\\
&\leq &{\bf P}_l ({\cal H}(D, Z_1, Z_2)\cap {\cal D}^+(D, \tau^+(D), {Z}_2),{\cal D}_1, {\cal D}_2,{\cal D}_3, {\cal D}_4, {\cal D}_5)+c_1 (nl^{-1})^{ c_2}l^{-2}.\hskip 3cm {(A.54)}
\end{eqnarray*}
By (A.54), 
\begin{eqnarray*}
&&{\bf P}_l ({\cal H}(D, Z_1, Z_2)\cap {\cal D}^+(D, \tau^+(D), {Z}_2),{\cal D}_1, {\cal D}_2,{\cal D}_3, {\cal D}_4, {\cal D}_5)\\
&=& \sum_{  G_1, G_2} {\bf P}_l ({\cal H}(D, Z_1, Z_2)\cap {\cal D}^+(D, \tau^+(D), {Z}_2), {\cal D}_1, {\cal D}_2,{\cal D}_3, {\cal D}_4, {\cal D}_5, \Gamma_1=G_1,\Gamma_2=G_2),\hskip 1.6cm {(A.55)}
\end{eqnarray*}
where the sum takes over all possible vertex sets $G_1\subset \bigtriangleup \alpha\beta\gamma$ and $G_2\subset \bigtriangleup \alpha\beta\gamma$. 
We denote by 
$$\sigma^\circ( G_1\cup G_2)=\{Z\in {\bf N}_{l^{-1}}\cap \bigtriangleup ABC: d(Z, G_1\cup G_2) \leq n l^{-1}\}.$$
 On $\{{\cal H}(D, Z_1, Z_2)\cap {\cal D}^+(D, \tau^+(D), {Z}_2), {\cal D}_1, {\cal D}_2,{\cal D}_3, {\cal D}_4, {\cal D}_5, \Gamma_1=G_1,\Gamma_2=G_2\}$, there are 
 two paths: an  open $\bullet$-path $\dot \beta_1$ from $\overline{BC} $ $\circ$-adjacent $u_1\in \sigma^\circ( G_1\cup G_2)$  and a closed $\circ$-path $\dot \beta_2$ from $\overline {AC}$ also $\circ$-adjacent to $u_2\in \sigma^\circ( G_1\cup G_2)$, where both paths do not use the vertices of $\sigma^\circ( G_1\cup G_2)$. Furthermore, $d( r_1,u_1)\leq nl^{-1} +1$  and  $d( r_2,u_2)\leq nl^{-1} +1$.
  We also require that $\dot \beta_1$ and $\dot \beta_2$ satisfy the requirements in ${\cal D}_i$ for $i=4,5$ (see Fig. 11). More precisely,  for  small $1/4> \epsilon>0$,  there are squares $T_1$ and $T_2$ with edge size $ (nl^{-1})^{\eta+\epsilon}$ and with the centers at $u_1$ and $u_2$, respectively such that
 $$ T_1\cap r_2=\emptyset, T_2\cap r_1=\emptyset,  d(u_1, Z_2) \geq (nl^{-1})^{\eta+\epsilon}, d(u_2, \nu_2)\geq (nl^{-1})^{\eta+\epsilon}, d(u_2, Z_2) \geq (nl^{-1})^{\eta+\epsilon}. \eqno{(A.56)}$$
Let 
${\cal E}({G_1, G_2})$ be  the event of existence of  $\dot \beta_1$ and  $\dot \beta_2$. By fixed $G_1$ and $G_2$,  we move $G_1$ and $G_2$ in $nl^{-1}$ units to the right horizontally, but $ \dot \beta_1$ and $\dot \beta_2$, $\lambda_1$ and $\lambda_2$ remain the same.
After moving, $Z_1$ and $Z_2$ are moved to be $Z_5$ and $ Z_6$ and $\{\Gamma_1=G_1, \Gamma_2=G_2\}$ will be $\tau^{n} (\Gamma_1=G_1, \Gamma_2=G_2)$.
 With this move, if ${\cal C}(\dot\beta_1)$ and ${\cal C}(\dot\beta_2)$ are the $\bullet$-cluster and the $\circ$-cluster containing
 $ \dot\beta_1$ and $\dot\beta_2$ outside $\sigma^\circ(G_1\cup G_2)$, respectively,  we can check whether ${\cal C}(\dot\beta_1)$ and ${\cal C}(\dot\beta_2)$ $\bullet$-adjacent to $\tau^{n}(r_1)$ and $\circ$-adjacent to $\tau^n(r_2)$. If they do,  then 
$\{{\cal H}(D, Z_5, Z_6)\cap {\cal D}^+(D, \tau^+(D), {Z}_6)\}$ occurs as we proved (A.43). If they do not, then the probability , similar as we estimated in (A.44),  is less than $c_1(nl^{-1})^{-c_2\epsilon } l^{-2}$. Similarly as (A.48) and (A. 49),
we show that 
$${\bf P}_l ({\cal H}(D, Z_1, Z_2)\cap {\cal D}^+(D, \tau^+(D), {Z}_2)\leq {\bf P}_l ({\cal H}(D, Z_5, Z_6)\cap {\cal D}^+(D, \tau^+(D), {Z}_6))+ c_1(nl^{-1})^{c_2\epsilon} l^{-2}.\eqno{ (A.57)}$$

By the same estimates above,  
but with shifting  the corresponding event $\tau^{-n} (\Gamma_1=G_1,\Gamma_2=G_2)$,  and by the same proof of (A.57), we have
$${\bf P}_l ({\cal H}(D, Z_5, Z_6)\cap {\cal D}^+(D, \tau^+(D), {Z}_6)\leq {\bf P}_l ({\cal H}(D, Z_1, Z_2)\cap {\cal D}^+(D, \tau^+(D), {Z}_2))+ c_1(nl^{-1})^{c_2\epsilon} l^{-2}.\eqno{ (A.58)}$$
By (A.57) and (A.58),   (A.1)  follows for any $n$ with $nl^{-1}$ small. 
$\blacksquare$\\

{\bf Acknowledgement.} The author wish to thank F. Wang, Wang, X. Wu for pointing a few mistakes and many suggestions for the first version of the paper. The author also wish to thank  X. Sun and X. Li for their many comments. 

\begin{center}{\large \bf References} \end{center}
Aizenman, M., Duplantier, B.  and Aharony, A. (1999).
 Path crossing exponents and the external perimeter in 2D percolation. {\em Phy. Rev. Let.} {\bf 83}
1359--1362.\\
Bollobas, B. and Riordan, O. (2006). {\em Percolation}. Cambridge University Press, New York. 
Broadbent, S. R. and  Hammersley, J. M. (1957). Percolation processes I. Crystals and mazes. {\em Proceedings of the Cambridge Philosophical Society} {\bf 53} 629--641.\\
Cardy, J. (1992). Critical percolation in finite geometries. {\em  J. of Phys A.} {\bf 25}  L201.\\
Grimmett, G.  (1999). {\em Percolation.} Springer-Verlag,Â  New York.\\
Higuchi, Y., Takei, M., Zhang, Y. (2012). Scaling relations for two-dimensional Ising percolation. {\em J. Statist. Phys.} {\bf 148} 777--799.\\
Kesten, H. (1982). {\em Percolation theory for mathematicians}. Birkhauser, Boston.\\
Kesten, H. (1987). Scaling relations for 2D-percolation. {\em Comm.
Math. Phys.} {\bf 109} 109--156.\\
Kesten, H., Sidoravicius, V. and Zhang, Y. (1998). Almost all words are seen in critical site percolation on the triangular lattice.  {\em Electron. J. Probab.} {\bf 3}  1-75.\\
Kesten, H. and Zhang, Y. (1987). Strict inequalities for some critical exponents in two-dimensional percolation.
{\em J. Statist.  Phys.}  {\bf  46} 1031--1055.\\
Langlands, R., Pouliot, P. and Saint-Aubin, Y. (1994). Conformal invariance in two-dimensional percolation. {\em Bull.  Am. Math. Soc.} {\bf 30}  1--61.\\
Smirnov, S. (2001). Critical percolation in the plane: Conformal invariance, Cardy's formula, scaling limits. {\em  C. R. Acad.  Sci. I-Math.} {\bf 333}  239--244.\\
Werner, W. (2008). Lectures on two-dimensional critical percolation, arXiv:0710.0856.
(math. PR.)\\

Yu Zhang\\
Department of Mathematics\\
University of Colorado\\
Colorado Springs, CO 80933\\
yzhang3@uccs.edu
\end{document}